\documentclass[a4paper,12pt]{amsart}

\usepackage{fullpage}

\calclayout

\makeindex

\usepackage{amsfonts,graphics,amsmath,amsthm,amsfonts,amscd,amssymb,amsmath,latexsym,euscript, enumerate}
\usepackage{epsfig}
\usepackage{flafter}
\usepackage[all,cmtip,line]{xy}
\usepackage{array}
\usepackage[english]{babel}
\usepackage{overpic}
\usepackage{subcaption}
\usepackage{multirow}
\usepackage{microtype}
\usepackage{wrapfig}

\usepackage{tikz}
\newcommand{\blowup}[1]{\tikz \node[draw,circle, inner sep=0pt, minimum size=.5mm]{#1};}
\usetikzlibrary{positioning}
\usepackage{tikz-cd}
\usepackage{float}
\usepackage{longtable}
\usepackage{pifont}
\allowdisplaybreaks
\usepackage{setspace}

\usepackage{graphicx}
\usepackage[hyphens]{url}
\usepackage[bookmarks=false]{hyperref}
\hypersetup{
    colorlinks=true,    
    linkcolor=blue,          
    citecolor=blue,      
    filecolor=blue,      
    urlcolor=blue        
}

\newtheorem{theorem}{Theorem}[section]
\newtheorem*{maintheorem}{Main Theorem}

\newtheorem{proposition}[theorem]{Proposition}
\newtheorem{conjecture}[theorem]{Conjecture}

\newtheorem*{theorem*}{Theorem}
\theoremstyle{plain}
\newtheorem{corollary}[theorem]{Corollary}

\theoremstyle{definition} 
\newtheorem*{ack}{Acknowledgement}
\newtheorem{definition}[theorem]{Definition}
\newtheorem{definition-lemma}[theorem]{Definition-Lemma}

\newtheorem{example}[theorem]{Example}

\newtheorem*{notation}{Notation}
\numberwithin{equation}{section}

\def\P{\mathbb{P}}

\newcommand{\qsim}{\sim_{\mathbb{Q}}}
\newcommand{\supp}{\operatorname{Supp}}

\DeclareMathOperator{\Amp}{Amp}

\DeclareMathOperator{\Pic}{Pic}
\DeclareMathOperator{\Pol}{Pol}

\def\Pic{\operatorname{Pic}}

\def\NE{\operatorname{NE}}

\def\amp{\operatorname{Amp}}
\def\camp{\operatorname{Amp^{\mathrm{cyl}}}}

\newcommand{\one}{\mbox{\tiny \ding{172}}}
\newcommand{\two}{\mbox{\tiny \ding{173}}}
\newcommand{\three}{\mbox{\tiny \ding{174}}}
\newcommand{\four}{\mbox{\tiny \ding{175}}}
\newcommand{\five}{\mbox{\tiny \ding{176}}}
\newcommand{\six}{\mbox{\tiny \ding{177}}}
\newcommand{\seven}{\mbox{\tiny \ding{178}}}
\newcommand{\eight}{\mbox{\tiny \ding{179}}}
\newcommand{\nine}{\mbox{\tiny \ding{180}}}
\newcommand{\ten}{\mbox{\tiny \ding{181}}}
\newcommand{\eleven}{\mbox{\miniscule \blowup{\textbf{11}}}}
\newcommand{\twelve}{\mbox{\miniscule \blowup{\textbf{12}}}}
\newcommand{\thirteen}{\mbox{\miniscule \blowup{\textbf{13}}}}

\usepackage{mathtools}

\DeclarePairedDelimiterX{\norm}[1]{\lVert}{\rVert}{#1}

\expandafter\def\expandafter\normalsize\expandafter{%
    \normalsize%
    \setlength\abovedisplayskip{4pt}%
    \setlength\belowdisplayskip{8pt}%
    \setlength\abovedisplayshortskip{-8pt}%
    \setlength\belowdisplayshortskip{2pt}%
}

\usepackage{lmodern}

\makeatletter
\ifcase \@ptsize \relax
  \newcommand{\miniscule}{\@setfontsize\miniscule{3}{7}}
    \newcommand{\stiny}{\@setfontsize\miniscule{5}{7}}
\or
  \newcommand{\miniscule}{\@setfontsize\miniscule{3}{7}}
   \newcommand{\stiny}{\@setfontsize\miniscule{5}{7}}
\or
  \newcommand{\miniscule}{\@setfontsize\miniscule{3}{7}}
    \newcommand{\stiny}{\@setfontsize\miniscule{5}{7}}
\fi
\makeatother
\title{Cylinders in Du Val del Pezzo surfaces\\ of degree one with Picard rank two}

\begin{document}

\author[J.~Kim]{Jaehyun Kim}
\author[D.-W.~Lee]{Dae-Won Lee}
\author[M.~Sawahara]{Masatomo Sawahara}
\address[Jaehyun Kim]{Department of Mathematics, Ewha Womans University, 52 Ewhayeodae-gil, Seodaemun-gu, Seoul 03760, Republic of Korea}
\email{kjh6691@ewha.ac.kr}
\address[Dae-Won Lee]{Department of Mathematics, Ewha Womans University, 52 Ewhayeodae-gil, Seodaemun-gu, Seoul 03760, Republic of Korea}
\email{daewonlee@ewha.ac.kr}
\address[Masatomo Sawahara]{Faculty of Education, Hirosaki University, Bunkyocho 1, Hirosaki-shi, Aomori 036-8560, Japan}
\email{sawahara.masatomo@gmail.com}

\subjclass[2010]{}

\begin{abstract}
We prove that for Du Val del Pezzo surfaces of degree one with Picard rank two, the existence of an anticanonical polar cylinder implies the ample polar cylindricity.
\end{abstract}


\maketitle

All varieties considered in this article are assumed to be defined over an algebraically closed field of characteristic zero. 

\section{Introduction}\label{sec:intro}
As we have seen in \cite{Arzhantsev2013, Arzhantsev2012}, the unipotent components of automorphism groups play an important in the study of algebraic varieties. In particular, earlier studies including \cite{Kishimoto2013} show that $\mathbb{G}_a$-actions on affine cones correspond to a geometric structure on the base variety, called a cylinder. A \emph{cylinder} in a normal projective variety is a Zariski open subset that is isomorphic to $\mathbb{A}^1\times Z$ for some affine variety $Z$. Such a cylinder is called \emph{$H$-polar} if the complement of the cylinder is the support of an effective $\mathbb{Q}$-divisor that is $\mathbb{Q}$-linearly equivalent to a given ample divisor $H$ on the variety. Building on this cylindrical perspective, $\mathbb{G}_a$-actions on polarized del Pezzo cones have been extensively investigated in \cite{Cheltsov2016, Cheltsov2017, Kishimoto2009, Kishimoto2014,  Sawahara2024, Sawahara2025}. Within this framework, special attention has been paid to the notion of the \emph{cylindrical ample cone} -- the set of ample divisors that admit polar cylinders -- with emphasis on the anticanonical class, as stated below.

\begin{conjecture}[\hspace{1sp}{\cite{Cheltsov2017}}]\label{conj}
Let $S$ be a log del Pezzo surface. Then the cylindrical ample cone $\camp(S)$ of $S$ contains the anticanonical class if and only if $\camp(S)= \amp(S)$.
\end{conjecture}

Several results concerning the cylindrical ample cones of smooth del Pezzo surfaces provide useful context for the present work. The anticanonical polar cylindricity for such surfaces is determined in \cite{Cheltsov2016, Kishimoto2009, Kishimoto2014}. In particular, every smooth del Pezzo surface of degree at least four admits an anticanonical polar cylinder, whereas no such structure exists in lower degrees. The polar cylindricity with respect to general ample divisors is studied in \cite{Cheltsov2017}, thereby confirming Conjecture \ref{conj} for smooth del Pezzo surfaces. Motivated by these results in the smooth case, it is natural to extend the study to singular del Pezzo surfaces, particularly those with Du Val singularities. The anticanonical polar cylindricity for such surfaces is first determined in \cite{Cheltsov2016s}.

\begin{theorem}[\hspace{1sp}{\cite{Cheltsov2016s}}]\label{thm:CPW}
Let $S_d$ be a Du Val del Pezzo surface of degree $d$. Then $S_d$ admits an anticanonical polar cylinder if and only if it does not fall into one of the following cases:
\begin{enumerate}
    \item $S_1$ allows only singular points of types $\mathsf{A_1}$, $\mathsf{A_2}$, $\mathsf{A_3}$, or $\mathsf{D_4}$ if any;
    \item $S_2$ allows only singular points of type $\mathsf{A_1}$ if any;
    \item $S_3$ allows no singular point.
\end{enumerate} 
\end{theorem}

\noindent
In addition, it is shown in \cite{Sawahara2024, Sawahara2025} that every Du Val del Pezzo surface of degree at least two that admits an anticanonical polar cylinder also allows a polar cylinder with respect to any ample divisor. These results establish the conjecture in those cases.

\begin{theorem}[\hspace{1sp}{\cite{Sawahara2024, Sawahara2025}}]\label{thm:Sawahara}
Let $S_d$ be a Du Val del Pezzo surface of degree $d$ at least $2$. Then $\camp(S_d)$ contains the anticanonical class if and only if $\camp(S_d)= \amp(S_d)$.
\end{theorem}

The purpose of this paper is to examine Conjecture \ref{conj} for Du Val del Pezzo surfaces of degree one. Among the nineteen singularity types that arise in the case of Picard rank two, we focus on the thirteen types for which the surface admits an anticanonical polar cylinder, listed below. For reference, Table \ref{table:cyl} in Section \ref{sec:appendix} provides the full list of singularity types.
\begin{equation*}
\begin{array}{cc}
 \mathsf{A_4+A_2+A_1}, \mathsf{A_4+A_3}, \mathsf{A_5+2A_1}, \mathsf{A_5+A_2}, \mathsf{A_6+A_1},\\
 \mathsf{A_7'}, \mathsf{A_7''}, \mathsf{D_5+2A_1}, \mathsf{D_5+A_2}, \mathsf{D_6+A_1}, \mathsf{D_7}, \mathsf{E_6+A_1}, \mathsf{E_7}.
  \end{array}
\end{equation*}

\noindent
This paper develops an alternative framework for understanding the polarity of cylinders. We consider the polarity cone determined by the boundary divisor of each cylinder, as introduced in \cite{Perepechko2021}, and related aspects were studied in \cite{Kim2021, Perepechko2013}.
 By comparing the ample cone with these polarity cones, the boundary condition of a cylinder can be reformulated as a covering problem among the associated cones. This approach is also expected to be applicable to cases of higher Picard ranks, which gives further significance to the present study. We conclude this section by stating our main result, which provides an answer to Conjecture \ref{conj}.

\begin{maintheorem}\label{thm:main}
Let $S_1$ be a Du Val del Pezzo surface of degree $1$ with Picard rank $2$. Then $\camp(S_1)$ contains the anticanonical class if and only if $\camp(S_1)= \amp(S_1)$.
\end{maintheorem}
 
The main result is established by Theorems \ref{thm:basic} and \ref{thm:supple}.
The rest of the paper is organized as follows. In Section \ref{sec:pre}, we introduce the necessary preliminaries used throughout the paper. Section \ref{sec:polarity} is devoted to the study of cones relevant to the polarity of the cylinders. Section \ref{sec:fibration} extends the preceding discussion via $\mathbb{P}^{1}$-fibrations. Finally, in Section \ref{sec:appendix}, we present tables that summarize the results obtained for each singularity type.

\section{Preliminary}\label{sec:pre}
This section begins with relevant definitions and examples that will be used throughout the paper, and then reviews known results on the cylindrical ample cones of del Pezzo surfaces.

\begin{definition}\label{def:cyl}
Let $X$ be a normal projective variety and let $H$ be an ample divisor on $X$. A Zariski open subset $U$ of $X$ is called $H$-polar cylinder if the following conditions hold
\begin{itemize}
 \item[$\bullet$] $U$ is a cylinder, that is, $U$ is isomorphic to $\mathbb{A}^1 \times Z$ for some affine variety $Z$;

 \item[$\bullet$] there exists an effective $\mathbb{Q}$-divisor $D \qsim H$ such that $U=X\setminus \supp(D)$.
 \end{itemize}
\end{definition}

\noindent
Building on this definition, we present examples of cylinders in $\P^2$ that will be used.

\begin{example}\label{ex:line}
Let $L_1$ and $L_2$ be two distinct lines on $\mathbb{P}^2$. For positive real numbers $a,b$ with $a+b=3$, the divisor $aL_1 + bL_2$ defines an anticanonical polar cylinder isomorphic to $\mathbb{A}^1 \times \mathbb{A}^1_{\ast}$, where $\mathbb{A}^1_{\ast}$ is the once-punctured affine line. 

Similarly, if three distinct lines $L_1, L_2$ and $L_3$ that meet at a single point on $\mathbb{P}^2$, then for positive real numbers $a,b,c$ with $a+b+c=3$, the divisor $aL_1 + bL_2 + cL_3$ defines an anticanonical polar cylinder isomorphic to $\mathbb{A}^1 \times \mathbb{A}^1_{\ast\ast}$, where $\mathbb{A}^1_{\ast\ast}$ is the twice-punctured affine line.
\end{example}

\begin{example}\label{ex:conic}
Let $L$ and $Q$ be a line and a conic on $\mathbb{P}^2$ that intersect tangentially at a point. For positive real numbers $a,b$ with $a+2b=3$, the divisor $aL + bQ$ defines an anticanonical polar cylinder isomorphic to $\mathbb{A}^1 \times \mathbb{A}^1_{\ast}$.
\end{example}

Based on the preceding examples, we briefly outline the method of cylinder construction in \cite{Cheltsov2016s}, which we adopt in this paper. The process begins with the construction of an anticanonical polar cylinder on $\mathbb{P}^2$. One then constructs a sequence of blow-ups corresponding to the morphisms $h$ and $g$ in the diagram \eqref{diagram:CPW}, where all exceptional curves arising in the process lie in the complement of the cylinder. Finally, by contracting all $(-2)$-curves in the minimal resolution $\tilde{S}$, one obtains an anticanonical polar cylinder on the resulting Du Val del Pezzo surface $S$. 
\begin{equation}\label{diagram:CPW}
\begin{tikzcd}[row sep=large, column sep=large,
               arrow style=tikz,
               >={Stealth[length=6pt, width=4pt]}]
& & (\check{S}, D_{\check{S}}) \arrow[ddl, "h"'] \arrow[dr, "g"] & & \\
& & & (\tilde{S}, D_{\tilde{S}}) \arrow[d, "f"] \\
& (\mathbb{P}^2, D_{\mathbb{P}^2}) & & (S, D_S) &
\end{tikzcd}
\end{equation}

For the reader's convenience, the complements of the cylinders in $\check{S}$ and $\tilde{S}$ provided in \cite{Cheltsov2016s} are listed in Table \ref{table:tiger}. Subscripts indicate the order of the blow-ups, and superscripts denote the self-intersection numbers.

\begin{definition}\label{def:cone}
Let $X$ be a normal projective variety and $\amp(X)$ the ample cone of $X$ in $\Pic(X)\otimes \mathbb{R}$. 
\begin{itemize}
 \item[$\bullet$] The \emph{cylindrical ample cone} $\camp(X)$ of $X$ is the subset of $\amp(X)$ consisting of all ample divisors on $X$ that admit polar cylinders. 

 \item[$\bullet$] For a cylinder $U$ in $X$, the \emph{polarity cone} $\Pol(U)$ of $U$ is the cone generated by the classes of all irreducible components of $X \setminus U$ in $\Pic(X)\otimes \mathbb{R}$.
 \end{itemize}
\end{definition}

\noindent
The known results on Conjecture \ref{conj} for Du Val del Pezzo surfaces are summarized as follows.

\begin{proposition}[\hspace{1sp}{\cite{Cheltsov2016s,Sawahara2024,Sawahara2025}}]\label{prop:cyl-Amp}
Let $S_d$ be a Du Val del Pezzo surface of degree d.  
\begin{enumerate}
\item[(1)] For $3 \leq d \leq 9$, we have $\camp(S_d) = \amp(S_d)$.

\item[(2)] Let $\mathsf{A}=\left \{\mathsf{6A_1},\mathsf{5A_1},\mathsf{(4A_1)'},\mathsf{(4A_1)''}, \mathsf{(3A_1)'},\mathsf{(3A_1)''},\mathsf{2A_1},\mathsf{A_1}\right \}$ be the set of singularity types of $S_2$ consisting only of $\mathsf{A_1}$-type. Then
\begin{itemize}
 \item[$\bullet$] $-K_{S_2} \not \in \camp(S_2)$, if $S_2$ allows one of the singularities in $\mathsf{A}$;

 \item[$\bullet$] $\camp(S_2)$ = $\amp(S_2)$, otherwise.
 \end{itemize}

\item[(3)] $-K_{S_1} \not \in \camp(S_1)$ if and only if $S_1$ allows one of the singularities in $\mathsf{B}$, where $\mathsf{B}$ denote the set of singularity types of $S_1$ consisting only of $\mathsf{A_1},\mathsf{A_2},\mathsf{A_3}$ or $\mathsf{D_4}$-types as follows: 
\begin{equation*}
   \mathsf{B}=\left \{
\begin{array}{cc}
\mathsf{2A_3 + 2A_1},\mathsf{2A_3 + A_1}, \mathsf{(2A_3)'}, \mathsf{(2A_3)''}, \mathsf{A_3 + A_2 + 2A_1}, \mathsf{A_3 + A_2 + A_1},\\ 
\mathsf{A_3 + A_2}, \mathsf{A_3 + 4A_1},\mathsf{A_3 + 3A_1}, \mathsf{(A_3 + 2A_1)'}, \mathsf{(A_3 + 2A_1)''}, \mathsf{A_3 + A_1}, \mathsf{A_3},\\ 
\mathsf{4A_2}, \mathsf{3A_2 + A_1}, \mathsf{3A_2}, \mathsf{2A_2 + 2A_1}, \mathsf{2A_2 + A_1},
\mathsf{2A_2}, \mathsf{A_2 + 4A_1}, \mathsf{A_2 + 3A_1},\\
\mathsf{A_2 + 2A_1}, \mathsf{A_2 + A_1}, \mathsf{A_2}, \mathsf{6A_1}, \mathsf{5A_1}, \mathsf{(4A_1)'},\mathsf{(4A_1)''}, \mathsf{3A_1},\mathsf{2A_1}, \mathsf{A_1}\\
\mathsf{2D_4}, \mathsf{D_4 + A_3}, \mathsf{D_4 + A_2}, \mathsf{D_4 + 3A_1}, \mathsf{D_4 + 2A_1}, \mathsf{D_4 + A_1}, \mathsf{D_4}
 \end{array}
 \right \}.
\end{equation*}
\end{enumerate}
\end{proposition}

It is worth noting that the cylinder $U$ is $H$-polar for any ample class $H$ in the relative interior of $\Pol(U)$. Consequently, if the ample cone of a variety is entirely contained in $\Pol(U)$, then the variety admits a polar cylinder with respect to every ample divisor. 

\begin{notation}\label{notation}
We conclude this section by summarizing the notation which will be used repeatedly.
\begin{itemize}
\item[$\bullet$] Let $S$ be a Du Val del Pezzo surface of degree $1$ with Picard rank $2$ and $U_0$ denotes the anticanonical polar cylinder in $S$ constructed in \cite{Cheltsov2016s}.

\item[$\bullet$]  Let $\overline{\NE}(S)$ and $\Amp(S)$ be the Mori cone and the ample cone of the surface $S$, respectively.
 
 \item[$\bullet$] In the diagram \eqref{diagram:CPW}, $\ell$ a general line on $\mathbb{P}^2$, $E_{i}$ the exceptional curves of $h$, $e_{i}$ the total transform of the exceptional curves on $\check{S}$, and $\overline{e}_{i}$ the image of $E_{i}$ on the surface $S$.
 
  \item[$\bullet$] For ADE-types, we use $\mathsf{A_{7}'} $ if there is a $(-1)$-curves intersecting the $(-2)$-curve corresponding to the central vertex, and $\mathsf{A_{7}''}$ if there is no such $(-1)$-curve. 
 \end{itemize}

\end{notation}  

It is well known in \cite{Lazarsfeld2004} that the boundary of the Mori cone of the minimal resolution of $S$ is generated by the $(-1)$-curves and $(-2)$-curves. This implies that the cone $\overline{\NE}(S)$ is generated by the pushforwards of the $(-1)$-curves.

\section{Polarity of Cylinders}\label{sec:polarity}
In this section, we discuss the polarity of cylinders by comparing the polarity cones of the anticanonical polar cylinders in \cite{Cheltsov2016s} with the ample cones. We first determine the region of the ample cone covered by a single such cylinder $U_0$ for each singularity type. Among the thirteen singularity types,
\begin{equation*}\label{eq:singularities}
  \left \{
\begin{array}{cc}
\mathsf{A_4+A_2+A_1}, \mathsf{A_4+A_3}, \mathsf{A_5+2A_1}, \mathsf{A_5+A_2}, \mathsf{A_6+A_1},\\ 
\mathsf{A_7'}, \mathsf{A_7''}, \mathsf{D_5+2A_1}, \mathsf{D_5+A_2}, \mathsf{D_6+A_1}, \mathsf{D_7}, \mathsf{E_6+A_1}, \mathsf{E_7}
 \end{array}
 \right \},
\end{equation*}
we find that the ample cone is entirely covered by the polarity cone of the single cylinder $U_0$ presented in \cite{Cheltsov2016s} for four types: $\mathsf{A_5+A_2}, \mathsf{A'_7}, \mathsf{D_6+A_1}, \mathsf{E_7}$. In addition, for the $\mathsf{A_5+2A_1}$ type, we instead use another cylinder $U_{1}$ to cover the ample cone. The polarity cones of $U_0$ and $U_1$ entirely covers the ample cone. The numerical classes of all lines on the minimal resolutions of Du Val del Pezzo surfaces are listed in Table \ref{table:line}.

\begin{theorem}\label{thm:basic}
Let $S$ be a Du Val del Pezzo surface of degree $1$ with Picard rank $2$. Assume that $S$ allows one of the following singularity types: 
\begin{equation*}
\left \{ \mathsf{A_5+A_2}, \mathsf{A'_7}, \mathsf{D_6+A_1}, \mathsf{E_7}, \mathsf{A_5+2A_1} \right \}.
\end{equation*}
Then $\camp(S)= \amp(S)$.
\end{theorem}

\begin{proof}
Based on the precise number of lines on the minimal resolution of each Du Val del Pezzo surface given in \cite{Sawahara2023}, we compute the Mori cone and the ample cone for each surface under consideration.\\

\noindent
\textbf{Case 1.} $\mathsf{A_5+A_2}$\\
The minimal resolution of the surface $S$ has exactly twelve lines:

\begin{equation*}
\begin{cases}
  \ell_{1}=e_4,\\
  \ell_{2}=e_5,\\
  \ell_{3}=e_{10},\\
  \ell_{4}=\ell-e_1-e_6,\\
\end{cases}
\qquad
\begin{cases}
  \ell_{5}=2\ell-e_1-e_2-e_6-e_7-e_8,\\
  \ell_{6}=3\ell-e_1-e_2-e_3-e_4-2e_6-e_7-e_8,\\
  \ell_{7}=3\ell-e_1-e_2-e_3-e_5-2e_6-e_7-e_8,\\
  \ell_{8}=4\ell-2e_1-e_2-e_3-e_4-2e_6-2e_7-e_8-e_9,\\
  \ell_{9}=4\ell-2e_1-e_2-e_3-e_5-2e_6-2e_7-e_8-e_9,\\
  \ell_{10}=6\ell-2e_1-2e_2-2e_3-e_4-e_5-3e_6-3e_7-2e_8-e_9,\\
  \ell_{11}=6\ell-2e_1-2e_2-2e_3-2e_4-3e_6-3e_7-e_8-e_9-e_{10},\\
  \ell_{12}=6\ell-2e_1-2e_2-2e_3-2e_5-3e_6-3e_7-e_8-e_9-e_{10}.
\end{cases}
\end{equation*}

\noindent
If we let $\Pic(S)=\mathbb{Z} \overline{\ell} \oplus \mathbb{Z}\overline{e}_{5}$, then 
the Mori cone $\overline{\NE}(S)$ of $S$ is generated by 
\begin{equation*}
\begin{cases}
\overline{\ell}_{1} \equiv \overline{\ell}_{7} \equiv \overline{\ell}_{9} 
   \equiv \overline{e}_{4} \equiv \tfrac{1}{3}\overline{\ell}-\overline{e}_{5}, \\ 
\overline{\ell}_{2} \equiv \overline{\ell}_{6} \equiv \overline{\ell}_{8} 
   \equiv \overline{e}_{5}, \\ 
\overline{\ell}_{3} \equiv \overline{\ell}_{4} \equiv \overline{\ell}_{5} 
   \equiv \overline{\ell}_{10} \equiv \overline{e}_{10} 
   \equiv \tfrac{1}{6}\overline{\ell},
\end{cases}
\qquad
\begin{cases}
\overline{\ell}_{11} \equiv -\tfrac{1}{6}\overline{\ell}+2\overline{e}_{5}, \\ 
\overline{\ell}_{12} \equiv \tfrac{1}{2}\overline{\ell}-2\overline{e}_{5}.
\end{cases}
\end{equation*}

\noindent
Furthermore, for the birational morphism $f \circ g\colon \check{S} \to S$ in \eqref{diagram:CPW}, we have
\begin{equation*}
\begin{cases}
(f \circ g)^{\ast}(\overline{\ell})= \ell +2E_1+4E_2+6E_3+14L_1+21L_2+6E_7+3E_6,\\
(f \circ g)^{\ast}(\overline{e}_{5})= e_{5} +\frac{1}{2}E_1+E_2+\frac{3}{2}E_3+\frac{5}{2}L_1+\frac{7}{2}L_2+E_7+\frac{1}{2}E_6,
\end{cases}
\end{equation*} 
with $(\overline{\ell})^2=36, (\overline{e}_{5})^2=\frac{1}{2}$, and $\overline{\ell}\cdot \overline{e}_{5}=6$. This implies that the ample cone $\Amp(S)$ of $S$ is the intersection of the regions
\begin{equation*}
\begin{cases}
6a+\tfrac{3}{2}b>0,\\ 
6a+\tfrac{1}{2}b>0,\\ 
6a+b>0,
\end{cases}
\qquad
\begin{cases}
6a>0,\\ 
2b>0.
\end{cases}
\end{equation*}

\noindent
Meanwhile, the polarity cone $\Pol(U_{0})$ of $U_{0}$ is generated by $\overline{\ell}_1, \overline{\ell}_{2}$ and $\overline{\ell}_{3}$. Hence, the ample cone is covered by the polarity cone as follows. 

\begin{figure}[H]
  \vspace{-1em}
\begin{subfigure}[t]{0.5\textwidth}\centering
\begin{tikzpicture}[scale=0.45, every node/.style={font=\scriptsize}, baseline=(origin)] 
  \coordinate (origin) at (0,0);
  \begin{scope}
    \clip (-5,-5) rectangle (5,5);

    \clip (0,0) -- (5,0) -- (5,5) -- (0,5) -- cycle;

    \clip (0,0) -- (0,5) -- (5,5) -- (5,0) -- cycle;

    \fill[yellow!100] (-5,-5) rectangle (5,5);
  \end{scope}

  \draw[->] (-5.5,0) -- (5.5,0) node[right] {$a$};
  \draw[->] (0,-5.5) -- (0,5.5) node[above] {$b$};

  \begin{scope}
    \clip (-5,-5) rectangle (5,5);
    \draw[dotted,very thick] (0,0) -- (0,5);   
    \draw[dotted,very thick] (0,0) -- (5,0); 
  \end{scope}
\end{tikzpicture}
\end{subfigure}\hfill%
\begin{subfigure}[t]{0.5\textwidth}\centering
\begin{tikzpicture}[scale=0.45, every node/.style={font=\scriptsize}, baseline=(origin)]
  \coordinate (origin) at (0,0);
    \begin{scope}
    \clip (-5,-5) rectangle (5,5);
    \fill[gray!40] (0,0) -- (0,5) -- (5,5) -- (5,-5) -- (1.67,-5) -- cycle;
    
  \end{scope}

  \draw[->] (-5.5,0) -- (5.5,0) node[right] {$a$};
  \draw[->] (0,-5.5) -- (0,5.5) node[above] {$b$};

  \draw[thick] (0,0) -- (0,5);    
  \draw[thick] (0,0) -- (1.67,-5);  

  \node[fill=white, inner sep=1pt, anchor=east] at (3,-5.5) {$b=-3a$};
\end{tikzpicture}
\end{subfigure}
\vspace{-0.5em}
\caption{$\Amp(S)$ and $\Pol(U_{0})$ for type $\mathsf{A_5+A_2}$}
\vspace{-1em}
\end{figure}
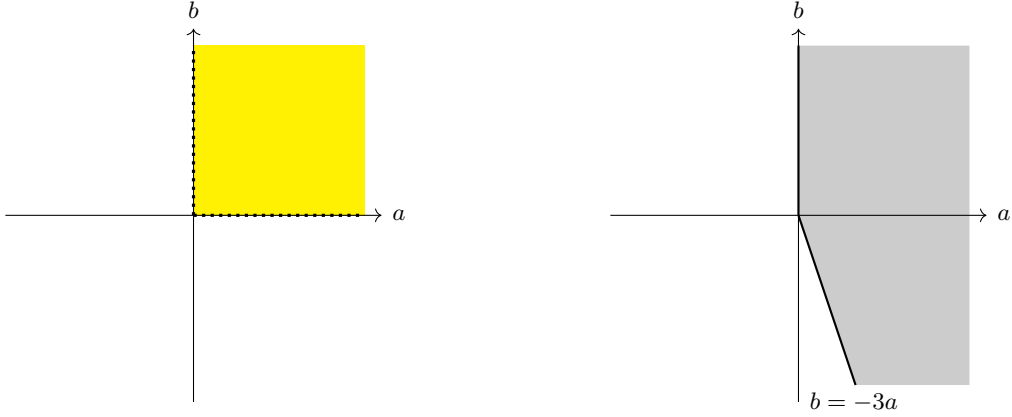

\noindent
\textbf{Case 2.} $\mathsf{A_7'}$\\
The minimal resolution of the surface $S$ has exactly eight lines:

\vspace{-1em}
\begin{equation*}
\begin{cases}
\ell_{1}=e_5,\\
\ell_{2}=e_8,\\
\ell_{3}=\ell-e_1-e_2,\\
\ell_{4}=\ell-e_1-e_6,\\
\end{cases}
\qquad
\begin{cases}
\ell_{5}=\ell-e_6-e_7,\\
\ell_{6}=2\ell-e_1-e_2-e_3-e_4-e_5,\\
\ell_{7}=2\ell-e_1-e_2-e_6-e_7-e_8,\\
\ell_{8}=4\ell-e_1-e_2-e_3-e_4-e_5-2e_6-2e_7-2e_8.
\end{cases}
\end{equation*}

\noindent
If we let $\Pic(S)=\mathbb{Z} \overline{\ell} \oplus \mathbb{Z}\overline{e}_{8}$, then 
the Mori cone $\overline{\NE}(S)$ of $S$ is generated by 
\begin{equation*}
\begin{cases}
\overline{\ell}_{1} \equiv \overline{\ell}_{4} \equiv \overline{e}_{5}
  \equiv \tfrac{1}{2}\overline{\ell}-\tfrac{1}{2}\overline{e}_{8},\\ 
\overline{\ell}_{2} \equiv \overline{\ell}_{3} \equiv \overline{e}_{8},\\ 
\overline{\ell}_{5} \equiv \overline{\ell}_{7} \equiv \overline{\ell}-2\overline{e}_{8},
\end{cases}
\qquad
\begin{cases}
\overline{\ell}_{6} \equiv -\tfrac{1}{2}\overline{\ell}+\tfrac{5}{2}\overline{e}_{8},\\ 
\overline{\ell}_{8} \equiv \tfrac{3}{2}\overline{\ell}-\tfrac{7}{2}\overline{e}_{8}.
\end{cases}
\end{equation*} 

\noindent
Furthermore, for the birational morphism $f\colon \tilde{S} \to S$ in \eqref{diagram:CPW}, we have
\begin{equation*}
\begin{cases}
f^{\ast}(\overline{\ell})= \ell +\frac{3}{4}E_1+\frac{3}{2}E_2+\frac{9}{4}E_3+3E_4+\frac{15}{4}Q+\frac{5}{2}E_7+\frac{5}{4}E_6,\\
f^{\ast}(\overline{e}_{8})= e_{8} +\frac{1}{4}E_1+\frac{1}{2}E_2+\frac{3}{4}E_3+E_4+\frac{5}{4}Q+\frac{3}{2}E_7+\frac{3}{4}E_6,
\end{cases}
\end{equation*} 
with $(\overline{\ell})^2=\frac{17}{2}, (\overline{e}_{8})^2=\frac{1}{2}$, and $\overline{\ell}\cdot \overline{e}_{8}=\frac{5}{2}$. This implies that the ample cone $\Amp(S)$ of $S$ is the intersection of the regions

\begin{equation*}
\begin{cases}
3a+b>0,\\
5a+b>0,\\
7a+3b>0,
\end{cases}
\qquad
\begin{cases}
2a>0,\\
4a+2b>0.
\end{cases}
\end{equation*} 

\noindent
Meanwhile, the polarity cone $\Pol(U_{0})$ of $U_{0}$ is generated by $\overline{\ell}_{1}$, $\overline{\ell}_{2}$ and $\overline{\ell}_{5}$. Hence, the ample cone is covered by the polarity cone as follows. 

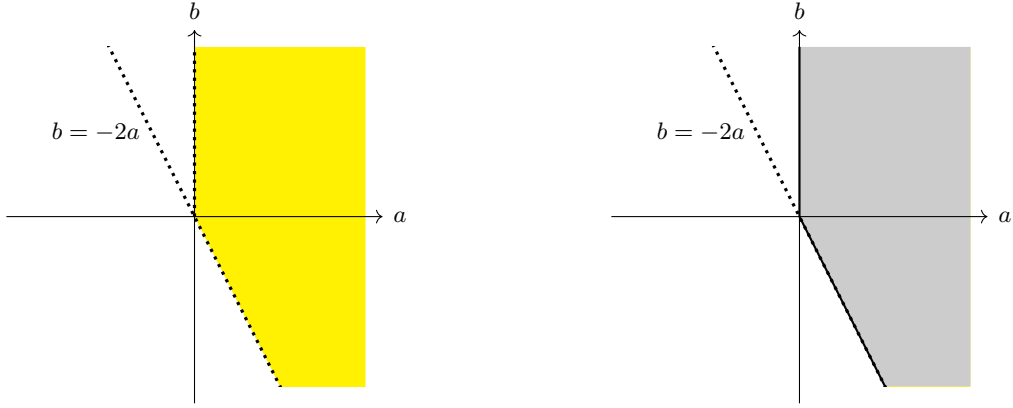
\begin{figure}[H]
  \vspace{-1em}
\begin{subfigure}[t]{0.5\textwidth}\centering
\begin{tikzpicture}[scale=0.45, every node/.style={font=\scriptsize}, baseline=(origin)] 
  \coordinate (origin) at (0,0);
  \begin{scope}
    \clip (-5,-5) rectangle (5,5);

    \clip (-5,10) -- (5,-10) -- (5,5) -- (-5,5) -- cycle;

    \fill[yellow!100] (0,-5) rectangle (5,5);
  \end{scope}

  \draw[->] (-5.5,0) -- (5.5,0) node[right] {$a$};
  \draw[->] (0,-5.5) -- (0,5.5) node[above] {$b$};

  \begin{scope}
    \clip (-5,-5) rectangle (5,5);
    \draw[dotted,very thick] (-5,10) -- (5,-10);    
    \draw[dotted,very thick] (0,0) -- (0,5); 
  \end{scope}

\node[fill=white, inner sep=1pt, anchor=east] at (-1.5,2.5) {$b=-2a$};
\end{tikzpicture}
\end{subfigure}\hfill%
\begin{subfigure}[t]{0.5\textwidth}\centering
\begin{tikzpicture}[scale=0.45, every node/.style={font=\scriptsize}, baseline=(origin)]
  \coordinate (origin) at (0,0);
  \begin{scope}
    \clip (-5,-5) rectangle (5,5);

    \clip (-5,10) -- (5,-10) -- (5,5) -- (-5,5) -- cycle;

    \fill[yellow!90] (0,-5) rectangle (5,5);
  \end{scope}

  \begin{scope}
    \clip (-5,-5) rectangle (5,5);
    \fill[gray!40] (0,0) -- (0,5) -- (5,5) -- (5,-5) -- (2.5,-5) -- cycle;
    
  \end{scope}

  \draw[->] (-5.5,0) -- (5.5,0) node[right] {$a$};
  \draw[->] (0,-5.5) -- (0,5.5) node[above] {$b$};

  \draw[thick] (0,0) -- (0,5);     
  \draw[thick] (0,0) -- (2.5,-5);  

  \begin{scope}
    \clip (-5,-5) rectangle (5,5);
    \draw[dotted,very thick] (-5,10) -- (5,-10); 
  \end{scope}

\node[fill=white, inner sep=1pt, anchor=east] at (-1.5,2.5) {$b=-2a$};
\end{tikzpicture}
\end{subfigure}
\vspace{-0.5em}
\caption{$\Amp(S)$ and $\Pol(U_{0})$ for type $\mathsf{A_7'}$}
\vspace{-1em}
\end{figure}

\noindent
\textbf{Case 3.} $\mathsf{D_6+A_1}$\\
The minimal resolution of the surface $S$ has exactly nine lines:

\vspace{-0.2em}
\begin{align*}
\begin{cases}
\ell_{1}= e_4,\\
\ell_{2}= e_6,\\
\ell_{3}= e_7,\\
\ell_{4}= e_8,\\
\end{cases}
\qquad
\begin{cases}
\ell_{5}= \ell-e_5-e_7,\\
\ell_{6}= \ell-e_5-e_8,\\
\ell_{7}= 3\ell-e_1-e_2-e_3-e_4-2e_5-e_7-e_8,\\
\ell_{8}= 3\ell-e_1-e_2-e_3-e_4-e_5-e_6-2e_7,\\
\ell_{9}= 3\ell-e_1-e_2-e_3-e_4-e_5-e_6-2e_8.
\end{cases}
\end{align*}

\noindent
If we let $\Pic(S)=\mathbb{Z} \overline{\ell} \oplus \mathbb{Z}\overline{e}_{7}$, then 
the Mori cone $\overline{\NE}(S)$ of $S$ is generated by 

\begin{equation*}
\begin{cases}
\overline{\ell}_{1} \equiv \overline{\ell}_{2} \equiv \overline{\ell}_{7} 
  \equiv \overline{e}_{4} \equiv \overline{e}_{6} \equiv \tfrac{1}{3}\overline{\ell},\\
\overline{\ell}_{3} \equiv \overline{\ell}_{6} \equiv \overline{e}_{7},\\
\overline{\ell}_{4} \equiv \overline{\ell}_{5} \equiv \overline{e}_{8} 
  \equiv \tfrac{2}{3}\overline{\ell}-\overline{e}_{7},
\end{cases}
\qquad
\begin{cases}
\overline{\ell}_{8} \equiv \overline{\ell}-2\overline{e}_{7},\\
\overline{\ell}_{9} \equiv -\tfrac{1}{3}\overline{\ell}+2\overline{e}_{7}.
\end{cases}
\end{equation*} 

\noindent
Furthermore, for the birational morphism $f\colon \tilde{S} \to S$ in \eqref{diagram:CPW}, we have
\begin{equation*}
\begin{cases}
f^{\ast}(\overline{\ell})= \ell +3L_1+5E_1+3L_2+4E_2+3E_3+2L_3,\\
f^{\ast}(\overline{e}_{7})= e_{7} +L_1+2E_1+\frac{3}{2}L_2+\frac{3}{2}E_2+E_3+\frac{1}{2}L_3,
\end{cases}
\end{equation*} 
with $(\overline{\ell})^2=9, (\overline{e}_{7})^2=\frac{1}{2}$, and $\overline{\ell}\cdot \overline{e}_{7}=3$. This implies that the ample cone $\Amp(S)$ of $S$ is the intersection of the regions

\begin{equation*}
\begin{cases}
3a+b>0,\\
3a+\tfrac{1}{2}b>0,\\
3a+\tfrac{3}{2}b>0,
\end{cases}
\qquad
\begin{cases}
3a+2b>0,\\
3a>0.
\end{cases}
\end{equation*}

\noindent
Meanwhile, the polarity cone $\Pol(U_{0})$ of $U_{0}$ is generated by $\overline{\ell}_{1}$, $\overline{\ell}_{3}$ and $\overline{\ell}_{4}$. Hence, the ample cone is covered by the polarity cone as follows. 

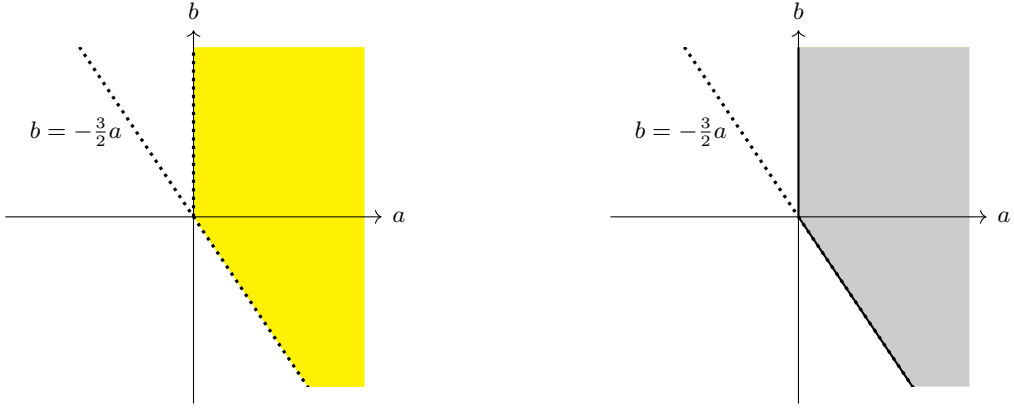
\begin{figure}[H]
  \vspace{-1em}
\begin{subfigure}[t]{0.5\textwidth}\centering
\begin{tikzpicture}[scale=0.45, every node/.style={font=\scriptsize}, baseline=(origin)] 
  \coordinate (origin) at (0,0);
  \begin{scope}
    \clip (-5,-5) rectangle (5,5);

    \clip (-5,7.5) -- (5,-7.5) -- (5,5) -- (-5,5) -- cycle;

    \fill[yellow!100] (0,5) rectangle (5,-7.5);
  \end{scope}

  \draw[->] (-5.5,0) -- (5.5,0) node[right] {$a$};
  \draw[->] (0,-5.5) -- (0,5.5) node[above] {$b$};

  \begin{scope}
    \clip (-5,-5) rectangle (5,5);
    \draw[dotted,very thick] (0,0) -- (0,5);
    \draw[dotted,very thick] (-5,7.5) -- (5,-7.5); 
  \end{scope}

\node[fill=white, inner sep=1pt, anchor=east] at (-2,2.5) {$b=-\frac{3}{2}a$};
\end{tikzpicture}
\end{subfigure}\hfill%
\begin{subfigure}[t]{0.5\textwidth}\centering
\begin{tikzpicture}[scale=0.45, every node/.style={font=\scriptsize}, baseline=(origin)]
  \coordinate (origin) at (0,0);
  \begin{scope}
    \clip (-5,-5) rectangle (5,5); 

    \clip (-5,7.5) -- (5,-7.5) -- (5,5) -- (-5,5) -- cycle;

    \fill[yellow!90] (0,5) rectangle (5,5);
  \end{scope}

  \begin{scope}
    \clip (-5,-5) rectangle (5,5);
    \fill[gray!40] (0,0) -- (0,5) -- (5,5) -- (5,-5) -- (3.33,-5) -- cycle;
    
  \end{scope}

  \draw[->] (-5.5,0) -- (5.5,0) node[right] {$a$};
  \draw[->] (0,-5.5) -- (0,5.5) node[above] {$b$};

  \draw[thick] (0,0) -- (0,5);     
  \draw[thick] (0,0) -- (3.33,-5);

  \begin{scope}
    \clip (-5,-5) rectangle (5,5);
    \draw[dotted,very thick] (-5,7.5) -- (5,-7.5); 
  \end{scope}

\node[fill=white, inner sep=1pt, anchor=east] at (-2,2.5) {$b=-\frac{3}{2}a$};
\end{tikzpicture}
\end{subfigure}
\vspace{-0.5em}
\caption{$\Amp(S)$ and $\Pol(U_{0})$ for type $\mathsf{D_6+A_1}$}
\vspace{-1em}
\end{figure}

\noindent
\textbf{Case 4.} $\mathsf{E_7}$\\
The minimal resolution of the surface $S$ has exactly five lines:
\begin{equation*}
\begin{cases}
\ell_{1}= e_6,\\
\ell_{2}= e_8,\\
\ell_{3}= \ell-e_7-e_8,\\
\end{cases}
\qquad
\begin{cases}
\ell_{4}= 2\ell-e_1-e_2-e_3-e_4-e_5,\\
\ell_{5}= 5\ell-2e_1-2e_2-2e_3-2e_4-2e_5-2e_6-e_7-e_8.
\end{cases}
\end{equation*}

\noindent
If we let $\Pic(S)=\mathbb{Z} \overline{\ell} \oplus \mathbb{Z}\overline{e}_{6}$, then 
the Mori cone $\overline{\NE}(S)$ of $S$ is generated by

\begin{equation*}
\begin{cases}
\overline{\ell}_{1} \equiv \overline{e}_{6},\\
\overline{\ell}_{2} \equiv \overline{e}_{8} \equiv \overline{\ell}-2\overline{e}_{6},\\
\overline{\ell}_{3} \equiv -\overline{\ell}+4\overline{e}_{6},\\
\end{cases}
\qquad
\begin{cases}
\overline{\ell}_{4} \equiv 2\overline{\ell}-5\overline{e}_{6},\\
\overline{\ell}_{5} \equiv 3\overline{\ell}-8\overline{e}_{6}.
\end{cases}
\end{equation*}

\noindent
with $\overline{\ell}_{4} \equiv \overline{Q} \equiv \overline{2\ell}-5\overline{e}_{6}$. Furthermore, for the birational morphism $f\colon \tilde{S} \to S$ in \eqref{diagram:CPW}, we have

\begin{equation*}
\begin{cases}
f^{\ast}(\overline{\ell})= \ell +3E_7+6L+8E_2+4E_1+6E_3+4E_4+2E_5,\\
f^{\ast}(\overline{e}_{6})= e_{6} +E_7+2L+3E_2+\frac{3}{2}E_1+\frac{5}{2}E_3+2E_4+\frac{3}{2}E_5,
\end{cases}
\end{equation*} 

\noindent
with $(\overline{\ell})^2=7, (\overline{e}_{6})^2=\frac{1}{2}$, and $\overline{\ell}\cdot \overline{e}_{6}=2$. This implies that the ample cone $\Amp(S)$ of $S$ is the intersection of the regions

\begin{equation*}
\begin{cases}
4a+b>0,\\
3a+b>0,\\
a>0,
\end{cases}
\qquad
\begin{cases}
4a+\tfrac{3}{2}b>0,\\
5a+2b>0.
\end{cases}
\end{equation*}

\noindent
Meanwhile, the polarity cone $\Pol(U_{0})$ of $U_{0}$ is generated by $\overline{\ell}_{1}$, $\overline{\ell}_{2}$ and $\overline{\ell}_{4}$. Hence, the ample cone is covered by the polarity cone as follows.

\begin{figure}[H]
  \vspace{-1em}
\begin{subfigure}[t]{0.5\textwidth}\centering
\begin{tikzpicture}[scale=0.45, every node/.style={font=\scriptsize}, baseline=(origin)] 
  \coordinate (origin) at (0,0);
  \begin{scope}
    \clip (-5,-5) rectangle (5,5); 

    \clip (-5,12.5) -- (5,-12.5) -- (5,5) -- (-5,5) -- cycle;

    \fill[yellow!100] (0,5) rectangle (5,-12.5);
  \end{scope}

  \draw[->] (-5.5,0) -- (5.5,0) node[right] {$a$};
  \draw[->] (0,-5.5) -- (0,5.5) node[above] {$b$};

  \begin{scope}
    \clip (-5,-5) rectangle (5,5);
    \draw[dotted,very thick] (-5,12.5) -- (5,-12.5);    
    \draw[dotted,very thick] (0,0) -- (0,5); 
  \end{scope}

\node[fill=white, inner sep=1pt, anchor=east] at (-1.5,2.5) {$b=-\frac{5}{2}a$};
\end{tikzpicture}
\end{subfigure}\hfill%
\begin{subfigure}[t]{0.5\textwidth}\centering
\begin{tikzpicture}[scale=0.45, every node/.style={font=\scriptsize}, baseline=(origin)]
  \coordinate (origin) at (0,0);
  \begin{scope}
    \clip (-5,-5) rectangle (5,5); 

    \clip (-5,12.5) -- (5,-12.5) -- (5,5) -- (-5,5) -- cycle;

    \fill[yellow!90] (0,5) rectangle (5,5);
  \end{scope}

  \begin{scope}
    \clip (-5,-5) rectangle (5,5);
    \fill[gray!40] (0,0) -- (0,5) -- (5,5) -- (5,-5) -- (2,-5) -- cycle;
    
  \end{scope}

  \draw[->] (-5.5,0) -- (5.5,0) node[right] {$a$};
  \draw[->] (0,-5.5) -- (0,5.5) node[above] {$b$};

  \draw[thick] (0,0) -- (0,5);    
  \draw[thick] (0,0) -- (2,-5); 

  \begin{scope}
    \clip (-5,-5) rectangle (5,5);
    \draw[dotted,very thick] (-5,12.5) -- (5,-12.5);
  \end{scope}

\node[fill=white, inner sep=1pt, anchor=east] at (-1.5,2.5) {$b=-\frac{5}{2}a$};
\end{tikzpicture}
\end{subfigure}
\vspace{-0.5em}
\caption{ $\Amp(S)$ and $\Pol(U_{0})$ for type $\mathsf{E_7}$}
\vspace{-1em}
\end{figure}
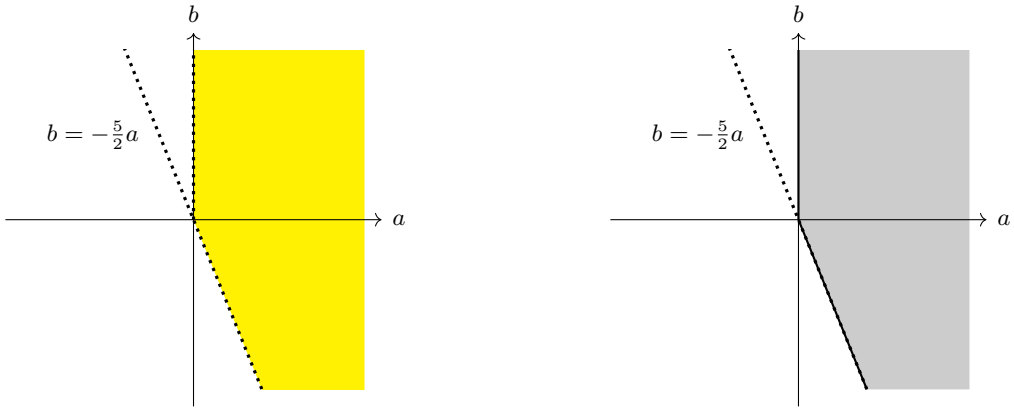

\noindent
\textbf{Case 5.} $\mathsf{A_5+2A_1}$\\
The minimal resolution of the surface $S$ has exactly fourteen lines:
\begin{equation*}
\begin{cases}
\ell_{1}= e_5,\\
\ell_{2}= e_9,\\
\ell_{3}= e_{10},\\
\ell_{4}= \ell-e_1-e_6,
\end{cases}
\qquad
\begin{cases}
\ell_{5}= 2\ell-e_1-e_2-e_6-e_7-e_8,\\
\ell_{6}= 2\ell-e_1-e_2-e_6-e_7-e_{10},\\
\ell_{7}= 3\ell-e_1-e_2-e_3-e_4-2e_6-e_7-e_8.\\
\ell_{8}= 3\ell-e_1-e_2-e_3-e_4-2e_6-e_7-e_{10},\\
\ell_{9}= 4\ell-2e_1-e_2-e_3-e_4-2e_6-2e_7-e_8-e_9,\\
\ell_{10}= 4\ell-2e_1-e_2-e_3-e_4-2e_6-2e_7-e_8-e_{10},\\
\ell_{11}= 6\ell-2e_1-2e_2-2e_3-e_4-e_5-3e_6-3e_7-2e_8-e_9,\\
\ell_{12}= 6\ell-2e_1-2e_2-2e_3-e_4-e_5-3e_6-3e_7-2e_8-e_{10},\\
\ell_{13}= 6\ell-2e_1-2e_2-2e_3-e_4-e_5-3e_6-3e_7-e_8-2e_{10},\\
\ell_{14}= 6\ell-2e_1-2e_2-2e_3-2e_4-3e_6-3e_7-e_8-e_9-e_{10}.
\end{cases}
\end{equation*}

\noindent
If we let $\Pic(S)=\mathbb{Z} \overline{\ell} \oplus \mathbb{Z}\overline{e}_{10}$, then the Mori cone $\overline{\NE}(S)$ of $S$ is generated by 

\begin{equation*}
\begin{cases}
\overline{\ell}_{1} \equiv \overline{\ell}_{4} \equiv \overline{\ell}_{12} \equiv \overline{\ell}_{14} 
   \equiv \overline{e}_{5} \equiv \tfrac{1}{6}\overline{\ell},\\
\overline{\ell}_{2} \equiv \overline{\ell}_{10} \equiv \overline{e}_{9} 
   \equiv \tfrac{1}{4}\overline{\ell}-\tfrac{1}{2}\overline{e}_{10},\\
\overline{\ell}_{3} \equiv \overline{\ell}_{9} \equiv \overline{e}_{10},\\
\overline{\ell}_{5} \equiv \overline{\ell}_{7} 
   \equiv \tfrac{1}{12}\overline{\ell}+\tfrac{1}{2}\overline{e}_{10},
\end{cases}
\qquad
\begin{cases}
\overline{\ell}_{6} \equiv \overline{\ell}_{8} 
   \equiv \tfrac{1}{3}\overline{\ell}-\overline{e}_{10},\\
\overline{\ell}_{11} \equiv -\tfrac{1}{12}\overline{\ell}+\tfrac{3}{2}\overline{e}_{10},\\
\overline{\ell}_{13} \equiv \tfrac{5}{12}\overline{\ell}-\tfrac{3}{2}\overline{e}_{10}.
\end{cases}
\end{equation*}

\noindent
Furthermore, for the birational morphism $f \circ g : \tilde{S} \to S$ in \eqref{diagram:CPW}, we have
\begin{equation*}
\begin{cases}
(f \circ g)^{\ast}(\overline{\ell})= \ell +2E_1+4E_2+6E_3+14L_1+21L_2+6E_7+3E_6,\\
(f \circ g)^{\ast}(\overline{e}_{10})= e_{10} +\frac{1}{3}E_1+\frac{2}{3}E_2+E_3+\frac{7}{3}L_1+\frac{11}{3}L_2+\frac{4}{3}E_7+\frac{2}{3}E_6,
\end{cases}
\end{equation*} 
with $(\overline{\ell})^2=36, (\overline{e}_{10})^2=\frac{1}{3}$, and $\overline{\ell}\cdot \overline{e}_{10}=6$. Hence, the ample cone $\Amp(S)$ of $S$ is the intersection of the regions

\begin{equation*}
\begin{cases}
6a+b>0,\\ 
6a+\tfrac{4}{3}b>0,\\
6a+\tfrac{1}{3}b>0,\\
6a+\tfrac{2}{3}b>0,
\end{cases}
\qquad
\begin{cases}
6a+\tfrac{5}{3}b>0,\\
6a>0,\\ 
6a+2b>0.
\end{cases}
\end{equation*}

\noindent
This implies that the region of ample cone covered by the polarity cone $\Pol(U_{0})$ which is generated by $\overline{\ell}_{1}$, $\overline{\ell}_{2}$ and $\overline{\ell}_{3}$ is as follows.
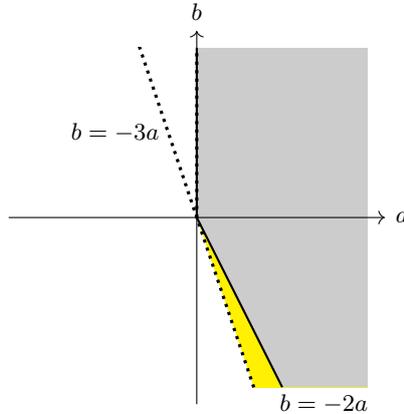
\begin{figure}[H]
  \vspace{-1em}
\begin{subfigure}[t]{0.5\textwidth}\centering
\begin{tikzpicture}[scale=0.45, every node/.style={font=\scriptsize}, baseline=(origin)] 
  \coordinate (origin) at (0,0);
  \begin{scope}
    \clip (-5,-5) rectangle (5,5); 

    \clip (-5,15) -- (5,-15) -- (5,5) -- (-5,5) -- cycle;

    \fill[yellow!100] (0,-5) rectangle (5,5);
  \end{scope}

  \begin{scope}
    \clip (-5,-5) rectangle (5,5);
    \fill[gray!40] (0,0) -- (0,5) -- (5,5) -- (5,-5) -- (2.5,-5) -- cycle;
    
  \end{scope}

  \draw[->] (-5.5,0) -- (5.5,0) node[right] {$a$};
  \draw[->] (0,-5.5) -- (0,5.5) node[above] {$b$};

  \draw[thick] (0,0) -- (0,5);    
  \draw[thick] (0,0) -- (2.5,-5);  

  \begin{scope}
    \clip (-5,-5) rectangle (5,5);
    \draw[dotted,very thick] (0,0) -- (0,5);    
    \draw[dotted,very thick] (-5,15) -- (5,-15); 
  \end{scope}

  \node[fill=white, inner sep=1pt, anchor=east] at (-1,2.5) {$b=-3a$};
  \node[fill=white, inner sep=1pt, anchor=east] at (5.1,-5.5) {$b=-2a$};
\end{tikzpicture}
\end{subfigure}
\vspace{-0.5em}
\caption{The region covered by $\Pol(U_{0})$ for type $\mathsf{A_5+2A_1}$}
\vspace{-1em}
\end{figure}

Now, let $U_{1}$ be a cylinder on $S$ defined by the image of the total transform of $\ell_{6}+L_1$ from $\mathbb{P}^{2}$. Then the ample cone is covered by the polarity cone $\Pol(U_{1})$ which is generated by  $\overline{\ell}_{1}, \overline{\ell}_{2}$, $\overline{\ell}_{3}$ and $\overline{\ell}_{6}$ as follows. Therefore, $S$ contains an ample polar cylinder for any ample divisor on the surface.
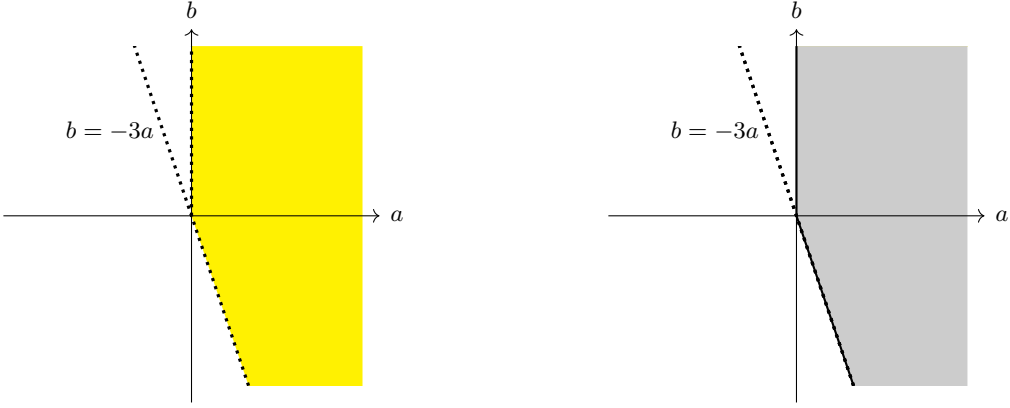
\begin{figure}[H]
  \vspace{-1em}
\begin{subfigure}[t]{0.5\textwidth}\centering
\begin{tikzpicture}[scale=0.45, every node/.style={font=\scriptsize}, baseline=(origin)] 
  \coordinate (origin) at (0,0);
  \begin{scope}
    \clip (-5,-5) rectangle (5,5); 

   \clip (-5,15) -- (5,-15) -- (5,5) -- (-5,5) -- cycle;

    \fill[yellow!100] (0,-5) rectangle (5,5);
  \end{scope}

  \draw[->] (-5.5,0) -- (5.5,0) node[right] {$a$};
  \draw[->] (0,-5.5) -- (0,5.5) node[above] {$b$};

  \begin{scope}
    \clip (-5,-5) rectangle (5,5);
    \draw[dotted,very thick] (0,0) -- (0,5);    
    \draw[dotted,very thick] (-5,15) -- (5,-15); 
  \end{scope}

  \node[fill=white, inner sep=1pt, anchor=east] at (-1,2.5) {$b=-3a$};
\end{tikzpicture}
\end{subfigure}\hfill%
\begin{subfigure}[t]{0.5\textwidth}\centering
\begin{tikzpicture}[scale=0.45, every node/.style={font=\scriptsize}, baseline=(origin)]
  \coordinate (origin) at (0,0);
  \begin{scope}
    \clip (-5,-5) rectangle (5,5); 

    \clip (-5,15) -- (5,-15) -- (5,5) -- (-5,5) -- cycle;

    \fill[yellow!90] (0,5) rectangle (5,5);
  \end{scope}
\begin{scope}
    \clip (-5,-5) rectangle (5,5); 

   \clip (-5,15) -- (5,-15) -- (5,5) -- (-5,5) -- cycle;

    \fill[gray!40] (0,-5) rectangle (5,5);
  \end{scope}

  \draw[->] (-5.5,0) -- (5.5,0) node[right] {$a$};
  \draw[->] (0,-5.5) -- (0,5.5) node[above] {$b$};

   \begin{scope}
    \clip (-5,-5) rectangle (5,5);
    \draw[thick] (0,0) -- (0,5);
    \draw[thick] (0,0) -- (5,-15); 
    \draw[dotted, very thick] (-5,15) -- (5,-15); 
  \end{scope}

  \begin{scope}
    \clip (-5,-5) rectangle (5,5);
    \draw[dotted,very thick] (-5,15) -- (5,-15);
  \end{scope}

\node[fill=white, inner sep=1pt, anchor=east] at (-1,2.5) {$b=-3a$};
\end{tikzpicture}
\end{subfigure}
\vspace{-0.5em}
\caption{ $\Amp(S)$ and $\Pol(U_{1})$ for type $\mathsf{A_{5}+2A_{1}}$}
\vspace{-1em}\hfill\qedhere
\end{figure}
\end{proof}

\section{$\mathbb{P}^1$-Fibrations on minimal resolutions}\label{sec:fibration}
In this section, we study the remaining eight singularity types among the thirteen cases that admit anticanonical polar cylinders. Based on explicit computations of the relative cones as in the preceding section, we first show that the ample cone cannot be covered by the single cylinder constructed in \cite{Cheltsov2016s}. Subsequently, by means of $\mathbb{P}^{1}$-fibrations on the minimal resolutions, we construct additional cylinders so that the  union of the  their polarity cones provides a covering of the ample cone.

\begin{theorem}\label{thm:supple}
Let $S$ be a Du Val del Pezzo surface of degree $1$ with Picard rank $2$. Assume that $S$ allows one of the following singularity types:
\begin{equation*}
\left \{\mathsf{A_6+A_1}, \mathsf{A_7''}, \mathsf{D_5+2A_1}, \mathsf{D_5+A_2}, \mathsf{D_7}, \mathsf{E_6+A_1}, \mathsf{A_4+A_3}, \mathsf{A_4+A_2+A_1}\right \}.
\end{equation*}
Then $\camp(S)= \amp(S)$. 
\end{theorem}

\begin{proof}
We first determine the Mori cones and the ample cones of the Du Val del Pezzo surfaces for each singularity type, and we examine the region of the ample cone covered by the polarity cone of the cylinder $U_0$.

\medskip
\noindent
\textbf{Case 1.} $\mathsf{A_6+A_1}$\\
If we let $\Pic(S)=\mathbb{Z} \overline{\ell} \oplus \mathbb{Z}\overline{e}_{5}$, then the Mori cone $\overline{\NE}(S)$ of $S$ is generated by 

\begin{equation*}
\begin{cases}
\overline{\ell}_{1} \equiv \overline{e}_{1} \equiv \tfrac{4}{3}\overline{\ell}-3\overline{e}_{5},\\
\overline{\ell}_{2} \equiv \overline{e}_{5},\\
\overline{\ell}_{3} \equiv \overline{e}_{8} \equiv \tfrac{1}{3}\overline{\ell},\\
\overline{\ell}_{4} \equiv \overline{\ell}-2\overline{e}_{5},\\
\overline{\ell}_{5} \equiv -\tfrac{1}{3}\overline{\ell}+2\overline{e}_{5},\\
\end{cases}
\qquad
\begin{cases}
\overline{\ell}_{6} \equiv -\tfrac{2}{3}\overline{\ell}+3\overline{e}_{5},\\
\overline{\ell}_{7} \equiv \overline{\ell}_{8} \equiv \tfrac{2}{3}\overline{\ell}-\overline{e}_{5},\\
\overline{\ell}_{9} \equiv \tfrac{5}{3}\overline{\ell}-4\overline{e}_{5},\\
\overline{\ell}_{10} \equiv 2\overline{\ell}-5\overline{e}_{5},
\end{cases}
\end{equation*}

\noindent
and we have 
\begin{equation*}
\begin{cases}
f^{\ast}(\overline{\ell})= \ell +\frac{6}{7}E_2+\frac{12}{7}E_3+\frac{18}{7}E_4+\frac{24}{7}Q+\frac{16}{7}E_7+\frac{8}{7}E_6+\frac{1}{2}L,\\
f^{\ast}(\overline{e}_{5})= e_{5} +\frac{4}{7}E_2+\frac{8}{7}E_3+\frac{12}{7}E_4+\frac{9}{7}Q+\frac{6}{7}E_7+\frac{3}{7}E_6,
\end{cases}
\end{equation*} 
with $(\overline{\ell})^2=\frac{117}{14}, (\overline{e}_{5})^2=\frac{5}{7}$, and $\overline{\ell}\cdot \overline{e}_{5}=\frac{18}{7}$. Hence, the ample cone $\Amp(S)$ of $S$ is the intersection of the regions

\begin{equation*}
\begin{cases}
8a+3b>0,\\
18a+5b>0,\\
13a+4b>0,\\
45a+16b>0,\\
33a+8b>0,
\end{cases}
\qquad
\begin{cases}
5a+b>0,\\
3a+b>0,\\
51a+20b>0,\\
27a+11b>0.
\end{cases}
\end{equation*}

\noindent
This implies the region of ample cone covered by $\Pol(U_{0})$  which is generated by  $\overline{\ell}_{1}$, $\overline{\ell}_{2}$ and $\overline{\ell}_{3}$ is as follows.
\begin{figure}[H]
  \vspace{-2.2em}
\begin{subfigure}[t]{0.5\textwidth}\centering
\begin{tikzpicture}[scale=0.45, every node/.style={font=\scriptsize}, baseline=(origin)] 
  \coordinate (origin) at (0,0);
  \begin{scope}
    \clip (-5,-5) rectangle (5,5);

    \clip (-5,25) -- (5,-25) -- (5,5) -- (-5,5) -- cycle;

    \clip (-5,12.27) -- (5,-12.27) -- (5,5) -- (-5,5) -- cycle;

    \fill[yellow!100] (-5,-5) rectangle (5,5);
  \end{scope}

  \begin{scope}
    \clip (-5,-5) rectangle (5,5);
    \fill[gray!40] (0,0) -- (0,5) -- (5,5) -- (5,-5) -- (2.22,-5) -- cycle;
    
  \end{scope}

  \draw[->] (-5.5,0) -- (5.5,0) node[right] {$a$};
  \draw[->] (0,-5.5) -- (0,5.5) node[above] {$b$};

  \draw[thick] (0,0) -- (0,5);     
  \draw[thick] (0,0) -- (2.22,-5); 

  \begin{scope}
    \clip (-5,-5) rectangle (5,5);
    \draw[dotted,very thick] (-5,25) -- (5,-25);    
    \draw[dotted,very thick] (-5,12.27) -- (5,-12.27); 
  \end{scope}

\node[fill=white, inner sep=1pt, anchor=east] at (-1.3,2.5) {$b=-\frac{27}{11}a$};
\node[fill=white, inner sep=1pt, anchor=east] at (1.8,-5.5) {$b=-5a$};
\node[fill=white, inner sep=1pt, anchor=east] at (5,-5.5) {$b=-\frac{9}{4}a$};
\end{tikzpicture}
\end{subfigure}
\vspace{-0.5em}
\caption{The region covered by $\Pol(U_{0})$ for type $\mathsf{A_{6}+A_{1}}$}
\vspace{-1em}
\end{figure}
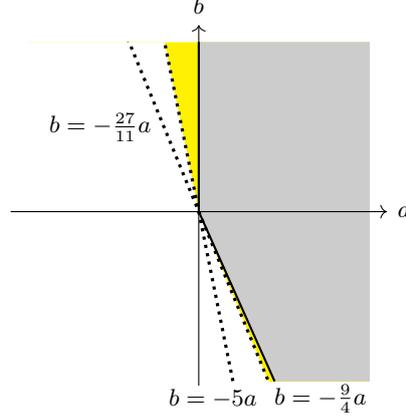

Now for a $\mathbb{P}^{1}$-fibration on the minimal resolution of the surface,

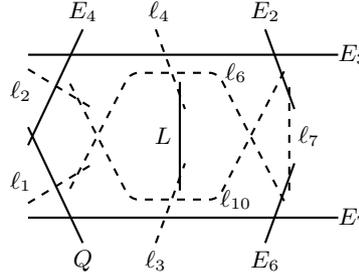
\begin{figure}[H]
  \vspace{-1.3em}
  \centering
\begin{subfigure}[t]{0.5\textwidth}\centering
  \begin{tikzpicture}[scale=0.48, every node/.style={font=\scriptsize}, baseline=(origin)]
\draw [thick] (-0.5,0)--(8,0); 
\draw [thick] (-0.5,4.5)--(8,4.5); 
\draw [thick] (-0.5,2)--(1,5.2); 
\draw [thick] (-0.5,2.5)--(1,-0.7); 
\draw [thick] (3.66,0.75)--(3.66,3.75); 
\draw [thick] (6,-0.7)--(6.8,1.5); 
\draw [thick] (6,5.2)--(6.8,3); 
\draw [thick,dashed] (3,5.2)--(3.8,3); 
\draw [thick,dashed] (3,-0.7)--(3.8,1.5); 
\draw [thick, dashed] (6.66,0.75)--(6.66,3.75); 
\draw [thick,dashed] (-0.5,4.1)--(1.3,3); 
\draw [thick,dashed] (-0.5,0.4)--(1.3,1.5); 
\draw [thick,dashed,rounded corners=3pt] (0.65,3.7)--(2.3,0.5)--(3,0.5)--(4.7,0.5)--(6.6,4.2); 
\draw [thick,dashed,rounded corners=3pt] (0.65,0.8)--(2.3,4)--(3,4)--(4.7,4)--(6.6,0.3); 
\node at (8.4,0) {$E_7$};
\node at (1,5.7) {$E_4$};
\node at (8.4,4.5) {$E_3$};
\node at (1,-1.2) {$Q$};
\node at (6,5.7) {$E_2$};
\node at (5.2,4) {$\ell_6$};
\node at (5.2,0.5) {$\ell_{10}$};
\node at (7.2,2.25) {$\ell_7$};
\node at (6,-1.2) {$E_6$};
\node at (3.2,2.25) {$L$};
\node at (3.1,5.7) {$\ell_4$};
\node at (3.,-1.2) {$\ell_3$};
\node at (-0.7,1) {$\ell_1$};
\node at (-0.7,3.5) {$\ell_2$};
\end{tikzpicture}
\end{subfigure}
\vspace{-1.3em}
\end{figure}

\noindent
we consider a birational morphism from $\tilde{S}$ to $\mathbb{F}_{1}$ obtained by contracting $\ell_{2}, E_{4}, E_{3}, E_{2}, \ell_{3}, L$ and $\ell_{6}$. Then the cylinder defined by the images of $Q, E_{7}, E_{6}$ in $\mathbb{F}_{1}$ implies that there exists a cylinder $U_{1} \cong \mathbb{A}^{1} \times \mathbb{A}^{1}_{\ast}$ in $S$ defined by $\overline{\ell_{2}}, \overline{\ell_{3}}, \overline{\ell_{6}}$. Similarly, a birational morphism obtained by contracting $ \ell_{1}, Q, E_{7}, E_{6}, \ell_{4}, L \ \text{and} \ \ell_{10}$ implies another cylinder $U_{2}\cong \mathbb{A}^{1} \times \mathbb{A}^{1}_{\ast}$ in $S$ defined by $\overline{\ell_{1}}, \overline{\ell_{4}}, \overline{\ell_{10}}$. Therefore, $\Amp(S)$ is covered by the  union of the  the relative interiors of $\Pol(U_{0}), \Pol(U_{1})$ and $\Pol(U_{2})$.\\

\begin{figure}[H]
  \vspace{-2.5em}
\begin{subfigure}[t]{0.5\textwidth}\centering
\begin{tikzpicture}[scale=0.45, every node/.style={font=\scriptsize}, baseline=(origin)] 
  \coordinate (origin) at (0,0);
  \begin{scope}
    \clip (-5,-5) rectangle (5,5);

    \clip (-5,25) -- (5,-25) -- (5,5) -- (-5,5) -- cycle;

    \clip (-5,12.27) -- (5,-12.27) -- (5,5) -- (-5,5) -- cycle;

    \fill[yellow!100] (-5,-5) rectangle (5,5);
  \end{scope}

  \begin{scope}
    \clip (-5,-5) rectangle (5,5);
    \fill[gray!20] (0,0) -- (-1.2,5) -- (5,5) -- (5,0) -- cycle;
   \fill[gray!60] (0,0) -- (1.8,-5) -- (2.5,-5) -- cycle; 
  \end{scope}

  \draw[->] (-5.5,0) -- (5.5,0) node[right] {$a$};
  \draw[->] (0,-5.5) -- (0,5.5) node[above] {$b$};

  \draw[thick] (0,0) -- (-1.2,5);     
  \draw[thick] (0,0) -- (1.8,-5);
  \draw[thick] (0,0) -- (2.5,-5);  
  \begin{scope}
    \clip (-5,-5) rectangle (5,5);
    \draw[dotted,very thick] (-5,25) -- (5,-25);    
    \draw[dotted,very thick] (-5,12.27) -- (5,-12.27); 
  \end{scope}

\node[fill=white, inner sep=1pt, anchor=east] at (-0.7,5.5) {$b=-\frac{9}{2}a$};
\node[fill=white, inner sep=1pt, anchor=east] at (2.3,-5.5) {$b=-\frac{5}{2}a$};
\node[fill=white, inner sep=1pt, anchor=east] at (4.9,-3.5) {$b=-2a$};
\end{tikzpicture}
\end{subfigure}
\vspace{-0.5em}
\caption{$\Pol(U_{1})$ and $\Pol(U_{2})$ for type $\mathsf{A_{6}+A_{1}}$}
\end{figure}

\noindent
\textbf{Case 2.} $\mathsf{A_7''}$\\
If we let $\Pic(S)=\mathbb{Z} \overline{\ell} \oplus \mathbb{Z}\overline{e}_{8}$, then the Mori cone $\overline{\NE}(S)$ of $S$ is generated by 
\begin{equation*}
\begin{cases}
\overline{\ell}_{1} \equiv \overline{e}_{4} \equiv \tfrac{2}{3}\overline{\ell}-\overline{e}_{8},\\
\overline{\ell}_{2} \equiv \overline{e}_{8},\\
\overline{\ell}_{3} \equiv -\tfrac{1}{3}\overline{\ell}+2\overline{e}_{8},\\
\overline{\ell}_{4} \equiv \overline{\ell}-2\overline{e}_{8},
\end{cases}
\qquad
\begin{cases}
\overline{\ell}_{5} \equiv \tfrac{1}{3}\overline{\ell},\\
\overline{\ell}_{6} \equiv \tfrac{4}{3}\overline{\ell}-3\overline{e}_{8},\\
\overline{\ell}_{7} \equiv -\tfrac{2}{3}\overline{\ell}+3\overline{e}_{8},
\end{cases}
\end{equation*}

\noindent
and we have 
\begin{equation*}
\begin{cases}
f^{\ast}(\overline{\ell})= \ell +E_1+2E_2+3E_3+4Q+3E_7+2E_6+E_5,\\
f^{\ast}(\overline{e}_{8})= e_{8} +\frac{3}{8}E_1+\frac{3}{4}E_2+\frac{9}{8}E_3+\frac{3}{2}Q+\frac{15}{8}E_7+\frac{5}{4}E_6+\frac{5}{8}E_5,
\end{cases}
\end{equation*} 

\noindent
with $(\overline{\ell})^2=9, (\overline{e}_{8})^2=\frac{7}{8}$, and $\overline{\ell}\cdot \overline{e}_{8}=3$. Hence, the ample cone $\Amp(S)$ of $S$ is the intersection of the regions 

\begin{equation*}
\begin{cases}
3a+\tfrac{9}{8}b>0,\\
3a+\tfrac{7}{8}b>0,\\
3a+\tfrac{6}{8}b>0,\\
3a+\tfrac{10}{8}b>0,
\end{cases}
\qquad
\begin{cases}
3a+b>0,\\
3a+\tfrac{11}{8}b>0,\\
3a+\tfrac{5}{8}b>0.
\end{cases}
\end{equation*}

\noindent
This implies that the region of ample cone covered by the polarity cone $\Pol(U_{0})$ which is generated by $\overline{\ell}_{1}$, $\overline{\ell}_{2}$ and $\overline{\ell}_{3}$ is as follows.

\begin{figure}[H]
  \vspace{-1 em}
\begin{subfigure}[t]{0.5\textwidth}\centering
\begin{tikzpicture}[scale=0.45, every node/.style={font=\scriptsize}, baseline=(origin)] 
  \coordinate (origin) at (0,0);
  \begin{scope}
    \clip (-5,-5) rectangle (5,5);

    \clip (-5,24) -- (5,-24) -- (5,5) -- (-5,5) -- cycle;

    \clip (-5,10.91) -- (5,-10.91) -- (5,5) -- (-5,5) -- cycle;

    \fill[yellow!100] (-5,-5) rectangle (5,5);
  \end{scope}

  \begin{scope}
    \clip (-5,-5) rectangle (5,5);
    \fill[gray!40] (0,0) -- (-0.83,5) -- (5,5) -- (5,-5) -- (3.33,-5) -- cycle;
  \end{scope}

  \draw[->] (-5.5,0) -- (5.5,0) node[right] {$a$};
  \draw[->] (0,-5.5) -- (0,5.5) node[above] {$b$};

  \draw[thick] (0,0) -- (-0.83,5);    
  \draw[thick] (0,0) -- (3.33,-5);  

  \begin{scope}
    \clip (-5,-5) rectangle (5,5);
    \draw[dotted,very thick] (-5,24) -- (5,-24);   
    \draw[dotted,very thick] (-5,10.91) -- (5,-10.91); 
  \end{scope}

\node[fill=white, inner sep=1pt, anchor=east] at (-1.6,2.5) {$b=-\frac{24}{11}a$};
  \node[fill=white, inner sep=1pt, anchor=east] at (1.8,-5.5) {$b=-\frac{24}{5}a$};
  \node[fill=white, inner sep=1pt, anchor=east] at (6,-5.5) {$b=-\frac{3}{2}a$};
   \node[fill=white, inner sep=1pt, anchor=east] at (2.4,3.5) {$b=-6a$};
\end{tikzpicture}
\end{subfigure}
\vspace{-0.5em}
\caption{The region covered by $\Pol(U_{0})$ for type $\mathsf{A_7''}$}
\vspace{-1em}
\end{figure}

Now for a $\mathbb{P}^{1}$-fibration on the minimal resolution of the surface,
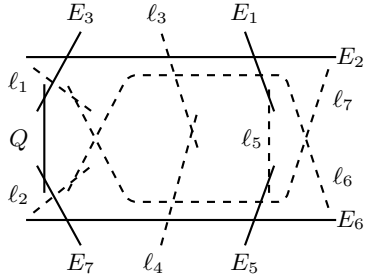
\begin{figure}[H]
  \vspace{-1em}
\begin{subfigure}[t]{0.5\textwidth}\centering
  \begin{tikzpicture}[scale=0.48, every node/.style={font=\scriptsize}, baseline=(origin)]
\draw [thick] (-0.5,0)--(8,0); 
\draw [thick] (-0.5,4.5)--(8,4.5); 
\draw [thick] (-0.2,3)--(1,5.2); 
\draw [thick] (-0.2,1.5)--(1,-0.7); 
\draw [thick] (0,0.75)--(0,3.75); 
\draw [thick] (5.5,-0.7)--(6.3,1.5); 
\draw [thick] (5.5,5.2)--(6.3,3); 
\draw [thick,dashed] (3.2,-0.7)--(4.2,3);
\draw [thick,dashed] (4.2,2)--(3.2,5.2);
\draw [thick, dashed] (6.16,0.75)--(6.16,3.75); 
\draw [thick,dashed] (-0.3,4.2)--(1.3,3); 
\draw [thick,dashed] (-0.3,0.2)--(1.3,1.5); 
\draw [thick,dashed,rounded corners=3pt] (0.65,3.7)--(2.3,0.5)--(3,0.5)--(6.6,0.5)--(7.8,4.2); 
\draw [thick,dashed,rounded corners=3pt] (0.65,0.8)--(2.3,4)--(3,4)--(6.6,4)--(7.8,0.3); 
\node at (8.4,0) {$E_6$};
\node at (1,5.7) {$E_3$};
\node at (8.4,4.5) {$E_2$};
\node at (1,-1.2) {$E_7$};
\node at (5.5,5.7) {$E_1$};
\node at (8.2,3.3) {$\ell_7$};
\node at (8.2,1.2) {$\ell_6$};
\node at (5.7,2.25) {$\ell_5$};
\node at (5.5,-1.2) {$E_5$};
\node at (3.1,5.7) {$\ell_3$};
\node at (3.,-1.2) {$\ell_4$};
\node at (-0.7,0.7) {$\ell_2$};
\node at (-0.7,3.8) {$\ell_1$};
\node at (-0.7,2.2) {$Q$};
\end{tikzpicture}
\end{subfigure}
\vspace{-1em}
\end{figure}

\noindent
we consider a birational morphism from $\tilde{S}$ to $\mathbb{P}^{1} \times \mathbb{P}^{1}$ obtained by contracting $\ell_{1}, E_{3}, \ell_{4}, E_{6}, E_{7}$, $\ell_{6}$ and $E_{1}$. Then the cylinder defined by the images of $Q, E_{2}, E_{5}$ in $\mathbb{P}^{1} \times \mathbb{P}^{1}$ implies  that there exists a cylinder $U_{1} \cong \mathbb{A}^{1} \times \mathbb{A}^{1}_{\ast}$ in $S$ defined by $\overline{\ell_{1}}, \overline{\ell_{4}}, \overline{\ell_{6}}$. Similarly, a cylinder in $\mathbb{F}_{1}$ obtained by contracting
$\ell_{2}, E_{7}, \ell_{3}, E_{2}, E_{1}, \ell_{7} \ \text{and} \ E_{5} \ 
(\text{resp}.\ \ell_{1}, E_{3}, E_{2}, \ell_{2}, \ell_{4}, E_{6} \ \text{and} \ E_{5})
$ implies  that there exists a cylinder $U_{2}$ (resp. $U_{3}$) in $S$ which is isomorphic to $\mathbb{A}^{1} \times \mathbb{A}^{1}_{\ast}$ defined by $\overline{\ell_{2}}, \overline{\ell_{3}}, \overline{\ell_{7}}$ (resp. $\overline{\ell_{1}},\overline{\ell_{2}}, \overline{\ell_{4}}$). Therefore, $\Amp(S)$ is covered by the  union of the  relative interiors of $\Pol(U_{0}), \Pol(U_{1}),\Pol(U_{2})$ and $\Pol(U_{3})$.\\

\begin{figure}[H]
  \vspace{-1.7em}
\begin{subfigure}[t]{0.5\textwidth}\centering
\begin{tikzpicture}[scale=0.45, every node/.style={font=\scriptsize}, baseline=(origin)] 
  \coordinate (origin) at (0,0);
  \begin{scope}
    \clip (-5,-5) rectangle (5,5);

    \clip (-5,24) -- (5,-24) -- (5,5) -- (-5,5) -- cycle;

    \clip (-5,10.91) -- (5,-10.91) -- (5,5) -- (-5,5) -- cycle;

    \fill[yellow!100] (-5,-5) rectangle (5,5);
  \end{scope}

  \begin{scope}
    \clip (-5,-5) rectangle (5,5);
    \fill[gray!20] (0,0) -- (-1.2,5) -- (0,5) -- cycle;
      \fill[gray!60] (0,0) -- (2.05,-5) -- (3.3,-5) -- cycle;
  \end{scope}

  \draw[->] (-5.5,0) -- (5.5,0) node[right] {$a$};
  \draw[->] (0,-5.5) -- (0,5.5) node[above] {$b$};

  \draw[thick] (0,0) -- (-1.2,5);    
  \draw[thick] (0,0) -- (0,5);   
\draw[thick] (0,0) -- (3.3,-5);  
\draw[thick] (0,0) -- (2.05,-5);  
  \begin{scope}
    \clip (-5,-5) rectangle (5,5);
    \draw[dotted,very thick] (-5,24) -- (5,-24);   
    \draw[dotted,very thick] (-5,10.91) -- (5,-10.91); 
  \end{scope}

\node[fill=white, inner sep=1pt, anchor=east] at (-0.7,5.5) {$b=-\frac{9}{2}a$};
  \node[fill=white, inner sep=1pt, anchor=east] at (2.5,-5.5) {$b=-\frac{9}{4}a$};
  \node[fill=white, inner sep=1pt, anchor=east] at (6,-5.5) {$b=-\frac{3}{2}a$};
\end{tikzpicture}
\end{subfigure}\hfill
\begin{subfigure}[t]{0.5\textwidth}\centering
\begin{tikzpicture}[scale=0.45, every node/.style={font=\scriptsize}, baseline=(origin)] 
  \coordinate (origin) at (0,0);
  \begin{scope}
    \clip (-5,-5) rectangle (5,5);

    \clip (-5,24) -- (5,-24) -- (5,5) -- (-5,5) -- cycle;

    \clip (-5,10.91) -- (5,-10.91) -- (5,5) -- (-5,5) -- cycle;

    \fill[yellow!100] (-5,-5) rectangle (5,5);
  \end{scope}

  \begin{scope}
    \clip (-5,-5) rectangle (5,5);
     \fill[gray!80] (0,0) -- (0,5) -- (5,5) -- (5,-5) -- (2.5,-5) -- cycle;
  \end{scope}

  \draw[->] (-5.5,0) -- (5.5,0) node[right] {$a$};
  \draw[->] (0,-5.5) -- (0,5.5) node[above] {$b$};

  \draw[thick] (0,0) -- (0,5);  
    \draw[thick] (0,0) -- (2.5,-5);  
  \begin{scope}
    \clip (-5,-5) rectangle (5,5);
    \draw[dotted,very thick] (-5,24) -- (5,-24);   
    \draw[dotted,very thick] (-5,10.91) -- (5,-10.91); 
  \end{scope}

  \node[fill=white, inner sep=1pt, anchor=east] at (5.1,-5.5) {$b=-2a$};
\end{tikzpicture}
\end{subfigure}
\vspace{-0.5em}
\caption{$\Pol(U_{1}),\Pol(U_{2})$, and $\Pol(U_{3})$ for type $\mathsf{A_7''}$}
\vspace{-1em}
\end{figure}

\medskip

\noindent
\textbf{Case 3.} $\mathsf{D_5+2A_1}$\\
If we let $\Pic(S)=\mathbb{Z} \overline{\ell} \oplus \mathbb{Z}\overline{e}_{8}$, then the Mori cone $\overline{\NE}(S)$ of $S$ is generated by 

\begin{equation*}
\begin{cases}
\overline{\ell}_{1} \equiv \overline{\ell}_{8} \equiv \overline{\ell}_{9} \equiv 
\overline{\ell}_{11} \equiv \overline{e}_{4} \equiv \tfrac{1}{3}\overline{\ell},\\
\overline{\ell}_{2} \equiv \overline{\ell}_{7} \equiv \overline{e}_{6} 
\equiv \tfrac{2}{3}\overline{\ell}-\overline{e}_{8},\\
\overline{\ell}_{3} \equiv \overline{\ell}_{6} \equiv \overline{e}_{8},
\end{cases}
\qquad
\begin{cases}
\overline{\ell}_{4} \equiv \overline{\ell}_{10} \equiv 
-\tfrac{1}{3}\overline{\ell}+2\overline{e}_{8},\\
\overline{\ell}_{5} \equiv \overline{\ell}_{12} \equiv 
\overline{\ell}-2\overline{e}_{8},
\end{cases}
\end{equation*}

\noindent
and we have 
\begin{equation*}
\begin{cases}
f^{\ast}(\overline{\ell})= \ell +2E_1+4E_2+3Q+3E_3+2L,\\
f^{\ast}(\overline{e}_{8})= e_{8} +\frac{3}{4}E_1+\frac{3}{2}E_2+\frac{5}{4}Q+E_3+\frac{1}{2}L+\frac{1}{2}E_7,
\end{cases}
\end{equation*}

\noindent
with $(\overline{\ell})^2=9, (\overline{e}_{8})^2=\frac{3}{4}$, and $\overline{\ell}\cdot \overline{e}_{8}=3$. Hence, the ample cone $\Amp(S)$ of $S$ is the intersection of the regions

\begin{equation*}
\begin{cases}
3a+b>0,\\ 
3a+\tfrac{5}{4}b>0,\\
3a+\tfrac{3}{4}b>0,
\end{cases}
\qquad
\begin{cases}
3a+\tfrac{1}{2}b>0,\\
3a+\tfrac{3}{2}b>0.
\end{cases}
\end{equation*}

\noindent
This implies that the region of ample cone covered by the polarity cone $\Pol(U_{0})$ which is generated by $\overline{\ell}_{1}$, $\overline{\ell}_{2}$ and $\overline{\ell}_{3}$ is as follows.
\begin{figure}[H]
  \vspace{-1.5em}
\begin{subfigure}[t]{0.5\textwidth}\centering
\begin{tikzpicture}[scale=0.45, every node/.style={font=\scriptsize}, baseline=(origin)] 
  \coordinate (origin) at (0,0);
  \begin{scope}
    \clip (-5,-5) rectangle (5,5); 

    \clip (-5,30) -- (5,-30) -- (5,5) -- (-5,5) -- cycle;

    \clip (-5,10) -- (5,-10) -- (5,5) -- (-5,5) -- cycle;

    \fill[yellow!100] (-5,-5) rectangle (5,5);
  \end{scope}

  \begin{scope}
    \clip (-5,-5) rectangle (5,5);
    \fill[gray!40] (0,0) -- (0,5) -- (5,5) -- (5,-5) -- (3.33,-5) -- cycle;
    
  \end{scope}

  \draw[->] (-5.5,0) -- (5.5,0) node[right] {$a$};
  \draw[->] (0,-5.5) -- (0,5.5) node[above] {$b$};

  \draw[thick] (0,0) -- (0,5);     
  \draw[thick] (0,0) -- (3.33,-5);  

  \begin{scope}
    \clip (-5,-5) rectangle (5,5);
    \draw[dotted,very thick] (-5,30) -- (5,-30);    
    \draw[dotted,very thick] (-5,10) -- (5,-10); 
  \end{scope}

\node[fill=white, inner sep=1pt, anchor=east] at (-1.5,2.5) {$b=-2a$};
  \node[fill=white, inner sep=1pt, anchor=east] at (1.8,-5.5) {$b=-6a$};
  \node[fill=white, inner sep=1pt, anchor=east] at (6,-5.5) {$b=-\frac{3}{2}a$};
\end{tikzpicture}
\end{subfigure}
\vspace{-0.5em}
\caption{The region covered by $\Pol(U_{0})$ for type $\mathsf{D_5+2A_1}$}
\vspace{-1em}
\end{figure}

Now for a $\mathbb{P}^{1}$-fibration on the minimal resolution of the surface,

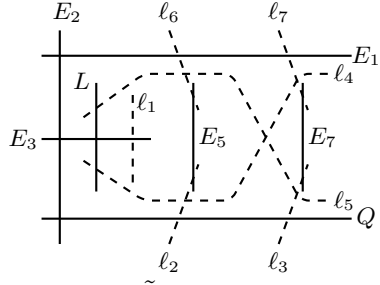
\begin{figure}[H]
  \vspace{-1em}
  \centering
\begin{subfigure}[t]{0.5\textwidth}\centering
  \begin{tikzpicture}[scale=0.48, every node/.style={font=\scriptsize}, baseline=(origin)]
\draw [thick] (-0.5,0)--(8,0); 
\draw [thick] (-0.5,4.5)--(8,4.5); 
\draw [thick] (0,-0.7)--(0,5.2); 
\draw [thick] (1,0.75)--(1,3.75); 
\draw [thick,dashed] (2,1)--(2,3.45); 
\draw [thick] (-0.5,2.2)--(2.5,2.2); 
\draw [thick] (3.66,0.75)--(3.66,3.75); 
\draw [thick,dashed] (6,-0.7)--(6.8,1.5); 
\draw [thick,dashed] (6,5.2)--(6.8,3); 
\draw [thick,dashed] (3,5.2)--(3.8,3); 
\draw [thick,dashed] (3,-0.7)--(3.8,1.5); 
\draw [thick] (6.66,0.75)--(6.66,3.75); 
\draw [thick,dashed,rounded corners=3pt] (0.65,1.6)--(2.3,0.5)--(3,0.5)--(4.7,0.5)--(6.6,4)--(7.5,4); 
\draw [thick,dashed,rounded corners=3pt] (0.65,2.8)--(2.3,4)--(3,4)--(4.7,4)--(6.6,0.5)--(7.5,0.5); 
\node at (8.4,0) {$Q$};
\node at (0.2,5.7) {$E_2$};
\node at (8.4,4.5) {$E_1$};
\node at (6,5.7) {$\ell_7$};
\node at (7.8,4) {$\ell_4$};
\node at (7.2,2.25) {$E_7$};
\node at (6,-1.2) {$\ell_3$};
\node at (4.2,2.25) {$E_5$};
\node at (7.8,0.5) {$\ell_5$};
\node at (3,-1.2) {$\ell_2$};
\node at (3,5.7) {$\ell_6$};
\node at (-1,2.2) {$E_3$};
\node at (0.6,3.8) {$L$};
\node at (2.4,3.2) {$\ell_1$};
\end{tikzpicture}
\vspace{-1.2em}
\end{subfigure} 
\end{figure}

\noindent
we consider a birational morphism from $\tilde{S}$ to $\mathbb{F}_{2}$ obtained by contracting $\ell_{1}, \ell_{3}, E_{7}, \ell_{4}, L, \ell_{6}$ and $E_{5}$. Then the cylinder defined by the images of $Q, E_{1}, E_{2}, E_{3}$ in $\mathbb{F}_{2}$ implies  that there exists a cylinder $U_{1} \cong \mathbb{A}^{1} \times \mathbb{A}^{1}_{\ast \ast}$ in $S$ defined by $\overline{\ell_{1}}, \overline{\ell_{3}}, \overline{\ell_{4}}, \overline{\ell_{6}}$. Similarly, a cylinder in $\mathbb{F}_{2}$ obtained by contracting
$\ell_{1}, \ell_{2}, E_{5}, \ell_{5}, L, \ell_{7} \ \text{and} \ E_{7}$ implies  that there exists a cylinder $U_{2}$ in $S$ which is isomorphic to $\mathbb{A}^{1} \times \mathbb{A}^{1}_{\ast \ast}$ defined by $\overline{\ell_{1}}, \overline{\ell_{2}}, \overline{\ell_{5}}, \overline{\ell_{7}}$. Therefore, $\Amp(S)$ is covered by the  union of the  relative interiors of $\Pol(U_{0}), \Pol(U_{1})$ and $\Pol(U_{2})$.

\begin{figure}[H]
  \vspace{-1em}
\begin{subfigure}[t]{0.5\textwidth}\centering
\begin{tikzpicture}[scale=0.45, every node/.style={font=\scriptsize}, baseline=(origin)] 
  \coordinate (origin) at (0,0);
  \begin{scope}
    \clip (-5,-5) rectangle (5,5); 

    \clip (-5,30) -- (5,-30) -- (5,5) -- (-5,5) -- cycle;

    \clip (-5,10) -- (5,-10) -- (5,5) -- (-5,5) -- cycle;

    \fill[yellow!100] (-5,-5) rectangle (5,5);
  \end{scope}

  \begin{scope}
    \clip (-5,-5) rectangle (5,5);
    \fill[gray!20] (0,0) -- (-0.83,5) -- (5,5) -- (5,-5) -- (3.3,-5) -- cycle;
    
  \end{scope}

  \draw[->] (-5.5,0) -- (5.5,0) node[right] {$a$};
  \draw[->] (0,-5.5) -- (0,5.5) node[above] {$b$};

  \draw[thick] (0,0) -- (-0.83,5);     
  \draw[thick] (0,0) -- (3.33,-5);  

  \begin{scope}
    \clip (-5,-5) rectangle (5,5);
    \draw[dotted,very thick] (-5,30) -- (5,-30);    
    \draw[dotted,very thick] (-5,10) -- (5,-10); 
  \end{scope}

  \node[fill=white, inner sep=1pt, anchor=east] at (-0.5,5.5) {$b=-6a$};
  \node[fill=white, inner sep=1pt, anchor=east] at (5.9,-5.5) {$b=-\frac{3}{2}a$};
\end{tikzpicture}
\end{subfigure}\hfill
\begin{subfigure}[t]{0.5\textwidth}\centering
\begin{tikzpicture}[scale=0.45, every node/.style={font=\scriptsize}, baseline=(origin)] 
  \coordinate (origin) at (0,0);
  \begin{scope}
    \clip (-5,-5) rectangle (5,5); 

    \clip (-5,30) -- (5,-30) -- (5,5) -- (-5,5) -- cycle;

    \clip (-5,10) -- (5,-10) -- (5,5) -- (-5,5) -- cycle;

    \fill[yellow!100] (-5,-5) rectangle (5,5);
  \end{scope}

  \begin{scope}
    \clip (-5,-5) rectangle (5,5);
    \fill[gray!60] (0,0) -- (5,0) -- (5,-5) -- (2.5,-5) -- cycle;
    
  \end{scope}

  \draw[->] (-5.5,0) -- (5.5,0) node[right] {$a$};
  \draw[->] (0,-5.5) -- (0,5.5) node[above] {$b$};

  \draw[thick] (0,0) -- (5,0);     
  \draw[thick] (0,0) -- (2.5,-5);  

  \begin{scope}
    \clip (-5,-5) rectangle (5,5);
    \draw[dotted,very thick] (-5,30) -- (5,-30);    
    \draw[dotted,very thick] (-5,10) -- (5,-10); 
  \end{scope}

  \node[fill=white, inner sep=1pt, anchor=east] at (4.8,-5.5) {$b=-2a$};
\end{tikzpicture}
\end{subfigure}
\vspace{-0.5em}
\caption{$\Pol(U_{1})$ and $\Pol(U_{2})$ for type $\mathsf{D_5+2A_1}$}
\vspace{-1em}
\end{figure}

\medskip

\noindent
\textbf{Case 4.} $\mathsf{D_5+A_2}$\\
If we let $\Pic(S)=\mathbb{Z} \overline{\ell} \oplus \mathbb{Z}\overline{e}_{8}$, then the Mori cone $\overline{\NE}(S)$ of $S$ is generated by 

\begin{equation*}
\begin{cases}
\overline{\ell}_{1} \equiv \overline{\ell}_{8} \equiv \overline{e}_{4} 
   \equiv \tfrac{1}{3}\overline{\ell},\\
\overline{\ell}_{2} \equiv \overline{e}_{5} 
   \equiv \tfrac{4}{3}\overline{\ell}-3\overline{e}_{8},\\
\overline{\ell}_{3} \equiv \overline{e}_{8},\\
\overline{\ell}_{4} \equiv -\tfrac{2}{3}\overline{\ell}+3\overline{e}_{8},\\
\overline{\ell}_{5} \equiv \tfrac{2}{3}\overline{\ell}-\overline{e}_{8},
\end{cases}
\qquad
\begin{cases}
\overline{\ell}_{6} \equiv -\tfrac{1}{3}\overline{\ell}+2\overline{e}_{8},\\
\overline{\ell}_{7} \equiv \overline{\ell}-2\overline{e}_{8},\\
\overline{\ell}_{9} \equiv \tfrac{5}{3}\overline{\ell}-4\overline{e}_{8},\\
\overline{\ell}_{10} \equiv -\overline{\ell}+4\overline{e}_{8}.
\end{cases}
\end{equation*}

\noindent
and we have 
\begin{equation*}
\begin{cases}
f^{\ast}(\overline{\ell})= \ell +2E_1+4E_2+3Q+3E_3+2L,\\
f^{\ast}(\overline{e}_{8})= e_{8} +\frac{3}{4}E_1+\frac{3}{2}E_2+\frac{5}{4}Q+E_3+\frac{1}{2}L+\frac{1}{3}E_6+\frac{2}{3}E_7,
\end{cases}
\end{equation*}

\noindent
with $(\overline{\ell})^2=9, (\overline{e}_{8})^2=\frac{11}{12}$, and $\overline{\ell}\cdot \overline{e}_{8}=3$. Hence, the ample cone $\Amp(S)$ of $S$ is the intersection of the regions

\begin{equation*}
\begin{cases}
3a+b>0,\\
3a+\tfrac{5}{4}b>0,\\
3a+\tfrac{11}{12}b>0,\\
3a+\tfrac{3}{4}b>0,\\
3a+\tfrac{13}{12}b>0,
\end{cases}
\qquad
\begin{cases}
3a+\tfrac{5}{6}b>0,\\
3a+\tfrac{7}{6}b>0,\\
3a+\tfrac{4}{3}b>0,\\
3a+\tfrac{2}{3}b>0.
\end{cases}
\end{equation*}

\noindent
This implies that the region of ample cone covered by the polarity cone $\Pol(U_{0})$ which is generated by $\overline{\ell}_{1}$, $\overline{\ell}_{2}$ and $\overline{\ell}_{3}$ is as follows.

\begin{figure}[H]
  \vspace{-1em}
\begin{subfigure}[t]{0.5\textwidth}\centering
\begin{tikzpicture}[scale=0.45, every node/.style={font=\scriptsize}, baseline=(origin)] 
  \coordinate (origin) at (0,0);
  \begin{scope}
    \clip (-5,-5) rectangle (5,5);

    \clip (-5,22.5) -- (5,-22.5) -- (5,5) -- (-5,5) -- cycle;

    \clip (-5,11.25) -- (5,-11.25) -- (5,5) -- (-5,5) -- cycle;

    \fill[yellow!100] (-5,-5) rectangle (5,5);
  \end{scope}

  \begin{scope}
    \clip (-5,-5) rectangle (5,5);
    \fill[gray!40] (0,0) -- (0,5) -- (5,5) -- (5,-5) -- (2.22,-5) -- cycle;
    
  \end{scope}

  \draw[->] (-5.5,0) -- (5.5,0) node[right] {$a$};
  \draw[->] (0,-5.5) -- (0,5.5) node[above] {$b$};

  \draw[thick] (0,0) -- (0,5);     
  \draw[thick] (0,0) -- (2.22,-5); 

  \begin{scope}
    \clip (-5,-5) rectangle (5,5);
    \draw[dotted,very thick] (-5,22.5) -- (5,-22.5);    
    \draw[dotted,very thick] (-5,11.25) -- (5,-11.25); 
  \end{scope}

\node[fill=white, inner sep=1pt, anchor=east] at (-1.5,2.5) {$b=-\frac{9}{4}a$};
  \node[fill=white, inner sep=1pt, anchor=east] at (1.8,-5.5) {$b=-\frac{9}{2}a$};
\end{tikzpicture}
\end{subfigure}
\vspace{-0.5em}
\caption{The region covered by $\Pol(U_{0})$ for type $\mathsf{D_5+A_2}$}
\vspace{-1em}
\end{figure}

Now for a $\mathbb{P}^{1}$-fibration on the minimal resolution of the surface,
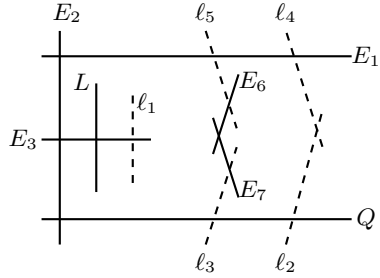
\begin{figure}[H]
  \vspace{-1em}
 \centering
\begin{subfigure}[t]{0.5\textwidth}\centering
  \begin{tikzpicture}[scale=0.48, every node/.style={font=\scriptsize}, baseline=(origin)]
\draw [thick] (-0.5,0)--(8,0); 
\draw [thick] (-0.5,4.5)--(8,4.5); 
\draw [thick] (0,-0.7)--(0,5.2); 
\draw [thick] (1,0.75)--(1,3.75); 
\draw [thick,dashed] (2,1)--(2,3.45); 
\draw [thick] (-0.5,2.2)--(2.5,2.2); 
\draw [thick] (4.9,4)--(4.2,1.8); 
\draw [thick] (4.2,2.8)--(4.9,0.6); 
\draw [thick,dashed] (4,5.2)--(4.85,2.51); 
\draw [thick,dashed] (4,-0.7)--(4.85,1.99); 
\draw [thick,dashed] (6.2,-0.7)--(7.2,3);
\draw [thick,dashed] (7.2,2)--(6.2,5.2);
\node at (8.4,0) {$Q$};
\node at (0.2,5.7) {$E_2$};
\node at (8.4,4.5) {$E_1$};
\node at (4,-1.2) {$\ell_3$};
\node at (6.2,-1.2) {$\ell_2$};
\node at (4,5.7) {$\ell_5$};
\node at (6.2,5.7) {$\ell_4$};
\node at (5.3,3.8) {$E_6$};
\node at (5.3,0.8) {$E_7$};
\node at (-1,2.2) {$E_3$};
\node at (0.6,3.8) {$L$};
\node at (2.4,3.2) {$\ell_1$};
\end{tikzpicture}
\end{subfigure}  
\vspace{-1em}
\end{figure}

\noindent
we consider a birational morphism from $\tilde{S}$ to $\mathbb{F}_{2}$ obtained by contracting $\ell_{1}, E_{3}, L, \ell_{4}, \ell_{5}, E_{6}$ and $E_{7}$. Then the cylinder defined by the images of $Q, E_{1}, E_{2}$ in $\mathbb{F}_{2}$ implies  that there exists a cylinder $U_{1} \cong \mathbb{A}^{1} \times \mathbb{A}^{1}_{\ast}$ in $S$ defined by $\overline{\ell_{1}}, \overline{\ell_{4}}, \overline{\ell_{5}}$. Therefore, $\Amp(S)$ is covered by the  union of the  relative interiors of $\Pol(U_{0})$ and $\Pol(U_{1})$.

\begin{figure}[H]
  \vspace{-1em}
\begin{subfigure}[t]{0.5\textwidth}\centering
\begin{tikzpicture}[scale=0.45, every node/.style={font=\scriptsize}, baseline=(origin)] 
  \coordinate (origin) at (0,0);
  \begin{scope}
    \clip (-5,-5) rectangle (5,5);

    \clip (-5,22.5) -- (5,-22.5) -- (5,5) -- (-5,5) -- cycle;

    \clip (-5,11.25) -- (5,-11.25) -- (5,5) -- (-5,5) -- cycle;

    \fill[yellow!100] (-5,-5) rectangle (5,5);
  \end{scope}

  \begin{scope}
    \clip (-5,-5) rectangle (5,5);
    \fill[gray!20] (0,0) -- (-1.1,5) -- (5,5) -- (5,-5) -- (3.3,-5) -- cycle;
    
  \end{scope}

  \draw[->] (-5.5,0) -- (5.5,0) node[right] {$a$};
  \draw[->] (0,-5.5) -- (0,5.5) node[above] {$b$};

  \draw[thick] (0,0) -- (-1.1,5);     
  \draw[thick] (0,0) -- (3.3,-5); 

  \begin{scope}
    \clip (-5,-5) rectangle (5,5);
    \draw[dotted,very thick] (-5,22.5) -- (5,-22.5);    
    \draw[dotted,very thick] (-5,11.25) -- (5,-11.25); 
  \end{scope}

\node[fill=white, inner sep=1pt, anchor=east] at (-0.7,5.5) {$b=-\frac{9}{2}a$};
  \node[fill=white, inner sep=1pt, anchor=east] at (6,-5.5) {$b=-\frac{3}{2}a$};
\end{tikzpicture}
\end{subfigure}
\vspace{-0.5em}
\caption{$\Pol(U_{1})$ for type $\mathsf{D_5+A_2}$}
\vspace{-1em}
\end{figure}

\medskip

\noindent
\textbf{Case 5.} $\mathsf{D_7}$\\
If we let $\Pic(S)=\mathbb{Z} \overline{\ell} \oplus \mathbb{Z}\overline{e}_{7}$, then the Mori cone $\overline{\NE}(S)$ of $S$ is generated by 

\begin{equation*}
\begin{cases}
\overline{\ell}_{1} \equiv \overline{e}_{7},\\
\overline{\ell}_{2} \equiv \overline{e}_{8} \equiv 2\overline{\ell}-5\overline{e}_{7},\\
\overline{\ell}_{3} \equiv \overline{L} \equiv \overline{\ell}-2\overline{e}_{7},
\end{cases}
\qquad
\begin{cases}
\overline{\ell}_{4} \equiv -\overline{\ell}+4\overline{e}_{7},\\
\overline{\ell}_{5} \equiv 3\overline{\ell}-8\overline{e}_{7},
\end{cases}
\end{equation*}

\noindent
and we have 
\begin{equation*}
\begin{cases}
f^{\ast}(\overline{\ell})= \ell +E_1+2E_2+3E_3+4E_4+5E_5+\frac{7}{2}Q+\frac{5}{2}E_6,\\
f^{\ast}(\overline{e}_{7})= e_{7} +\frac{1}{2}E_1+E_2+\frac{3}{2}E_3+2E_4+\frac{5}{2}E_5+\frac{5}{4}Q+\frac{7}{4}E_6,
\end{cases}
\end{equation*} 

\noindent
with $(\overline{\ell})^2=8, (\overline{e}_{7})^2=\frac{3}{4}$, and $\overline{\ell}\cdot \overline{e}_{7}=\frac{5}{2}$. Hence, the ample cone $\Amp(S)$ of $S$ is the intersection of the regions

\begin{equation*}
\begin{cases}
10a+3b>0,\\
14a+5b>0,\\
3a+b>0,
\end{cases}
\qquad
\begin{cases}
2a+\tfrac{1}{2}b>0,\\
4a+\tfrac{3}{2}b>0.
\end{cases}
\end{equation*}

\noindent
This implies the region of ample cone covered by $\Pol(U_{0})$ which is generated by $\overline{\ell}_{1} $, $\overline{\ell}_{2}$ and $\overline{\ell}_{3}$ is as follows.
\begin{figure}[H]
  \vspace{-1em}
\begin{subfigure}[t]{0.5\textwidth}\centering
\begin{tikzpicture}[scale=0.45, every node/.style={font=\scriptsize}, baseline=(origin)] 
  \coordinate (origin) at (0,0);
  \begin{scope}
    \clip (-5,-5) rectangle (5,5);

    \clip (-5,13.33) -- (5,-13.33) -- (5,5) -- (-5,5) -- cycle;

    \clip (-5,20) -- (5,-20) -- (5,5) -- (-5,5) -- cycle;

    \fill[yellow!100] (-5,-5) rectangle (5,5);
  \end{scope}

  \begin{scope}
    \clip (-5,-5) rectangle (5,5);
    \fill[gray!40] (0,0) -- (0,5) -- (5,5) -- (5,-5) -- (2.1,-5) -- cycle;
    
  \end{scope}

  \draw[->] (-5.5,0) -- (5.5,0) node[right] {$a$};
  \draw[->] (0,-5.5) -- (0,5.5) node[above] {$b$};

  \draw[thick] (0,0) -- (0,5);    
  \draw[thick] (0,0) -- (2.1,-5); 

  \begin{scope}
    \clip (-5,-5) rectangle (5,5);
    \draw[dotted,very thick] (-5,13.33) -- (5,-13.33);   
    \draw[dotted,very thick] (-5,20) -- (5,-20); 
  \end{scope}

\node[fill=white, inner sep=1pt, anchor=east] at (-1.5,2.5) {$b=-\frac{8}{3}a$};
  \node[fill=white, inner sep=1pt, anchor=east] at (1.5,-5.5) {$b=-4a$};
  \node[fill=white, inner sep=1pt, anchor=east] at (4.8,-5.5) {$b=-\frac{5}{2}a$};
\end{tikzpicture}
\end{subfigure}
\vspace{-0.5em}
\caption{The region covered by $\Pol(U_{0})$ for type $\mathsf{D_7}$}
\vspace{-1em}
\end{figure}

Now for a $\mathbb{P}^{1}$-fibration on the minimal resolution of the surface,
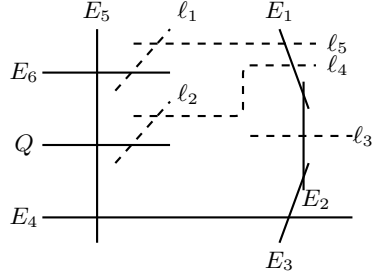
\begin{figure}[H]
  \vspace{-1em}
  \centering
\begin{subfigure}[t]{0.5\textwidth}\centering
  \begin{tikzpicture}[scale=0.48, every node/.style={font=\scriptsize}, baseline=(origin)]
\draw [thick] (-0.5,0)--(8,0); 
\draw [thick] (1,-0.7)--(1,5.2); 
\draw [thick] (6,-0.7)--(6.8,1.5); 
\draw [thick] (6,5.2)--(6.8,3); 
\draw [thick] (6.66,0.75)--(6.66,3.75); 
\draw [thick] (-0.5,4)--(3,4); 
\draw [thick] (-0.5,2)--(3,2); 
\draw [thick,dashed] (5.2,2.25)--(8,2.25); 
\draw [thick,dashed] (1.5,3.5)--(3,5.2); 
\draw [thick,dashed] (1.5,1.5)--(3,3.2); 
\draw [thick,dashed] (2,4.8)--(7,4.8); 
\draw [thick,dashed,rounded corners=3pt] (2,2.8)--(5,2.8)--(5,4.2)--(7,4.2); 
\node at (-1,0) {$E_4$};
\node at (1,5.7) {$E_5$};
\node at (-1,4) {$E_6$};
\node at (-1,2) {$Q$};
\node at (6,5.7) {$E_1$};
\node at (7.6,4.8) {$\ell_5$};
\node at (7.6,4.2) {$\ell_4$};
\node at (8.3,2.25) {$\ell_3$};
\node at (6,-1.2) {$E_3$};
\node at (7,0.5) {$E_2$};
\node at (3.5,5.7) {$\ell_1$};
\node at (3.5,3.4) {$\ell_2$};
\end{tikzpicture}
\vspace{-1em}
\end{subfigure}
\end{figure}

\noindent
we consider a birational morphism from $\tilde{S}$ to $\mathbb{F}_{2}$ obtained by contracting $\ell_{1}, E_{6}, E_{5}, Q, \ell_{4}, \ell_{3}$ and $E_{1}$. Then the cylinder defined by the images of $E_{2}, E_{3}, E_{4}$ in $\mathbb{F}_{2}$ implies  that there exists a cylinder $U_{1} \cong \mathbb{A}^{1} \times \mathbb{A}^{1}_{\ast}$ in $S$ defined by $\overline{\ell_{1}}, \overline{\ell_{3}}, \overline{\ell_{4}}$. Similarly, a cylinder in $\mathbb{F}_{2}$ obtained by contracting
$\ell_{2}, Q, E_{5}, E_{6}, \ell_{3}, \ell_{5} \ \text{and} \ E_{1}$ implies  that there exists a cylinder $U_{2}$  in $S$ which is isomorphic to $\mathbb{A}^{1} \times \mathbb{A}^{1}_{\ast}$ defined by $\overline{\ell_{2}}, \overline{\ell_{3}}, \overline{\ell_{5}}$. Therefore, $\Amp(S)$ is covered by the  union of the  relative interiors of $\Pol(U_{0}), \Pol(U_{1})$ and $\Pol(U_{2})$.

\begin{figure}[H]
  \vspace{-1em}
\begin{subfigure}[t]{0.5\textwidth}\centering
\begin{tikzpicture}[scale=0.45, every node/.style={font=\scriptsize}, baseline=(origin)] 
  \coordinate (origin) at (0,0);
  \begin{scope}
    \clip (-5,-5) rectangle (5,5);

    \clip (-5,13.33) -- (5,-13.33) -- (5,5) -- (-5,5) -- cycle;

    \clip (-5,20) -- (5,-20) -- (5,5) -- (-5,5) -- cycle;

    \fill[yellow!100] (-5,-5) rectangle (5,5);
  \end{scope}

  \begin{scope}
    \clip (-5,-5) rectangle (5,5);
    \fill[gray!20] (0,0) -- (-1.25,5) -- (5,5) -- (5,-5) -- (2.5,-5) -- cycle;
    
  \end{scope}

  \draw[->] (-5.5,0) -- (5.5,0) node[right] {$a$};
  \draw[->] (0,-5.5) -- (0,5.5) node[above] {$b$};

  \draw[thick] (0,0) -- (-1.25,5);    
  \draw[thick] (0,0) -- (2.5,-5); 

  \begin{scope}
    \clip (-5,-5) rectangle (5,5);
    \draw[dotted,very thick] (-5,13.33) -- (5,-13.33);   
    \draw[dotted,very thick] (-5,20) -- (5,-20); 
  \end{scope}

  \node[fill=white, inner sep=1pt, anchor=east] at (-0.7,5.5) {$b=-4a$};
  \node[fill=white, inner sep=1pt, anchor=east] at (5.0,-5.5) {$b=-2a$};
\end{tikzpicture}
\end{subfigure}\hfill
\begin{subfigure}[t]{0.5\textwidth}\centering
\begin{tikzpicture}[scale=0.45, every node/.style={font=\scriptsize}, baseline=(origin)] 
  \coordinate (origin) at (0,0);
  \begin{scope}
    \clip (-5,-5) rectangle (5,5);

    \clip (-5,13.33) -- (5,-13.33) -- (5,5) -- (-5,5) -- cycle;

    \clip (-5,20) -- (5,-20) -- (5,5) -- (-5,5) -- cycle;

    \fill[yellow!100] (-5,-5) rectangle (5,5);
  \end{scope}

  \begin{scope}
    \clip (-5,-5) rectangle (5,5);
    \fill[gray!60] (0,0) -- (1.875,-5) -- (2.5,-5) -- cycle;
    
  \end{scope}

  \draw[->] (-5.5,0) -- (5.5,0) node[right] {$a$};
  \draw[->] (0,-5.5) -- (0,5.5) node[above] {$b$};

  \draw[thick] (0,0) -- (1.875,-5);    
  \draw[thick] (0,0) -- (2.5,-5); 

  \begin{scope}
    \clip (-5,-5) rectangle (5,5);
    \draw[dotted,very thick] (-5,13.33) -- (5,-13.33);   
    \draw[dotted,very thick] (-5,20) -- (5,-20); 
  \end{scope}

\node[fill=white, inner sep=1pt, anchor=east] at (2.5,-5.5) {$b=-\frac{8}{3}a$};
  \node[fill=white, inner sep=1pt, anchor=east] at (4.9,-3.5) {$b=-2a$};
\end{tikzpicture}
\end{subfigure}
\vspace{-0.5em}
\caption{$\Pol(U_{1})$ and $\Pol(U_{2})$ for type $\mathsf{D_7}$}
\vspace{-1em}
\end{figure}

\medskip

\noindent
\textbf{Case 6.} $\mathsf{E_6+A_1}$\\
If we let $\Pic(S)=\mathbb{Z} \overline{\ell} \oplus \mathbb{Z}\overline{e}_{9}$ with $\overline{e}_{10} \equiv \frac{3}{5}\overline{\ell}$, then the Mori cone $\overline{\NE}(S)$ of $S$ is generated by 

\begin{equation*}
\begin{cases}
\overline{\ell}_{1} \equiv \overline{\ell}_{7} \equiv \overline{e}_{6} 
   \equiv \tfrac{1}{5}\overline{\ell},\\
\overline{\ell}_{2} \equiv \overline{e}_{7} 
   \equiv \tfrac{3}{5}\overline{\ell}-2\overline{e}_{9},\\
\overline{\ell}_{3} \equiv \overline{e}_{9},\\
\overline{\ell}_{4} \equiv \tfrac{2}{5}\overline{\ell}-\overline{e}_{9},
\end{cases}
\qquad
\begin{cases}
\overline{\ell}_{5} \equiv -\tfrac{1}{5}\overline{\ell}+2\overline{e}_{9},\\
\overline{\ell}_{6} \equiv -\tfrac{2}{5}\overline{\ell}+3\overline{e}_{9},\\
\overline{\ell}_{8} \equiv \tfrac{4}{5}\overline{\ell}-3\overline{e}_{9},
\end{cases}
\end{equation*}

\noindent
and we have 
\begin{equation*}
\begin{cases}
(f \circ g)^{\ast}(\overline{\ell})= \ell +5L_1+14E_1+15L_2+8E_2+10E_4+5E_5+7E_3+4L_3,\\
(f \circ g)^{\ast}(\overline{e}_{9})= e_{9} +\frac{4}{3}L_1+3E_1+3L_2+\frac{5}{3}E_2+2E_4+E_5+\frac{4}{3}E_3+\frac{2}{3}L_3+\frac{1}{2}E_8,
\end{cases}
\end{equation*}

\noindent
with $(\overline{\ell})^2=25, (\overline{e}_{9})^2=\frac{5}{6}$, and $\overline{\ell}\cdot \overline{e}_{9}=5$. Hence, the ample cone $\Amp(S)$ of $S$ is the intersection of the regions

\begin{equation*}
\begin{cases}
5a+b>0,\\
5a+\tfrac{4}{3}b>0,\\
5a+\tfrac{5}{6}b>0,\\
5a+\tfrac{7}{6}b>0,
\end{cases}
\qquad
\begin{cases}
5a+\tfrac{2}{3}b>0,\\
5a+\tfrac{1}{2}b>0,\\
5a+\tfrac{3}{2}b>0.
\end{cases}
\end{equation*}

\noindent
This implies that the region of ample cone covered by $\Pol(U_{0})$ which is generated by $\overline{\ell}_{1}$, $\overline{\ell}_{2}$, $\overline{\ell}_{3}$ and $\overline{e}_{10}$ is as follows.
\begin{figure}[H]
  \vspace{-1em}
\begin{subfigure}[t]{0.5\textwidth}\centering
\begin{tikzpicture}[scale=0.45, every node/.style={font=\scriptsize}, baseline=(origin)] 
  \coordinate (origin) at (0,0);
  \begin{scope}
    \clip (-5,-5) rectangle (5,5);

    \clip (-5,16.67) -- (5,-16.67) -- (5,5) -- (-5,5) -- cycle;

    \clip (-5,50) -- (5,-50) -- (5,5) -- (-5,5) -- cycle;

    \fill[yellow!100] (-5,-5) rectangle (5,5);
  \end{scope}

  \begin{scope}
    \clip (-5,-5) rectangle (5,5);
    \fill[gray!40] (0,0) -- (0,5) -- (5,5) -- (5,-5) -- (1.5,-5) -- cycle;
    
  \end{scope}

  \draw[->] (-5.5,0) -- (5.5,0) node[right] {$a$};
  \draw[->] (0,-5.5) -- (0,5.5) node[above] {$b$};

  \draw[thick] (0,0) -- (0,5);    
  \draw[thick] (0,0) -- (1.5,-5);  

  \begin{scope}
    \clip (-5,-5) rectangle (5,5);
    \draw[dotted,very thick] (-5,16.67) -- (5,-16.67);    
    \draw[dotted,very thick] (-5,50) -- (5,-50); 
  \end{scope}

\node[fill=white, inner sep=1pt, anchor=east] at (-1.2,2.5) {$b=-\frac{10}{3}a$};
  \node[fill=white, inner sep=1pt, anchor=east] at (1.2,-5.5) {$b=-10a$};
\end{tikzpicture}
\end{subfigure}
\vspace{-0.5em}
\caption{The region covered by $\Pol(U_{0})$ for type $\mathsf{E_6+A_1}$}
\vspace{-1em}
\end{figure}

Now for a $\mathbb{P}^{1}$-fibration on the minimal resolution of the surface,
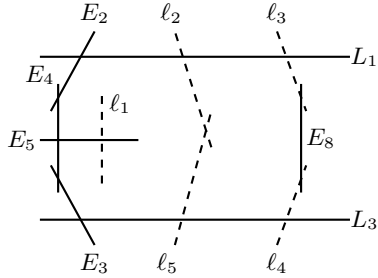
\begin{figure}[H]
  \vspace{-1em}
 \centering
\begin{subfigure}[t]{0.5\textwidth}\centering
  \begin{tikzpicture}[scale=0.48, every node/.style={font=\scriptsize}, baseline=(origin)]
\draw [thick] (-0.5,0)--(8,0); 
\draw [thick] (-0.5,4.5)--(8,4.5); 
\draw [thick] (-0.2,3)--(1,5.2); 
\draw [thick] (-0.2,1.5)--(1,-0.7); 
\draw [thick] (0,0.75)--(0,3.75); 
\draw [thick] (-0.5,2.2)--(2.2,2.2); 
\draw [thick,dashed] (6,-0.7)--(6.8,1.5); 
\draw [thick,dashed] (6,5.2)--(6.8,3); 
\draw [thick,dashed] (1.2,1)--(1.2,3.45); 
\draw [thick,dashed] (3.2,-0.7)--(4.2,3);
\draw [thick,dashed] (4.2,2)--(3.2,5.2);
\draw [thick] (6.66,0.75)--(6.66,3.75); 
\node at (8.4,0) {$L_3$};
\node at (1,5.7) {$E_2$};
\node at (8.4,4.5) {$L_1$};
\node at (1,-1.2) {$E_3$};
\node at (6,5.7) {$\ell_3$};
\node at (7.2,2.25) {$E_8$};
\node at (6,-1.2) {$\ell_4$};
\node at (3.1,5.7) {$\ell_2$};
\node at (3.,-1.2) {$\ell_5$};
\node at (-1,2.2) {$E_5$};
\node at (1.7,3.2) {$\ell_1$};
\node at (-0.5,4) {$E_4$};
\end{tikzpicture}
\end{subfigure} 
\vspace{-1em}
\end{figure}

\noindent
we consider a birational morphism from $\tilde{S}$ to $\mathbb{F}_{2}$ obtained by contracting $\ell_{1}, E_{5}, E_{4}, E_{3}, \ell_{4}, E_{8}$ and $\ell_{5}$. Then the cylinder defined by the images of $E_{2}, L_{1}, L_{3}$ in $\mathbb{F}_{2}$ implies  that there exists a cylinder $U_{1} \cong \mathbb{A}^{1} \times \mathbb{A}^{1}_{\ast}$ in $S$ defined by $\overline{\ell_{1}}, \overline{\ell_{4}}, \overline{\ell_{5}}$. Therefore, $\Amp(S)$ is covered by the  union of the  relative interiors of $\Pol(U_{0})$ and $\Pol(U_{1})$.

\begin{figure}[H]
  \vspace{-1em}
\begin{subfigure}[t]{0.5\textwidth}\centering
\begin{tikzpicture}[scale=0.46, every node/.style={font=\scriptsize}, baseline=(origin)] 
  \coordinate (origin) at (0,0);
  \begin{scope}
    \clip (-5,-5) rectangle (5,5);

    \clip (-5,16.67) -- (5,-16.67) -- (5,5) -- (-5,5) -- cycle;

    \clip (-5,50) -- (5,-50) -- (5,5) -- (-5,5) -- cycle;

    \fill[yellow!100] (-5,-5) rectangle (5,5);
  \end{scope}

  \begin{scope}
    \clip (-5,-5) rectangle (5,5);
    \fill[gray!20] (0,0) -- (-0.5,5) -- (5,5) -- (5,-5) -- (2,-5) -- cycle;
    
  \end{scope}

  \draw[->] (-5.5,0) -- (5.5,0) node[right] {$a$};
  \draw[->] (0,-5.5) -- (0,5.5) node[above] {$b$};

  \draw[thick] (0,0) -- (-0.5,5);    
  \draw[thick] (0,0) -- (2,-5);  

  \begin{scope}
    \clip (-5,-5) rectangle (5,5);
    \draw[dotted,very thick] (-5,16.67) -- (5,-16.67);    
    \draw[dotted,very thick] (-5,50) -- (5,-50); 
  \end{scope}

\node[fill=white, inner sep=1pt, anchor=east] at (-0.3,5.5) {$b=-10a$};
  \node[fill=white, inner sep=1pt, anchor=east] at (4.7,-5.5) {$b=-\frac{5}{2}a$};
\end{tikzpicture}
\end{subfigure}
\vspace{-0.5em}
\caption{$\Pol(U_{1})$ for type $\mathsf{E_6+A_1}$}
\vspace{-1em}
\end{figure}

\medskip

\noindent
\textbf{Case 7.} $\mathsf{A_4+A_3}$\\
If we let $\Pic(S)=\mathbb{Z} \overline{\ell} \oplus \mathbb{Z}\overline{e}_{13}$ with $\overline{e}_5 \equiv \frac{1}{5}\overline{\ell}$ and $\overline{e}_{9} \equiv \frac{1}{3}\overline{\ell}-4\overline{e}_{13}$, then the Mori cone $\overline{\NE}(S)$ of $S$ is generated by  

\begin{equation*}
\begin{cases}
\overline{\ell}_{1} \equiv \overline{e}_{9} \equiv \tfrac{1}{3}\overline{\ell}-4\overline{e}_{13},\\
\overline{\ell}_{2} \equiv \overline{e}_{13},\\
\overline{\ell}_{3} \equiv \overline{\ell}_{12} \equiv \tfrac{1}{15}\overline{\ell},\\
\overline{\ell}_{4} \equiv -\tfrac{1}{5}\overline{\ell}+4\overline{e}_{13},\\
\overline{\ell}_{5} \equiv \tfrac{2}{15}\overline{\ell}-\overline{e}_{13},\\
\overline{\ell}_{6} \equiv -\tfrac{2}{15}\overline{\ell}+3\overline{e}_{13},
\end{cases}
\qquad
\begin{cases}
\overline{\ell}_{7} \equiv \tfrac{1}{5}\overline{\ell}-2\overline{e}_{13},\\
\overline{\ell}_{8} \equiv -\tfrac{1}{15}\overline{\ell}+2\overline{e}_{13},\\
\overline{\ell}_{9} \equiv \tfrac{4}{15}\overline{\ell}-3\overline{e}_{13},\\
\overline{\ell}_{10} \equiv \tfrac{6}{15}\overline{\ell}-5\overline{e}_{13},\\
\overline{\ell}_{11} \equiv -\tfrac{4}{15}\overline{\ell}+5\overline{e}_{13}.
\end{cases}
\end{equation*}

\noindent
In addition, we have $(\overline{\ell})^2=225, (\overline{e}_{13})^2=\frac{19}{20}$ and $\overline{\ell}\cdot \overline{e}_{13}=15$ with 
\begin{equation*}
\begin{cases}
(f \circ g)^{\ast}(\overline{\ell})= \ell +9E_1+27E_3+45E_4+63E_2+144L_1+80L_2+15E_8+10E_7+5E_6,\\
(f \circ g)^{\ast}(\overline{e}_{13})= e_{13} +\frac{3}{5}E_1+\frac{9}{5}E_3+3E_4+\frac{21}{5}E_2+\frac{48}{5}L_1+\frac{27}{5}L_2+\frac{6}{5}E_8+\frac{4}{5}E_7+\frac{2}{5}E_6\\
\quad+\frac{1}{4}E_{10}+\frac{1}{2}E_{11}+\frac{3}{4}E_{12}.\\
\end{cases}
\end{equation*} 

\noindent
Hence, the ample cone of $S$ is the intersection of the regions

\begin{equation*}
\begin{cases}
15a+\tfrac{6}{5}b>0,\\
15a+\tfrac{19}{20}b>0,\\
15a+b>0,\\
15a+\tfrac{4}{5}b>0,\\
15a+\tfrac{21}{20}b>0,\\
15a+\tfrac{17}{20}b>0,
\end{cases}
\qquad
\begin{cases}
15a+\tfrac{11}{10}b>0,\\
15a+\tfrac{9}{10}b>0,\\
15a+\tfrac{23}{20}b>0,\\
15a+\tfrac{5}{4}b>0,\\
15a+\tfrac{3}{4}b>0.
\end{cases}
\end{equation*}

\noindent
This implies that the region of ample cone covered by $\Pol(U_{0})$ which is generated by is generated by $\overline{\ell}_{1}$, $\overline{\ell}_{2}$ and $\overline{e}_5$ is as follows.
 \begin{figure}[H]
  \vspace{-1em}
\begin{subfigure}[t]{0.5\textwidth}\centering
\begin{tikzpicture}[scale=0.45, every node/.style={font=\scriptsize}, baseline=(origin)] %
  \begin{scope}
    \clip (-5,-5) rectangle (5,5);

    \clip (-5,100) -- (5,-100) -- (5,5) -- (-5,5) -- cycle;

    \clip (-5,60) -- (5,-60) -- (5,5) -- (-5,5) -- cycle;

    \fill[yellow!100] (-5,-5) rectangle (5,5);
  \end{scope}

  \begin{scope}
    \clip (-5,-5) rectangle (5,5);
    \fill[gray!40] (0,0) -- (0,5) -- (5,5) -- (5,-5) -- (0.41,-5) -- cycle;
    
  \end{scope}

  \draw[->] (-5.5,0) -- (5.5,0) node[right] {$a$};
  \draw[->] (0,-5.5) -- (0,5.5) node[above] {$b$};

  \draw[thick] (0,0) -- (0,5);     
  \draw[thick] (0,0) -- (0.41,-5);  

  \begin{scope}
    \clip (-5,-5) rectangle (5,5);
    \draw[dotted,very thick] (-5,100) -- (5,-100);    
    \draw[dotted,very thick] (-5,60) -- (5,-60); 
  \end{scope}

  \node[fill=white, inner sep=1pt, anchor=east] at (-0.8,2.5) {$b=-12a$};
  \node[fill=white, inner sep=1pt, anchor=east] at (-0.5,-2.5) {$b=-20a$};
\end{tikzpicture}
\end{subfigure}
\caption{The region covered by $\Pol(U_{0})$ for type $\mathsf{A_4+A_3}$}
\vspace{-1em}
\end{figure}

Now for a $\mathbb{P}^{1}$-fibration on the minimal resolution of the surface,

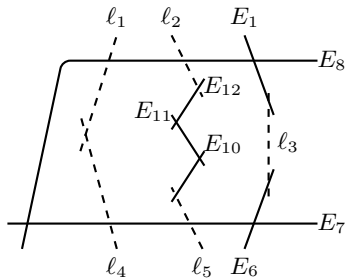
\begin{figure}[H]
  \vspace{-1em}
  \centering
\begin{subfigure}[t]{0.5\textwidth}\centering
  \begin{tikzpicture}[scale=0.48, every node/.style={font=\scriptsize}, baseline=(origin)]
\draw [thick] (-0.5,0)--(8,0); 
\draw [thick,rounded corners=3pt] (-0.1,-0.7)--(1,4.5)--(8,4.5); 
\draw [thick] (6,-0.7)--(6.8,1.5); 
\draw [thick] (6,5.2)--(6.8,3); 
\draw [thick] (4.9,4)--(4,2.6); 
\draw [thick] (4,3)--(4.9,1.6); 
\draw [thick] (4.9,2)--(4,0.6); 
\draw [thick,dashed] (4,5.2)--(4.85,3.5); 
\draw [thick,dashed] (4.85,-0.7)--(4,1); 
\draw [thick,dashed] (2.5,-0.7)--(1.5,3);
\draw [thick,dashed] (1.5,2)--(2.5,5.2);
\draw [thick, dashed] (6.66,0.75)--(6.66,3.75); 
\node at (8.4,0) {$E_7$};
\node at (8.4,4.5) {$E_8$};
\node at (4.85,-1.2) {$\ell_5$};
\node at (2.5,-1.2) {$\ell_4$};
\node at (4,5.7) {$\ell_2$};
\node at (2.5,5.7) {$\ell_1$};
\node at (5.4,3.8) {$E_{12}$};
\node at (3.5,3.1) {$E_{11}$};
\node at (5.4,2.2) {$E_{10}$};
\node at (6,5.7) {$E_1$};
\node at (7.2,2.25) {$\ell_3$};
\node at (6,-1.2) {$E_6$};
\end{tikzpicture}
\end{subfigure} 
\vspace{-1em}
\end{figure}

\noindent
We consider a birational map from $\tilde{S}$ to $\mathbb{F}_{2}$ obtained by a sequence of blow-ups and downs as follows.
The marked dots in 
\tikz[baseline=(char.base)]{\node[draw,circle,inner sep=0.7pt](char){\scriptsize 2};}, 
\tikz[baseline=(char.base)]{\node[draw,circle,inner sep=0.7pt](char){\scriptsize 4};}, 
and 
\tikz[baseline=(char.base)]{\node[draw,circle,inner sep=0.7pt](char){\scriptsize 6};}
indicate the centers of the blow-ups.

\begin{figure}[H]
  \centering
\begin{subfigure}[t]{0.28\textwidth}\centering
\begin{tikzpicture}[scale=0.32, every node/.style={font=\scriptsize}, baseline=(origin)]
\draw [thick] (-0.5,0)--(8,0); 
\draw [thick,rounded corners=3pt] (-0.1,-0.7)--(1,4.5)--(8,4.5); 
\draw [thick] (6,-0.7)--(6.8,1.5); 
\draw [thick] (6,5.2)--(6.8,3); 
\draw [thick] (4.9,4)--(4,2.6); 
\draw [thick] (4,3)--(4.9,1.6); 
\draw [thick] (4.9,2)--(4,0.6); 
\draw [thick] (0.05,-0.7)--(0.05,5.2); 
\draw [thick,dashed] (4.85,-0.7)--(4,1); 
\draw [thick,dashed] (2.5,-0.7)--(1.5,3);
\node[draw, circle, inner sep=0.5pt] at (-0.5,6){\scriptsize 1};
\node at (-0.4,5){$0$};
\end{tikzpicture}
\end{subfigure}\hfill
\begin{subfigure}[t]{0.28\textwidth}\centering
\begin{tikzpicture}[scale=0.32, every node/.style={font=\scriptsize}, baseline=(origin)]
\draw [thick] (-0.5,0)--(8,0); 
\draw [thick,rounded corners=3pt] (-0.1,-0.7)--(1,4.5)--(8,4.5); 
\draw [thick] (6,-0.7)--(6.8,1.5); 
\draw [thick] (6,5.2)--(6.8,3); 
\draw [thick] (0.05,-0.7)--(0.05,5.2); 
\node[draw, circle, inner sep=0.5pt] at (-0.5,6){\scriptsize 2};
\node at (8.7,0) {$3$};
\node at (8.7,4.5) {$-2$};
\node at (6,5.7) {$-2$};
\node at (6,-1.2) {$-2$};
\node at (-0.4,5){$0$};
\node at (0.05,0){$\bullet$};
\end{tikzpicture}
\end{subfigure}\hfill
\begin{subfigure}[t]{0.28\textwidth}\centering
\begin{tikzpicture}[scale=0.32, every node/.style={font=\scriptsize}, baseline=(origin)]
\draw [thick] (-0.5,0)--(8,0); 
\draw [thick,rounded corners=3pt] (-0.5,4.5)--(8,4.5); 
\draw [thick] (6,-0.7)--(6.8,1.5); 
\draw [thick] (-0.5,2.5)--(3,2.5); 
\draw [thick] (6,5.2)--(6.8,3); 
\draw [thick] (0.05,-0.7)--(0.05,5.2); 
\node[draw, circle, inner sep=0.5pt] at (-0.5,6){\scriptsize 3};
\node at (8.7,0) {$2$};
\node at (8.7,4.5) {$-3$};
\node at (2.8,2) {$-1$};
\node at (6,5.7) {$-2$};
\node at (6,-1.2) {$-2$};
\node at (-0.8,5){$-1$};
\end{tikzpicture}
\end{subfigure}\hfill
\begin{subfigure}[t]{0.28\textwidth}\centering
\begin{tikzpicture}[scale=0.32, every node/.style={font=\scriptsize}, baseline=(origin)]
\draw [thick] (-0.5,0)--(8,0); 
\draw [thick,rounded corners=3pt] (-0.5,4.5)--(8,4.5); 
\draw [thick] (6,-0.7)--(6.8,1.5); 
\draw [thick] (6,5.2)--(6.8,3); 
\draw [thick] (0.05,-0.7)--(0.05,5.2); 
\node[draw, circle, inner sep=0.5pt] at (-0.5,6){\scriptsize 4};
\node at (8.7,0) {$2$};
\node at (8.7,4.5) {$-3$};
\node at (6,5.7) {$-2$};
\node at (6,-1.2) {$-2$};
\node at (-0.4,5){$0$};
\node at (0.05,-0.02){$\bullet$};
\end{tikzpicture}
\end{subfigure}\hfill
\begin{subfigure}[t]{0.28\textwidth}\centering
\begin{tikzpicture}[scale=0.32, every node/.style={font=\scriptsize}, baseline=(origin)]
\draw [thick] (0.5,0)--(8,0); 
\draw [thick,rounded corners=3pt] (-0.5,4.5)--(8,4.5); 
\draw [thick] (6,-0.7)--(6.8,1.5); 
\draw [thick] (6,5.2)--(6.8,3); 
\draw [thick] (0.05,0.8)--(0.05,5.2); 
\draw [thick] (-0.5,2.5)--(1.3,-0.7); 
\node[draw, circle, inner sep=0.5pt] at (-0.5,6){\scriptsize 5};
\node at (8.7,0) {$1$};
\node at (8.7,4.5) {$-3$};
\node at (6,5.7) {$-2$};
\node at (6,-1.2) {$-2$};
\node at (-0.8,5){$-1$};
\node at (2,-0.7){$-1$};
\end{tikzpicture}
\end{subfigure}\hfill
\begin{subfigure}[t]{0.28\textwidth}\centering
\begin{tikzpicture}[scale=0.32, every node/.style={font=\scriptsize}, baseline=(origin)]
\draw [thick] (-0.5,0)--(8,0); 
\draw [thick,rounded corners=3pt] (-0.5,4.5)--(8,4.5); 
\draw [thick] (6,-0.7)--(6.8,1.5); 
\draw [thick] (6,5.2)--(6.8,3); 
\draw [thick] (0.05,-0.7)--(0.05,5.2); 
\node[draw, circle, inner sep=0.5pt] at (-0.5,6){\scriptsize 6};
\node at (8.7,0) {$1$};
\node at (8.7,4.5) {$-2$};
\node at (6,5.7) {$-2$};
\node at (6,-1.2) {$-2$};
\node at (-0.4,5){$0$};
\node at (0.05,-0.02){$\bullet$};
\end{tikzpicture}
\end{subfigure}\hfill
\begin{subfigure}[t]{0.28\textwidth}\centering
\begin{tikzpicture}[scale=0.32, every node/.style={font=\scriptsize}, baseline=(origin)]
\draw [thick] (0.5,0)--(8,0); 
\draw [thick,rounded corners=3pt] (-0.5,4.5)--(8,4.5); 
\draw [thick] (6,-0.7)--(6.8,1.5); 
\draw [thick] (6,5.2)--(6.8,3); 
\draw [thick] (0.05,0.8)--(0.05,5.2); 
\draw [thick] (-0.5,2.5)--(1.3,-0.7); 
\node[draw, circle, inner sep=0.5pt] at (-0.5,6){\scriptsize 7};
\node at (8.7,0) {$0$};
\node at (8.7,4.5) {$-2$};
\node at (6,5.7) {$-2$};
\node at (6,-1.2) {$-2$};
\node at (-0.8,5){$-1$};
\node at (2,-0.7){$-1$};
\end{tikzpicture}
\end{subfigure}\hfill
\begin{subfigure}[t]{0.28\textwidth}\centering
\begin{tikzpicture}[scale=0.32, every node/.style={font=\scriptsize}, baseline=(origin)]
\draw [thick] (-0.5,0)--(8,0); 
\draw [thick,rounded corners=3pt] (-0.5,4.5)--(8,4.5); 
\draw [thick] (6,-0.7)--(6.8,1.5); 
\draw [thick] (6,5.2)--(6.8,3); 
\draw [thick] (0.05,-0.7)--(0.05,5.2); 
\node[draw, circle, inner sep=0.5pt] at (-0.5,6){\scriptsize 8};
\node at (8.7,0) {$0$};
\node at (8.7,4.5) {$-1$};
\node at (6,5.7) {$-2$};
\node at (6,-1.2) {$-2$};
\node at (-0.4,5){$0$};
\end{tikzpicture}
\end{subfigure}\hfill
\begin{subfigure}[t]{0.28\textwidth}\centering
\begin{tikzpicture}[scale=0.32, every node/.style={font=\scriptsize}, baseline=(origin)]
\draw [thick] (-0.5,0)--(8,0); 
\draw [thick,rounded corners=3pt] (-0.5,4.5)--(8,4.5); 
\draw [thick] (6,-0.7)--(6.8,1.5); 
\draw [thick] (0.05,-0.7)--(0.05,5.2); 
\node[draw, circle, inner sep=0.5pt] at (-0.5,6){\scriptsize 9};
\node at (8.7,0) {$0$};
\node at (8.7,4.5) {$-1$};
\node at (6,-1.2) {$-2$};
\node at (-0.4,5){$1$};
\end{tikzpicture}
\end{subfigure}
\vspace{-0.5em}
\caption{Construction of $U_{1}$ for type $\mathsf{A_{4}+A_{3}}$}
\vspace{-1em}
\end{figure}

\noindent
Then, by contracting the $(-1)$-curve in the ninth configuration, we obtain a cylinder in $\mathbb{F}_{2}$ defined by two sections with self-intersection numbers $-2$ and $2$, respectively, and one fiber. This implies  that there exists a cylinder $U_{1}$ in $S$ defined by $\overline{\ell_{4}}, \overline{\ell_{5}}$ and $\overline{F} \equiv \frac{2}{15} \overline{\ell}$, which is the image of the fiber $F$ on $S$. Therefore, $\Amp(S)$ is covered by the  union of the  relative interiors of $\Pol(U_{0})$ and $\Pol(U_{1})$.

 \begin{figure}[H]
  \vspace{-1em}
\begin{subfigure}[t]{0.5\textwidth}\centering
\begin{tikzpicture}[scale=0.45, every node/.style={font=\scriptsize}, baseline=(origin)] %
  \begin{scope}
    \clip (-5,-5) rectangle (5,5);

    \clip (-5,100) -- (5,-100) -- (5,5) -- (-5,5) -- cycle;

    \clip (-5,60) -- (5,-60) -- (5,5) -- (-5,5) -- cycle;

    \fill[yellow!100] (-5,-5) rectangle (5,5);
  \end{scope}

  \begin{scope}
    \clip (-5,-5) rectangle (5,5);
    \fill[gray!20] (0,0) -- (-0.25,5) -- (5,5) -- (5,-5) -- (0.66,-5) -- cycle;
    
  \end{scope}

  \draw[->] (-5.5,0) -- (5.5,0) node[right] {$a$};
  \draw[->] (0,-5.5) -- (0,5.5) node[above] {$b$};

  \draw[thick] (0,0) -- (-0.25,5);     
  \draw[thick] (0,0) -- (0.66,-5);  

  \begin{scope}
    \clip (-5,-5) rectangle (5,5);
    \draw[dotted,very thick] (-5,100) -- (5,-100);    
    \draw[dotted,very thick] (-5,60) -- (5,-60); 
  \end{scope}

  \node[fill=white, inner sep=1pt, anchor=east] at (3.3,2.5) {$b=-20a$};
  \node[fill=white, inner sep=1pt, anchor=east] at (3.9,-2.5) {$b=-\frac{15}{2}a$};
\end{tikzpicture}
\end{subfigure}
\vspace{-0.5em}
\caption{$\Pol(U_{1})$ for type $\mathsf{A_4+A_3}$}
\vspace{-1.5em}
\end{figure}

\medskip

\noindent
\textbf{Case 8.} $\mathsf{A_4+A_2+A_1}$\\
If we let $\Pic(S)=\mathbb{Z} \overline{\ell} \oplus \mathbb{Z}\overline{e}_{13}$ with $\overline{e}_5 \equiv \frac{1}{5}\overline{\ell}$ and $\overline{e}_{11} \equiv \frac{1}{9}\overline{\ell}-\frac{2}{3}\overline{e}_{13}$, then the Mori cone $\overline{\NE}(S)$  of $S$ is generated by   
\begin{equation*}
\begin{cases}
\overline{\ell}_{1} \equiv \tfrac{1}{9}\overline{\ell}-\tfrac{2}{3}\overline{e}_{13},\\
\overline{\ell}_{2} \equiv \overline{e}_{13},\\
\overline{\ell}_{3} \equiv \overline{\ell}_{13} \equiv \overline{\ell}_{14} 
   \equiv \tfrac{1}{15}\overline{\ell},\\
\overline{\ell}_{4} \equiv \tfrac{1}{45}\overline{\ell}+\tfrac{2}{3}\overline{e}_{13},\\
  \overline{\ell}_{5} \equiv \tfrac{2}{15}\overline{\ell}-\overline{e}_{13},\\ 

\end{cases}
\qquad
\begin{cases}
\overline{\ell}_{6} \equiv \tfrac{4}{45}\overline{\ell}-\tfrac{1}{3}\overline{e}_{13},\\
\overline{\ell}_{7} \equiv -\tfrac{1}{45}\overline{\ell}+\tfrac{4}{3}\overline{e}_{13},\\
\overline{\ell}_{8} \equiv \tfrac{1}{5}\overline{\ell}-2\overline{e}_{13},\\
\overline{\ell}_{9} \equiv -\tfrac{1}{15}\overline{\ell}+2\overline{e}_{13},\\
\overline{\ell}_{10} \equiv \tfrac{2}{45}\overline{\ell}+\tfrac{1}{3}\overline{e}_{13},\\
\end{cases}
\qquad
\begin{cases}
\overline{\ell}_{11} \equiv \tfrac{7}{45}\overline{\ell}-\tfrac{4}{3}\overline{e}_{13},\\
\overline{\ell}_{12} \equiv -\tfrac{2}{45}\overline{\ell}+\tfrac{5}{3}\overline{e}_{13},\\
\overline{\ell}_{15} \equiv \tfrac{8}{45}\overline{\ell}-\tfrac{5}{3}\overline{e}_{13}.
\end{cases}
\end{equation*}

\noindent
In addition, we have $(\overline{\ell})^2=225, (\overline{e}_{13})^2=\frac{7}{10}$ and $\overline{\ell}\cdot \overline{e}_{13}=15$ with 
\begin{equation*}
\begin{cases}
(f \circ g)^{\ast}(\overline{\ell})= \ell +9E_1+27E_3+45E_4+63E_2+144L_1+80L_2+15E_8+10E_7+5E_6,\\
(f \circ g)^{\ast}(\overline{e}_{13})= e_{13} +\frac{3}{5}E_1+\frac{9}{5}E_3+3E_4+\frac{21}{5}E_2+\frac{48}{5}L_1+\frac{27}{5}L_2+\frac{6}{5}E_8+\frac{4}{5}E_7+\frac{2}{5}E_6+\frac{1}{2}E_{12}.\\
\end{cases}
\end{equation*} 

\noindent
Hence, the ample cone of $S$ is the intersection of the regions

\begin{equation*}
\begin{cases}
15a+\tfrac{6}{5}b>0,\\
15a+\tfrac{7}{10}b>0,\\
15a+b>0,\\ 
15a+\tfrac{4}{5}b>0,\\
15a+\tfrac{13}{10}b>0,\\
15a+\tfrac{11}{10}b>0,\\
15a+\tfrac{3}{5}b>0,
\end{cases}
\qquad
\begin{cases}
15a+\tfrac{8}{5}b>0,\\
15a+\tfrac{2}{5}b>0,\\
15a+\tfrac{9}{10}b>0,\\
15a+\tfrac{7}{5}b>0,\\
15a+\tfrac{1}{2}b>0,\\
15a+\tfrac{3}{2}b>0.
\end{cases}
\end{equation*}

\noindent
This implies that the region of ample cone covered by $\Pol(U_{0})$ which is generated by $\overline{\ell}_{2}$, $\overline{e}_5$ and $\overline{e}_{11}$ is as follows.

\begin{figure}[H]
  \vspace{-2em}
\begin{subfigure}[t]{0.5\textwidth}\centering
\begin{tikzpicture}[scale=0.45, every node/.style={font=\scriptsize}, baseline=(origin)] %
  \begin{scope}
    \clip (-5,-5) rectangle (5,5);

    \clip (-5,187.5) -- (5,-187.5) -- (5,5) -- (-5,5) -- cycle;

    \clip (-5,46.875) -- (5,-46.875) -- (5,5) -- (-5,5) -- cycle;

    \fill[yellow!100] (-5,-5) rectangle (5,5);
  \end{scope}

  \begin{scope}
    \clip (-5,-5) rectangle (5,5);
    \fill[gray!40] (0,0) -- (0,5) -- (5,5) -- (5,-5) -- (0.83,-5) -- cycle;
    
  \end{scope}

  \draw[->] (-5.5,0) -- (5.5,0) node[right] {$a$};
  \draw[->] (0,-5.5) -- (0,5.5) node[above] {$b$};

  \draw[thick] (0,0) -- (0,5);     
  \draw[thick] (0,0) -- (0.83,-5);  

  \begin{scope}
    \clip (-5,-5) rectangle (5,5);
    \draw[dotted,very thick] (-5,187.5) -- (5,-187.5);    
    \draw[dotted,very thick] (-5,46.875) -- (5,-46.875); 
  \end{scope}

  \node[fill=white, inner sep=1pt, anchor=east] at (-0.8,2.5) {$b=-\tfrac{75}{8}\,a$};
  \node[fill=white, inner sep=1pt, anchor=east] at (-0.5,-2.5) {$b=-\tfrac{75}{2}\,a$};
  \node[fill=white, inner sep=1pt, anchor=east] at (2.6,-5.5) {$b=-6a$};
\end{tikzpicture}
\end{subfigure}
\vspace{-0.5em}
\caption{The region covered by $\Pol(U_{0})$ for type $\mathsf{A_4+A_2+A_1}$}
\vspace{-1em}
\end{figure}

Now for a $\mathbb{P}^{1}$-fibration on the minimal resolution of the surface,
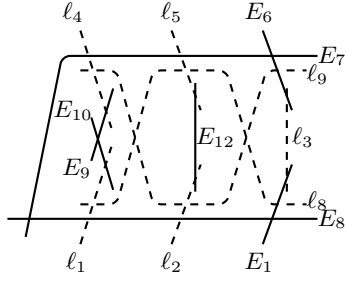
\begin{figure}[H]
  \vspace{-1em}
  \centering
\begin{subfigure}[t]{0.5\textwidth}\centering
  \begin{tikzpicture}[scale=0.48, every node/.style={font=\scriptsize}, baseline=(origin)]
\draw [thick] (-0.5,0)--(8,0); 
\draw [thick,rounded corners=3pt] (0,-0.5)--(1,4.5)--(8,4.5); 
\draw [thick] (6.5,-0.7)--(7.3,1.5); 
\draw [thick] (6.5,5.2)--(7.3,3); 
\draw [thick,dashed] (4,5.2)--(4.8,3); 
\draw [thick,dashed] (4,-0.7)--(4.8,1.5); 
\draw [thick] (4.65,0.75)--(4.65,3.75); 
\draw [thick] (2.4,3.6)--(1.8,1.6); 
\draw [thick] (1.8,2.8)--(2.4,0.8); 
\draw [thick,dashed] (1.5,5.2)--(2.35,2.51); 
\draw [thick,dashed] (1.5,-0.7)--(2.35,1.99); 
\draw [thick, dashed] (7.16,0.75)--(7.16,3.75); 
\draw [thick,dashed,rounded corners=3pt] (1.5,0.4)--(2.5,0.4)--(3.5,4.1)--(5.5,4.1)--(6.6,0.4)--(7.7,0.4); 
\draw [thick,dashed,rounded corners=3pt] (1.5,4.1)--(2.5,4.1)--(3.5,0.4)--(5.5,0.4)--(6.6,4.1)--(7.7,4.1); 
\node at (8.4,0) {$E_8$};
\node at (8.4,4.5) {$E_7$};
\node at (8,4) {$\ell_9$};
\node at (8,0.5) {$\ell_8$};
\node at (4,-1.2) {$\ell_2$};
\node at (1.4,-1.2) {$\ell_1$};
\node at (4,5.7) {$\ell_5$};
\node at (1.3,5.7) {$\ell_4$};
\node at (5.2,2.25) {$E_{12}$};
\node at (6.4,5.7) {$E_6$};
\node at (7.6,2.25) {$\ell_3$};
\node at (6.4,-1.2) {$E_1$};
\node at (1.4,1.4) {$E_9$};
\node at (1.3,3) {$E_{10}$};
\end{tikzpicture}
\vspace{-1em}
\end{subfigure} 
\end{figure}

\noindent
we consider a birational morphism from $\tilde{S}$ to $\mathbb{F}_{2}$ obtained by contracting $\ell_{2}, E_{12}, \ell_{4}, E_{9}, E_{10}$, $\ell_{9}$ and $E_{6}$ (resp. $\ell_{1}, E_{10}, E_{9}, \ell_{5}, E_{12}$, $\ell_{8}$ and $E_{1}$). Then the cylinder defined by the images of $E_{1}, E_{7}, E_{8}$ in $\mathbb
{F}_{2}$ implies  that there exists a cylinder $U_{1}$ (resp. $U_{2}) \cong \mathbb{A}^{1} \times \mathbb{A}^{1}_{\ast}$ in $S$ defined by $\overline{\ell_{2}}, \overline{\ell_{4}}, \overline{\ell_{9}}$ (resp. $\overline{\ell_{1}}, \overline{\ell_{5}}, \overline{\ell_{8}}$). Furthermore, for a general fiber $F$ of the $\mathbb{P}^{1}$-fibration, by contracting the negative curves on the lower left, we obtain the configuration on the lower right.
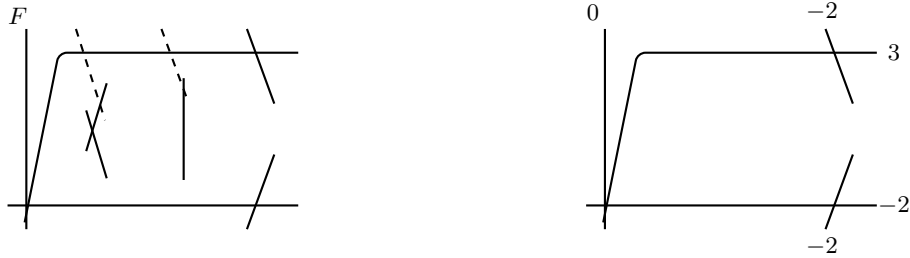
\begin{figure}[H]
  \centering
  \vspace{-1em}
\begin{subfigure}[t]{0.5\textwidth}\centering
  \begin{tikzpicture}[scale=0.45, every node/.style={font=\scriptsize}, baseline=(origin)]
\draw [thick] (-0.5,0)--(8,0); 
\draw [thick,rounded corners=3pt] (0,-0.5)--(1,4.5)--(8,4.5); 
\draw [thick] (6.5,-0.7)--(7.3,1.5); 
\draw [thick] (6.5,5.2)--(7.3,3); 
\draw [thick,dashed] (4,5.2)--(4.8,3); 
\draw [thick] (4.65,0.75)--(4.65,3.75); 
\draw [thick] (2.4,3.6)--(1.8,1.6); 
\draw [thick] (1.8,2.8)--(2.4,0.8); 
\draw [thick,dashed] (1.5,5.2)--(2.35,2.51); 
\draw [thick] (0.05,-0.7)--(0.05,5.2); 
\node at (-0.2,5.6){$F$};
\end{tikzpicture}
\end{subfigure}\hfill 
\begin{subfigure}[t]{0.5\textwidth}\centering
  \begin{tikzpicture}[scale=0.45, every node/.style={font=\scriptsize}, baseline=(origin)]
\draw [thick] (-0.5,0)--(8,0); 
\draw [thick,rounded corners=3pt] (0,-0.5)--(1,4.5)--(8,4.5); 
\draw [thick] (6.5,-0.7)--(7.3,1.5); 
\draw [thick] (6.5,5.2)--(7.3,3); 
\draw [thick] (0.05,-0.7)--(0.05,5.2); 
\node at (8.5,0) {$-2$};
\node at (8.5,4.5) {$3$};
\node at (6.4,5.7) {$-2$};
\node at (6.4,-1.2) {$-2$};
\node at (-0.3,5.7){$0$};
\end{tikzpicture}
\vspace{-2em}
\end{subfigure} 
\end{figure}

\noindent
Using the same argument as in the case of type $\mathsf{A_{4}+A_{3}}$, we obtain another cylinder $U_{3} \cong \mathbb{A}^{1} \times \mathbb{A}^{1}_{\ast}$ defined by $\overline{\ell}_{4}$, $\overline{\ell}_{5}$ and $\overline{F}\equiv \frac{2}{15}\overline{\ell}$. Therefore, $\Amp(S)$ is covered by the  union of the  relative interiors of $\Pol(U_{0}), \Pol(U_{1}), \Pol(U_{2})$ and $\Pol(U_{3})$.

\begin{figure}[H]
  \vspace{-1em}
\begin{subfigure}[t]{0.5\textwidth}\centering
\begin{tikzpicture}[scale=0.45, every node/.style={font=\scriptsize}, baseline=(origin)] %
  \begin{scope}
    \clip (-5,-5) rectangle (5,5);

    \clip (-5,187.5) -- (5,-187.5) -- (5,5) -- (-5,5) -- cycle;

    \clip (-5,46.875) -- (5,-46.875) -- (5,5) -- (-5,5) -- cycle;

    \fill[yellow!100] (-5,-5) rectangle (5,5);
  \end{scope}

  \begin{scope}
    \clip (-5,-5) rectangle (5,5);
    \fill[gray!20] (0,0) -- (-0.25,5) -- (0.25,5) -- cycle;
    \fill[gray!60] (0,0) -- (0.35,-5) -- (0.83,-5)  -- cycle;
  \end{scope}

  \draw[->] (-5.5,0) -- (5.5,0) node[right] {$a$};
  \draw[->] (0,-5.5) -- (0,5.5) node[above] {$b$};

  \draw[thick] (0,0) -- (-0.25,5);     
  \draw[thick] (0,0) -- (0.25,5);  
 \draw[thick] (0,0) -- (0.35,-5);  
  \draw[thick] (0,0) -- (0.83,-5);  
  \begin{scope}
    \clip (-5,-5) rectangle (5,5);
    \draw[dotted,very thick] (-5,187.5) -- (5,-187.5);    
    \draw[dotted,very thick] (-5,46.875) -- (5,-46.875); 
  \end{scope}

  \node[fill=white, inner sep=1pt, anchor=east] at (2.9,2.5) {$b=30a$};
  \node[fill=white, inner sep=1pt, anchor=east] at (-0.2,5.5) {$b=-30a$};
  \node[fill=white, inner sep=1pt, anchor=east] at (3.5,-2.5) {$b=-6a$};
  \node[fill=white, inner sep=1pt, anchor=east] at (0.7,-5.5) {$b=-10a$};
\end{tikzpicture}
\end{subfigure}\hfill
\begin{subfigure}[t]{0.5\textwidth}\centering
\begin{tikzpicture}[scale=0.45, every node/.style={font=\scriptsize}, baseline=(origin)] %
  \begin{scope}
    \clip (-5,-5) rectangle (5,5);

    \clip (-5,187.5) -- (5,-187.5) -- (5,5) -- (-5,5) -- cycle;

    \clip (-5,46.875) -- (5,-46.875) -- (5,5) -- (-5,5) -- cycle;

    \fill[yellow!100] (-5,-5) rectangle (5,5);
  \end{scope}

  \begin{scope}
    \clip (-5,-5) rectangle (5,5);
    \fill[gray!80] (0,0) -- (0.25,5) -- (5,5) -- (5,-5) -- (0.66,-5) -- cycle;
    
  \end{scope}

  \draw[->] (-5.5,0) -- (5.5,0) node[right] {$a$};
  \draw[->] (0,-5.5) -- (0,5.5) node[above] {$b$};

  \draw[thick] (0,0) -- (0.25,5);     
  \draw[thick] (0,0) -- (0.66,-5);  

  \begin{scope}
    \clip (-5,-5) rectangle (5,5);
    \draw[dotted,very thick] (-5,187.5) -- (5,-187.5);    
    \draw[dotted,very thick] (-5,46.875) -- (5,-46.875); 
  \end{scope}

  \node[fill=white, inner sep=1pt, anchor=east] at (2.9,2.5) {$b=30a$};
  \node[fill=white, inner sep=1pt, anchor=east] at (3.8,-2.5) {$b=-\frac{15}{2}a$};
\end{tikzpicture}
\end{subfigure}
\vspace{-0.5em}
\caption{$\Pol(U_{1}), \Pol(U_{2})$ and $\Pol(U_{3})$ for type $\mathsf{A_4+A_2+A_1}$}
\vspace{-1em}\hfill\qedhere
\end{figure}
\end{proof}

Therefore, by Theorems \ref{thm:basic} and \ref{thm:supple} we obtain the main result.
\begin{corollary}\label{cor:main}
For Du Val del Pezzo surfaces of degree $1$ with Picard rank $2$, the Conjecture \ref{conj} holds.
\end{corollary}

\newpage
\section{Appendix}\label{sec:appendix}
This section presents tables that summarize the results for each singularity type. The tigers in \cite{Cheltsov2016s} are listed in Table \ref{table:tiger} as follows.

\begin{center}
    \begin{longtable}{c|c|cc}
    \caption{Complements of Cylinders in Weak del Pezzo Surfaces}\label{table:tiger}\\
    \hline\hline
    No. & \begin{minipage}[m]{.15\linewidth}\setlength{\unitlength}{.10mm}
\begin{center}
\medskip

Singularity 

Type
\medskip
\end{center}
\end{minipage}  & 
    \begin{minipage}[m]{.5\linewidth}\setlength{\unitlength}{.25mm}
      \begin{center}
      Complement 
      
(divisor contracted if any)
      \end{center}
    \end{minipage} 
    & \begin{minipage}[m]{.2\linewidth}\setlength{\unitlength}{.10mm}
      \begin{center}
        Construction

      \end{center}
    \end{minipage} \\[0.6 ex]
    \hline
    \emph{1} & $\mathsf{A_4+A_2+A_1}$  & 
    \begin{minipage}[m]{.5\linewidth}\setlength{\unitlength}{.25mm}
      \begin{center}
      $\begin{array}{cc}
    &\frac{19}{30}E^{-3}_{1} + \frac{1}{6}E^{-1}_{5} + \frac{11}{30}E_6 + \frac{11}{15}E_7 + \frac{11}{10}E^{-6}_{8} \\ 
    &+ \frac{1}{10}E_9 + \frac{1}{5}E_{10} + \frac{3}{10}E^{-1}_{11} + \frac{1}{10}E_{12} + \frac{1}{5}E^{-1}_{13}\\
    &\left(\frac{19}{15}E^{-3}_{2} + \frac{9}{10}E_3 + \frac{7}{6}E_4 + \frac{49}{30}L^{-1}_{1} + \frac{41}{30}L_2 \right)
    \end{array}$ 
  \end{center}
  \end{minipage} & 
    \begin{minipage}[m]{.3\linewidth}\setlength{\unitlength}{.25mm}
\begin{center}
\begin{picture}(400,100)(0,-12)

\put(90,25){\line(1,0){30}}
\put(90,35){\line(1,0){30}}

\qbezier[25](20,12)(55,7)(75,-2)

\qbezier[25](55,-2)(85,7)(120,12)
\qbezier[20](30,5)(45,13)(63,20) \qbezier[25](55,73)(55,43)(55,13)
\qbezier[15](30,62)(45,62)(60,62)

 \put(7,12){\mbox{\tiny $L_1$}}
\put(40,-5){\mbox{\tiny $L_2$}}
 \put(3,47){\mbox{\tiny $\one$}}
\put(21,0){\mbox{\tiny $\two$}} \put(22,60){\mbox{\tiny $\three$}}
\put(52,73){\mbox{\tiny $\four$}} \put(28,27){\mbox{\tiny$\five$}}
  \put(141,56){\mbox{\tiny$\six$}}
  \put(121,73){\mbox{\tiny $\seven$}}
  \put(97,76){\mbox{\tiny$\eight$}}
  \put(159,33){\mbox{\tiny $\nine$}}
 \put(139,18){\mbox{\tiny $\ten$}}
 \put(79,32){\mbox{\tiny $\eleven$}}
 \put(139,8){\mbox{\tiny $\twelve$}}
 \put(79,22){\mbox{\tiny $\thirteen$}}

 \thicklines

\put(11,50){\line(2,1){30}}
\put(90,60){\line(2,1){30}}
\put(110,75){\line(2,-1){30}}
\put(100,75){\line(0,-1){80}}
\put(108,27){\line(2,-1){30}}
\put(108,37){\line(2,-1){30}}
\put(128,22){\line(2,1){30}}

\put(89,35){\circle*{3}}\put(89,25){\circle*{3}}

\dottedline[\tiny $\sim$]{6.5}(40,30)(70,30)

\end{picture}
\end{center}
\end{minipage}
  \\[5ex]
    \hline

     \emph{2} & $\mathsf{A_4+A_3}$  &  \begin{minipage}[m]{.5\linewidth}\setlength{\unitlength}{.25mm}
      $\begin{array}{cc}
    &\frac{19}{30}E^{-3}_{1} + \frac{1}{6}E^{-1}_{5} + \frac{11}{30}E_6 + \frac{11}{15}E_7 + \frac{11}{10}E^{-6}_{8}\\
    &+ \frac{1}{10}E^{-1}_{9} + \frac{1}{10}E_{10} + \frac{1}{5}E_{11} + \frac{3}{10}E_{12} + \frac{2}{5}E^{-1}_{13}\\
    &\left(\frac{19}{15}E^{-3}_{2} + \frac{9}{10}E_3 + \frac{7}{6}E_4 + \frac{49}{30}L^{-1}_{1} + \frac{41}{30}L_2 \right)
    \end{array}$
  \end{minipage} &
    \begin{minipage}[m]{.30\linewidth}\setlength{\unitlength}{.25mm}
\begin{center}
\begin{picture}(400,100)(0,-12)

\put(90,25){\line(1,0){30}}
\put(90,45){\line(1,0){30}}

\qbezier[25](20,12)(55,7)(75,-2)

\qbezier[25](55,-2)(85,7)(120,12)
\qbezier[20](30,5)(45,13)(63,20) \qbezier[25](55,73)(55,43)(55,13)
\qbezier[15](30,62)(45,62)(60,62)

 \put(7,12){\mbox{\tiny $L_1$}}
\put(40,-5){\mbox{\tiny $L_2$}}
 \put(3,47){\mbox{\tiny $\one$}}
\put(21,0){\mbox{\tiny $\two$}} \put(22,60){\mbox{\tiny $\three$}}
\put(52,73){\mbox{\tiny $\four$}} \put(28,27){\mbox{\tiny$\five$}}
  \put(141,56){\mbox{\tiny$\six$}}
  \put(121,73){\mbox{\tiny $\seven$}}
  \put(97,76){\mbox{\tiny$\eight$}}
  \put(121,42){\mbox{\tiny $\nine$}}
 \put(179,8){\mbox{\tiny $\ten$}}
 \put(159,25){\mbox{\tiny $\eleven$}}
 \put(139,8){\mbox{\tiny $\twelve$}}
 \put(121,22){\mbox{\tiny $\thirteen$}}

 \thicklines
\put(11,50){\line(2,1){30}}
\put(90,60){\line(2,1){30}}
\put(110,75){\line(2,-1){30}}
\put(100,75){\line(0,-1){80}}
\put(108,27){\line(2,-1){30}}
\put(128,12){\line(2,1){30}}
\put(148,27){\line(2,-1){30}}

\put(89,45){\circle*{3}}\put(89,25){\circle*{3}}

\dottedline[\tiny $\sim$]{6.5}(40,30)(70,30)

\end{picture}
\end{center}

\end{minipage}\\[4ex]
      \hline

     \emph{3} & $\mathsf{A_5+2A_1}$  & \begin{minipage}[m]{.5\linewidth}\setlength{\unitlength}{.25mm}
      \begin{center}
      $\begin{array}{cc}
    &\frac{2}{5}E_1 + \frac{4}{5}E_2 + \frac{6}{5}E^{-3}_3 + \frac{1}{5}E_4 + \frac{2}{5}E^{-1}_{5}\\ 
    &+ \frac{3}{5}E_6 + \frac{6}{5}E^{-4}_{7} + \frac{1}{5}E_8 + \frac{2}{5}E^{-1}_{9} + \frac{1}{5}E^{-1}_{10} \\
    &\left( \frac{7}{5}L_1 + \frac{8}{5}L^{-1}_2 \right)
    \end{array}$
  \end{center}
  \end{minipage}  &
    \begin{minipage}[m]{.3\linewidth}\setlength{\unitlength}{.25mm}
\begin{center}
\begin{picture}(400,75)(-10,-12)

\put(20,25){\line(1,0){30}}

\put(80,25){\line(1,0){30}}
\put(80,15){\line(1,0){30}}

\qbezier[20](30,10)(55,7)(75,-2)
\put(114,3){\mbox{\tiny $L_2$}}
\qbezier[20](55,-2)(75,7)(110,8)
\put(16,4){\mbox{\tiny $L_1$}}

\put(-9,24){\mbox{\tiny $\one$}} \put(11,45){\mbox{\tiny $\two$}}

\put(38,47){\mbox{\tiny $\three$}}\put(-4,7){\mbox{\tiny $\four$}}
\put(12,22){\mbox{\tiny $\five$}}

 \put(112,45){\mbox{\tiny
$\six$}} \put(88,47){\mbox{\tiny $\seven$}}
\put(127,37){\mbox{\tiny $\eight$}} \put(68,22){\mbox{\tiny
$\nine$}} \put(112,12){\mbox{\tiny $\ten$}}

 \thicklines
\put(4,13){\line(2,1){30}}
\put(0,30){\line(2,1){30}}
\put(20,45){\line(2,-1){30}}
\put(80,30){\line(2,1){30}}
\put(96,22){\line(2,1){30}}
\put(40,45){\line(0,-1){50}}
\put(90,45){\line(0,-1){50}}

\put(51,25){\circle*{3}}\put(79,25){\circle*{3}}\put(79,15){\circle*{3}}

\end{picture}
\end{center}

\end{minipage}\\[3ex]
      \hline

     \emph{4} & $\mathsf{A_5+A_2}$  & \begin{minipage}[m]{.5\linewidth}\setlength{\unitlength}{.25mm}
      \begin{center}
      $\begin{array}{cc}
    &\frac{2}{5}E_1 + \frac{4}{5}E_2 + \frac{6}{5}E^{-3}_3 + \frac{1}{5}E^{-1}_{4} + \frac{1}{5}E^{-1}_{5}\\ 
    &+ \frac{3}{5}E_6 + \frac{6}{5}E^{-4}_{7} + \frac{1}{5}E_8 + \frac{2}{5}E_9 + \frac{3}{5}E^{-1}_{10} \\
    &\left(\frac{7}{5}L_1 + \frac{8}{5}L^{-1}_2 \right)
    \end{array}$
  \end{center}
  \end{minipage}   &
    \begin{minipage}[m]{.30\linewidth}\setlength{\unitlength}{.25mm}
\begin{center}
\begin{picture}(400,75)(-10,-12)

\put(30,25){\line(1,0){30}}
\put(30,20){\line(1,0){30}}
\put(80,25){\line(1,0){30}}

\qbezier[20](30,10)(55,7)(75,-2)
\put(114,3){\mbox{\tiny $L_2$}}
\qbezier[20](55,-2)(75,7)(110,8)
\put(16,4){\mbox{\tiny $L_1$}}

\put(-9,24){\mbox{\tiny $\one$}} \put(11,45){\mbox{\tiny $\two$}}

\put(38,47){\mbox{\tiny $\three$}}\put(22,24){\mbox{\tiny
$\four$}} \put(22,15){\mbox{\tiny $\five$}}

 \put(112,45){\mbox{\tiny
$\six$}} \put(88,47){\mbox{\tiny $\seven$}}
\put(127,37){\mbox{\tiny $\nine$}} \put(68,22){\mbox{\tiny
$\ten$}} \put(148,20){\mbox{\tiny $\eight$}}

 \thicklines
\put(0,30){\line(2,1){30}}
\put(20,45){\line(2,-1){30}}
\put(80,30){\line(2,1){30}}
\put(96,22){\line(2,1){30}}
\put(116,37){\line(2,-1){30}}
\put(40,45){\line(0,-1){50}}
\put(90,45){\line(0,-1){50}}

\put(61,25){\circle*{3}}\put(79,25){\circle*{3}}
\put(61,20){\circle*{3}}

\end{picture}
\end{center}

\end{minipage}\\[3ex]
      \hline

     \emph{5} & $\mathsf{A_6+A_1}$ & \begin{minipage}[m]{.5\linewidth}\setlength{\unitlength}{.25mm}
      \begin{center}
      $\begin{array}{cc}
    &\frac{5}{12}E^{-1}_{1} + \frac{5}{12}E_2 + \frac{5}{6}E_3 + \frac{5}{4}E_4 + \frac{1}{4}E^{-1}_{5}\\ 
    &+ \frac{7}{12}E_6 + \frac{7}{6}E_7 + \frac{1}{3}E^{-1}_{8} + \frac{1}{6}L + \frac{17}{12}Q
    \end{array}$
  \end{center}
  \end{minipage}  &
    \begin{minipage}[m]{.30\linewidth}\setlength{\unitlength}{.25mm}
\begin{center}
\begin{picture}(400,80)(-10,-5)
\put(84,56){\line(1,0){30}}
\put(111,8){\eight}\put(38,12){\line(1,0){70}}
\put(111,28){\five} \put(75,32){\line(1,0){30}}

 \thicklines
\put(76,53){\one}
\put(10,45){\two}
\put(70,45){\three}
\put(50,45){\four}
\put(-10,5){\six}
\put(10,25){\seven}
\put(106,60){\mbox{\tiny $Q$}}
\put(82,2){\mbox{\tiny $L$}}

\put(20,45){\line(2,-1){30}}
\put(40,30){\line(2,1){30}}
 \put(60,45){\line(2,-1){30}}
\put(0,10){\line(2,1){30}}
\put(20,25){\line(2,-1){30}}
\qbezier(103,63)(70,30)(10,12)
 \qbezier(80,2)(61,12)(82,25)

\put(107,32){\circle*{3}}\put(107,12){\circle*{3}}
\put(116,56){\circle*{3}}

\end{picture}

\end{center}
\end{minipage}\\[3ex]
      \hline

     \emph{6} & $\mathsf{A_7'}$ &  \begin{minipage}[m]{.5\linewidth}\setlength{\unitlength}{.25mm}
      \begin{center}
      $\begin{array}{cc}
    &\frac{1}{3}E_1 + \frac{2}{3}E_2 + E_3 + \frac{4}{3}E_4 + \frac{1}{3}E^{-1}_{5}\\ 
    &+ \frac{2}{3}E_6 + \frac{4}{3}E_7 + \frac{1}{3}E^{-1}_{8} + \frac{1}{3}L^{-1} + \frac{4}{3}Q
    \end{array}$
  \end{center}
  \end{minipage} &
    \begin{minipage}[m]{.30\linewidth}\setlength{\unitlength}{.25mm}
\begin{center}
\begin{picture}(400,70)(-10,-10)

\put(112,8){\eight}\put(38,12){\line(1,0){70}}
\put(110,28){\five} \put(75,32){\line(1,0){30}}
\put(65,22){\mbox{\tiny $L$}} \qbezier(28,12)(46,17)(62,25)

 \thicklines
\put(-10,25){\one}
\put(10,45){\two}\put(70,45){\three}
\put(50,45){\four} \put(-10,5){\six}\put(10,25){\seven}
\put(93,47){\mbox{\tiny $Q$}}

\put(0,30){\line(2,1){30}}
\put(20,45){\line(2,-1){30}}
\put(40,30){\line(2,1){30}}
\put(60,45){\line(2,-1){30}}
\put(0,10){\line(2,1){30}}
\put(20,25){\line(2,-1){30}}
\qbezier(90,50)(70,30)(10,12)

\put(107,32){\circle*{3}}\put(109,12){\circle*{3}}
\put(28,12){\circle*{3}}

\end{picture}

\end{center}
\end{minipage}\\[3ex]
      \hline

     \emph{7} & $\mathsf{A_7''}$  & \begin{minipage}[m]{.5\linewidth}\setlength{\unitlength}{.25mm}
      \begin{center}
      $\begin{array}{cc}
    &\frac{7}{12}E_1 + \frac{7}{6}E_2 + \frac{19}{12}E_3 + \frac{7}{12}E^{-1}_{4} + \frac{5}{12}E_5\\ 
    &+ \frac{5}{6}E_6 + \frac{5}{4}E_7 + \frac{1}{4}E^{-1}_{8} + \frac{1}{6}L^{-1} + \frac{17}{12}Q
    \end{array}$
  \end{center}
  \end{minipage}  &
    \begin{minipage}[m]{.30\linewidth}\setlength{\unitlength}{.25mm}
\begin{center}
\begin{picture}(400,70)(-10,-5)

\put(58,22){\line(1,0){30}}
\put(58,42){\line(1,0){30}}
\qbezier(45,55)(35,45)(25,35)
\put(50,46){\mbox{\tiny $Q$}}

\put(46,56){\mbox{\tiny $L$}}

\put(91,18){\eight}

\put(-10,25){\one} \put(10,45){\two} \put(70,45){\three}
\put(91,38){\four}

\put(-10,5){\five} \put(10,25){\six} \put(30,5){\seven}
 \thicklines
\put(0,30){\line(2,1){30}}
\put(20,45){\line(2,-1){30}}
\put(40,30){\line(2,1){30}}
\put(0,10){\line(2,1){30}}
\put(20,25){\line(2,-1){30}}
\put(40,10){\line(2,1){30}}
\qbezier(54,43)(54,33)(54,2)

\put(87,42){\circle*{3}}\put(87,22){\circle*{3}}
\put(25,35){\circle*{3}}

\end{picture}

\end{center}
\end{minipage}\\[3ex]
      \hline

    \emph{8} & $\mathsf{D_5+2A_1}$ &  
    \begin{minipage}[m]{.5\linewidth}\setlength{\unitlength}{.25mm}
      \begin{center}
      $\begin{array}{cc}
    &\frac{4}{5}E_1 + \frac{8}{5}E_2 + \frac{6}{5}E_3 + \frac{1}{5}E^{-1}_{4} + \frac{1}{5}E_5\\ 
    &+ \frac{2}{5}E^{-1}_{6} + \frac{1}{5}E_7 + \frac{2}{5}E^{-1}_{8} + \frac{3}{5}L + \frac{6}{5}Q
    \end{array}$
  \end{center}
  \end{minipage}  &
    \begin{minipage}[m]{.30\linewidth}\setlength{\unitlength}{.25mm}
\begin{center}
\begin{picture}(400,65)(-10,-2)
\put(18,9){\line(1,0){35}}
\put(18,19){\line(1,0){30}}
\put(58,42){\line(1,0){30}}

\put(80,25){\mbox{\tiny $L$}}

\put(48,55){\mbox{\tiny $Q$}}

\put(-9,25){\mbox{\tiny $\one$}} \put(71,45){\mbox{\tiny
$\three$}}

\put(7,16){\mbox{\tiny $\six$}}\put(70,21){\mbox{\tiny $\seven$}}

\put(31,10){\mbox{\tiny $\five$}}\put(7,5){\mbox{\tiny $\eight$}}

\put(11,40){\mbox{\tiny $\two$}} \put(92,38){\mbox{\tiny $\four$}}

\thicklines
\put(0,30){\line(2,1){30}}
\put(20,45){\line(2,-1){30}}
\put(40,30){\line(2,1){30}}
\qbezier(45,43)(60,38)(76,33)
\qbezier(45,52)(25,25)(25,2)
\put(38,6){\line(2,1){30}}
\put(38,16){\line(2,1){30}}

\put(89,42){\circle*{3}}\put(18,9){\circle*{3}}
\put(18,19){\circle*{3}}

\end{picture}

\end{center}
\end{minipage}\\[3ex]
      \hline

     \emph{9} & $\mathsf{D_5+A_2}$ & \begin{minipage}[m]{.5\linewidth}\setlength{\unitlength}{.25mm}
      \begin{center}
      $\begin{array}{cc}
    &\frac{4}{5}E_1 + \frac{8}{5}E_2 + \frac{6}{5}E_3 + \frac{1}{5}E^{-1}_{4} + \frac{1}{5}E^{-1}_{5}\\ 
    &+ \frac{1}{5}E_6 + \frac{2}{5}E_7 + \frac{3}{5}E^{-1}_{8} + \frac{3}{5}L + \frac{6}{5}Q
    \end{array}$
  \end{center}
  \end{minipage} &
    \begin{minipage}[m]{.30\linewidth}\setlength{\unitlength}{.25mm}
\begin{center}
\begin{picture}(400,70)(-10,-5)
\put(18,24){\line(1,0){30}}
\put(18,9){\line(1,0){35}}
\put(58,42){\line(1,0){30}}

\put(80,25){\mbox{\tiny $L$}}

\put(48,55){\mbox{\tiny $Q$}}

\put(-9,25){\mbox{\tiny $\one$}} \put(71,45){\mbox{\tiny
$\three$}}

\put(7,21){\mbox{\tiny $\five$}}\put(70,21){\mbox{\tiny
$\seven$}}

\put(91,3){\mbox{\tiny $\six$}}\put(7,5){\mbox{\tiny $\eight$}}

\put(11,40){\mbox{\tiny $\two$}} \put(92,38){\mbox{\tiny $\four$}}

\thicklines
\put(0,30){\line(2,1){30}}
\put(20,45){\line(2,-1){30}}
\put(40,30){\line(2,1){30}}
\qbezier(45,43)(60,38)(76,33)
\qbezier(45,52)(25,25)(25,2)
\put(38,6){\line(2,1){30}}
\put(58,21){\line(2,-1){30}}

\put(89,42){\circle*{3}}\put(18,9){\circle*{3}}
\put(18,24){\circle*{3}}

\end{picture}

\end{center}
\end{minipage}  \\[3ex]
      \hline

     \emph{10} & $\mathsf{D_6+A_1}$ & \begin{minipage}[m]{.5\linewidth}\setlength{\unitlength}{.25mm}
      \begin{center}
      $\begin{array}{cc}
    &2E_1 + \frac{8}{5}E_2 + \frac{6}{5}E_3 + \frac{1}{5}E^{-1}_{4} + \frac{1}{5}E_5 + \frac{2}{5}E^{-1}_{6} \\
    &+ \frac{1}{5}E^{-1}_{7} + \frac{1}{5}E^{-1}_{8} + \frac{6}{5}L_1 + \frac{6}{5}L_2 + \frac{3}{5}L_3
    \end{array}$
  \end{center}
  \end{minipage}  &
    \begin{minipage}[m]{.30\linewidth}\setlength{\unitlength}{.25mm}
\begin{center}
\begin{picture}(200,95)(-10,-15)

\put(87,32){\mbox{\tiny $L_3$}}

\put(15,63){\mbox{\tiny $L_2$}}

\put(0,50){\mbox{\tiny $L_1$}}

\put(-5,7){\line(1,0){30}}

\put(30,20){\line(1,0){30}}
\put(30,25){\line(1,0){30}}

\put(78,58){\line(1,0){30}}

\put(-9,25){\mbox{\tiny $\one$}} \put(50,41){\mbox{\tiny
$\three$}}

\put(46,-8){\mbox{\tiny $\five$}}\put(63,23){\mbox{\tiny
$\seven$}}

\put(-15,2){\mbox{\tiny $\six$}}\put(63,15){\mbox{\tiny $\eight$}}

\put(34,62){\mbox{\tiny $\two$}} \put(112,54){\mbox{\tiny
$\four$}}

\thicklines

\put(0,30){\line(2,1){60}}
\put(40,60){\line(2,-1){30}}
\put(60,45){\line(2,1){30}}
\put(13,10){\line(2,-1){30}}
\qbezier(65,60)(75,50)(85,40)
\qbezier(20,60)(40,30)(40,12)
\qbezier(7,45)(7,25)(7,0)

\put(109,58){\circle*{3}}\put(30,20){\circle*{3}}
\put(-5,7){\circle*{3}}\put(30,25){\circle*{3}}

\end{picture}
\end{center}
\end{minipage}\\[3ex]
      \hline

     \emph{11} & $\mathsf{D_7}$ & 
     \begin{minipage}[m]{.5\linewidth}\setlength{\unitlength}{.25mm}
      \begin{center}
      $\begin{array}{cc}
    &\frac{3}{4}E_1 + \frac{3}{2}E_2 + \frac{7}{4}E_3 + 2E_4 + \frac{9}{4}E_5\\ 
    &+ \frac{5}{2}E_6 + \frac{1}{4}E^{-1}_{7} + \frac{1}{4}E^{-1}_{8} + \frac{1}{2}L^{-1} + \frac{5}{4}Q
    \end{array}$ 
  \end{center}
  \end{minipage} &
    \begin{minipage}[m]{.30\linewidth}\setlength{\unitlength}{.25mm}
\begin{center}
\begin{picture}(400,60)(-10,5)

\put(101,23){\line(1,0){30}}

\put(117,33){\line(1,0){30}}

\put(41,48){\mbox{\tiny $L$}}

\qbezier(41,46)(35,36)(27,32)

\put(77,45){\mbox{\tiny $Q$}}

\put(-9,25){\mbox{\tiny $\one$}} \put(31,25){\mbox{\tiny
$\three$}}\put(71,25){\mbox{\tiny $\five$}}

\put(11,40){\mbox{\tiny $\two$}} \put(51,40){\mbox{\tiny
$\four$}}\put(91,43){\mbox{\tiny $\six$}}\put(150,29){\mbox{\tiny
$\seven$}} \put(93,20){\mbox{\tiny $\eight$}}

\thicklines
\put(0,30){\line(2,1){30}}
\put(20,45){\line(2,-1){30}}
\put(40,30){\line(2,1){30}}
\put(60,45){\line(2,-1){30}}
\put(80,30){\line(2,1){30}}
\put(100,45){\line(2,-1){30}}
\qbezier(85,42)(111,36)(111,10)

\put(146,33){\circle*{3}}\put(132,23){\circle*{3}}
\put(27,32){\circle*{3}}

\end{picture}

\end{center}
\end{minipage}\\[3ex]
      \hline

     \emph{12} & $\mathsf{E_6+A_1}$ &  \begin{minipage}[m]{.5\linewidth}\setlength{\unitlength}{.25mm}
      \begin{center}
      $\begin{array}{cc}
    &\frac{12}{7}E^{-3}_{2} + \frac{10}{7}E_3 + \frac{15}{7}E_4 + \frac{8}{7}E_5 + \frac{1}{7}E^{-1}_{6}\\
    & + \frac{1}{7}E^{-1}_{7} + \frac{1}{7}E_8 + \frac{2}{7}E^{-1}_{9} + \frac{1}{7}E^{-1}_{10}\\ 
    &+ \frac{8}{7}L^{-3}_{1} + \frac{5}{7}L_3 \left(\frac{8}{7}L^{-1}_{2} + 2E_1 \right)
    \end{array}$ 
  \end{center}
  \end{minipage}&
    \begin{minipage}[m]{.30\linewidth}\setlength{\unitlength}{.25mm}
\begin{center}
\begin{picture}(200,85)(-10,-15)

\put(127,-6){\mbox{\tiny $L_3$}}

\qbezier[25](30,60)(40,30)(40,12)
\put(37,2){\mbox{\tiny $L_2$}}
\qbezier[30](0,30)(30,45)(60,60)

\put(1,50){\mbox{\tiny $L_1$}}

\put(-5,17){\line(1,0){30}}
\put(-5,7){\line(1,0){30}}

\put(105,33){\line(1,0){30}}
\put(-9,25){\mbox{\tiny $\one$}} \put(76,-4){\mbox{\tiny
$\three$}}

\put(46,-9){\mbox{\tiny $\eight$}} \put(-15,4){\mbox{\tiny
$\nine$}}\put(65,27){\mbox{\tiny $\ten$}} \put(-15,14){\mbox{\tiny
$\seven$}}

\put(34,62){\mbox{\tiny $\two$}} \put(87,47){\mbox{\tiny $\four$}}
 \put(139,29){\mbox{\tiny $\six$}}\put(76,16){\mbox{\tiny
$\five$}}

\thicklines
\qbezier(7,45)(7,25)(7,-5)
\qbezier(103,20)(113,10)(123,0)
\put(13,10){\line(2,-1){30}}
\put(40,60){\line(2,-1){60}}
\put(92,45){\line(0,-1){55}}
\put(87,20){\line(2,1){30}}
\put(87,0){\line(2,1){30}}

\put(136,33){\circle*{3}}
\put(-4,7){\circle*{3}}\put(-4,17){\circle*{3}}

\dottedline[\tiny $\sim$]{6.5}(30,30)(60,30)

\end{picture}
\end{center}
\end{minipage}\\[3ex]
      \hline

     \emph{13} & $\mathsf{E_7}$ & \begin{minipage}[m]{.5\linewidth}\setlength{\unitlength}{.25mm}
      \begin{center}
      $\begin{array}{cc}
&\frac{5}{3}E_1 + \frac{10}{3}E_2 + \frac{8}{3}E_3 + 2E_4 + \frac{4}{3}E_5 \\
&+\frac{1}{3}E^{-1}_{6} + \frac{4}{3}E_7 + \frac{1}{3}E^{-1}_{8} + \frac{1}{3}Q^{-1} + \frac{7}{3}L\\
\end{array}$
\end{center}
\end{minipage}  & 
    \begin{minipage}[m]{.30\linewidth}\setlength{\unitlength}{.25mm}
    \begin{center}
\begin{picture}(400,70)(-10,-10)

\qbezier[25](40,50)(5,20)(50,0)
\put(52,-5){\mbox{\tiny $L$}}

\put(112,10){\mbox{\tiny $Q$}}
\qbezier(76,50)(96,35)(116,20)

\put(-9,25){\mbox{\tiny $\one$}} \put(31,25){\mbox{\tiny
$\three$}}\put(71,25){\mbox{\tiny $\five$}}

\put(11,40){\mbox{\tiny $\two$}} \put(51,40){\mbox{\tiny
$\four$}}\put(144,40){\mbox{\tiny $\six$}}\put(12,-5){\mbox{\tiny
$\seven$}}\put(72,8){\mbox{\tiny $\eight$}}

\put(38,12){\line(1,0){30}}
\put(95,43){\line(1,0){45}}

\thicklines
\put(0,30){\line(2,1){30}}
\put(20,45){\line(2,-1){30}}
\put(40,30){\line(2,1){30}}
\put(60,45){\line(2,-1){30}}
\put(80,30){\line(2,1){35}}
\put(20,0){\line(2,1){30}}
\qbezier(40,50)(5,20)(50,0)

\put(139,43){\circle*{3}}\put(76,51){\circle*{3}}
\put(69,12){\circle*{3}}

\end{picture}
\end{center}
\end{minipage}\\[3ex]
    \hline\hline
    \end{longtable}
    \end{center}

The list of all $(-1)$-curves in minimal resolutions of Du Val del Pezzo surfaces constructed from the curves in $\mathbb{P}^{2}$ is presented in Table \ref{table:line}.  

\begin{center}
    \begin{longtable}{c|c|cc}
    \caption{Lines on minimal resolutions of Du Val del Pezzo surfaces}\label{table:line}\\
    \hline\hline
    No. & \begin{minipage}[m]{.15\linewidth}\setlength{\unitlength}{.10mm}
\begin{center}
\medskip
Singularity

Type
\medskip
\end{center}
\end{minipage}  & 
    \begin{minipage}[m]{.6\linewidth}\setlength{\unitlength}{.25mm}
      \begin{center}
    Lines
      \end{center}
    \end{minipage} 
    & \begin{minipage}[m]{.1\linewidth}\setlength{\unitlength}{.10mm}
      \begin{center}
      \end{center}
    \end{minipage} \\[0.6 ex]
    \hline
    \emph{1} & $\mathsf{A_4+A_2+A_1}$  & 
    \begin{minipage}[m]{.65\linewidth}\setlength{\unitlength}{.25mm}
      \begin{center}
      \footnotesize{$\begin{array}{cc}
        e_{11},e_{13},\ell-e_1-e_6,\\
          5\ell-3e_1-2e_2-e_3-e_4-e_5-2e_6-2e_7-e_8-e_9,\\
          5\ell-3e_1-2e_2-e_3-e_4-e_5-2e_6-2e_7-e_8-e_{12},\\
          6\ell-4e_1-2e_2-e_3-e_4-e_5-2e_6-2e_7-2e_8-e_9-e_{12},\\
          6\ell-4e_1-2e_2-e_3-e_4-e_5-2e_6-2e_7-2e_8-e_9-e_{10},\\
          6\ell-4e_1-2e_2-e_3-e_4-e_5-2e_6-2e_7-2e_8-e_{12}-e_{13},\\
          10\ell-6e_1-4e_2-2e_3-2e_4-2e_5-4e_6-3e_7-3e_8-e_9-e_{10}-e_{11},\\
          10\ell-6e_1-4e_2-2e_3-2e_4-2e_5-4e_6-3e_7-3e_8-e_9-e_{10}-e_{12},\\
          10\ell-6e_1-4e_2-2e_3-2e_4-2e_5-4e_6-3e_7-3e_8-e_9-e_{12}-e_{13},\\
          15\ell-9e_1-6e_2-3e_3-3e_4-3e_5-5e_6-5e_7-5e_8-2e_9-e_{10}-e_{11}-e_{12},\\
          15\ell-9e_1-6e_2-3e_3-3e_4-3e_5-5e_6-5e_7-5e_8-2e_9-e_{10}-e_{12}-e_{13},\\
          15\ell-9e_1-6e_2-3e_3-3e_4-3e_5-5e_6-5e_7-5e_8-e_9-e_{10}-e_{11}-2e_{12},\\
		      15\ell-9e_1-6e_2-3e_3-3e_4-3e_5-5e_6-5e_7-5e_8-e_9-e_{10}-2e_{12}-e_{13}.\\
        \end{array}$ } 
  \end{center}
  \end{minipage} & 
\begin{minipage}[m]{.1\linewidth}\setlength{\unitlength}{.25mm}
\end{minipage}
  \\[16ex]
    \hline

     \emph{2} & $\mathsf{A_4+A_3}$  &  \begin{minipage}[m]{.65\linewidth}\setlength{\unitlength}{.25mm}
      \begin{center}
      \footnotesize{$\begin{array}{cc}
        e_9,e_{13},\ell-e_1-e_6,\\
          5\ell-3e_1-2e_2-e_3-e_4-e_5-2e_6-2e_7-e_8-e_9,\\
          5\ell-3e_1-2e_2-e_3-e_4-e_5-2e_6-2e_7-e_8-e_{10},\\
          6\ell-4e_1-2e_2-e_3-e_4-e_5-2e_6-2e_7-2e_8-e_9-e_{10},\\
          6\ell-4e_1-2e_2-e_3-e_4-e_5-2e_6-2e_7-2e_8-e_{10}-e_{11},\\
          10\ell-6e_1-4e_2-2e_3-2e_4-2e_5-4e_6-3e_7-3e_8-e_9-e_{10}-e_{11},\\
          10\ell-6e_1-4e_2-2e_3-2e_4-2e_5-4e_6-3e_7-3e_8-e_{10}-e_{11}-e_{12},\\
          15\ell-9e_1-6e_2-3e_3-3e_4-3e_5-5e_6-5e_7-5e_8-2e_{10}-e_{11}-e_{12}-e_{13},\\
          15\ell-9e_1-6e_2-3e_3-3e_4-3e_5-5e_6-5e_7-5e_8-2e_9-e_{10}-e_{11}-e_{12},\\
          15\ell-9e_1-6e_2-3e_3-3e_4-3e_5-5e_6-5e_7-5e_8-e_9-2e_{10}-e_{11}-e_{12}.\\
        \end{array}$ }
      \end{center}
  \end{minipage} &
\begin{minipage}[m]{.1\linewidth}\setlength{\unitlength}{.25mm}
\end{minipage}\\[12ex]
\hline

     \emph{3} & $\mathsf{A_5+2A_1}$  & \begin{minipage}[m]{.65\linewidth}\setlength{\unitlength}{.25mm}
      \begin{center}
      \footnotesize{$\begin{array}{cc}
        e_5,e_9,e_{10},\ell-e_1-e_6,\\
          2\ell-e_1-e_2-e_6-e_7-e_8,\\
          2\ell-e_1-e_2-e_6-e_7-e_{10},\\
          3\ell-e_1-e_2-e_3-e_4-2e_6-e_7-e_8,\\
          3\ell-e_1-e_2-e_3-e_4-2e_6-e_7-e_{10},\\
          4\ell-2e_1-e_2-e_3-e_4-2e_6-2e_7-e_8-e_9,\\
          4\ell-2e_1-e_2-e_3-e_4-2e_6-2e_7-e_8-e_{10},\\
          6\ell-2e_1-2e_2-2e_3-e_4-e_5-3e_6-3e_7-2e_8-e_9,\\
          6\ell-2e_1-2e_2-2e_3-e_4-e_5-3e_6-3e_7-2e_8-e_{10},\\
          6\ell-2e_1-2e_2-2e_3-e_4-e_5-3e_6-3e_7-e_8-2e_{10},\\
          6\ell-2e_1-2e_2-2e_3-2e_4-3e_6-3e_7-e_8-e_9-e_{10}.
        \end{array}$}
  \end{center}
  \end{minipage}  &
    \begin{minipage}[m]{.1\linewidth}\setlength{\unitlength}{.25mm}
\end{minipage}\\[14ex]
      \hline

     \emph{4} & $\mathsf{A_5+A_2}$  & \begin{minipage}[m]{.65\linewidth}\setlength{\unitlength}{.25mm}
      \begin{center}
      
    \footnotesize{$\begin{array}{cc}
        e_4,e_5,e_{10},\ell-e_1-e_6,\\
          2\ell-e_1-e_2-e_6-e_7-e_8,\\
          3\ell-e_1-e_2-e_3-e_4-2e_6-e_7-e_8,\\
          3\ell-e_1-e_2-e_3-e_5-2e_6-e_7-e_8,\\
          4\ell-2e_1-e_2-e_3-e_4-2e_6-2e_7-e_8-e_9,\\
          4\ell-2e_1-e_2-e_3-e_5-2e_6-2e_7-e_8-e_9,\\
          6\ell-2e_1-2e_2-2e_3-e_4-e_5-3e_6-3e_7-2e_8-e_9,\\
          6\ell-2e_1-2e_2-2e_3-2e_4-3e_6-3e_7-e_8-e_9-e_{10},\\
          6\ell-2e_1-2e_2-2e_3-2e_5-3e_6-3e_7-e_8-e_9-e_{10}.\\
        \end{array}$ }
  \end{center}
  \end{minipage}   &
    \begin{minipage}[m]{.1\linewidth}\setlength{\unitlength}{.25mm}

\end{minipage}\\[11ex]
      \hline

     \emph{5} & $\mathsf{A_6+A_1}$ & \begin{minipage}[m]{.65\linewidth}\setlength{\unitlength}{.25mm}
      \begin{center}
      \footnotesize{$ \begin{array}{cc}
        e_1,e_5,e_8,\ell-e_2-e_3,\ell-e_1-e_2,\\
        \ell-e_1-e_6,L-e_2-e_6,\\
        2\ell-e_1-e_2-e_3-e_4-e_5,\\
        2\ell-e_2-e_3-e_4-e_5-e_6,\\
        3\ell-2e_2-e_3-e_4-e_5-e_6-e_7-e_8.
        \end{array}$}
  \end{center}
  \end{minipage}  &
    \begin{minipage}[m]{.1\linewidth}\setlength{\unitlength}{.25mm}

\end{minipage}\\[6ex]
      \hline

     \emph{6} & $\mathsf{A_7'}$ &  \begin{minipage}[m]{.65\linewidth}\setlength{\unitlength}{.25mm}
      \begin{center}
      \footnotesize{$\begin{array}{cc}
       e_5,e_8, \ell-e_1-e_2,\ell-e_1-e_6,\ell-e_6-e_7,\\
		  2\ell-e_1-e_2-e_3-e_4-e_5,\\
		  2\ell-e_1-e_2-e_6-e_7-e_8,\\
		  4\ell-e_1-e_2-e_3-e_4-e_5-2e_6-2e_7-2e_8.
        \end{array}$}
  \end{center}
  \end{minipage} &
    \begin{minipage}[m]{.1\linewidth}\setlength{\unitlength}{.25mm}

\end{minipage}\\[4ex]
      \hline

     \emph{7} & $\mathsf{A_7''}$  & \begin{minipage}[m]{.65\linewidth}\setlength{\unitlength}{.25mm}
      \begin{center}
      \footnotesize{$\begin{array}{cc}
       e_4,e_8,\ell-e_1-e_2, \ell-e_5-e_6,\ell-e_1-e_5,\\
        2\ell-e_1-e_5-e_6-e_7-e_8,\\
        2\ell-e_1-e_2-e_3-e_4-e_5.
        \end{array}$ }
  \end{center}
  \end{minipage}  &
    \begin{minipage}[m]{.1\linewidth}\setlength{\unitlength}{.25mm}

\end{minipage}\\[3ex]
      \hline

    \emph{8} & $\mathsf{D_5+2A_1}$ &  
    \begin{minipage}[m]{.65\linewidth}\setlength{\unitlength}{.25mm}
      \begin{center}
     \footnotesize{$\begin{array}{cc}
       e_4,e_6,e_8,\ell-e_5-e_6,\ell-e_7-e_8,\\
          \ell-e_1-e_5,\ell-e_1-e_7,\ell-e_5-e_7,\\
          3\ell-e_1-e_2-e_3-e_4-2e_5-e_7-e_8,\\
          3\ell-e_1-e_2-e_3-e_4-2e_5-e_6-e_7,\\
          3\ell-e_1-e_2-e_3-e_4-e_5-e_6-2e_7,\\
          3\ell-e_1-e_2-e_3-e_4-e_5-2e_7-e_8.
        \end{array}$  }
  \end{center}
  \end{minipage}  &
    \begin{minipage}[m]{.1\linewidth}\setlength{\unitlength}{.25mm}

\end{minipage}\\[7ex]
      \hline

     \emph{9} & $\mathsf{D_5+A_2}$ & \begin{minipage}[m]{.65\linewidth}\setlength{\unitlength}{.25mm}
      \begin{center}
      \footnotesize{$\begin{array}{cc}
       e_4,e_5,e_8,\ell-e_1-e_5,\ell-e_1-e_6,\\
          \ell-e_5-e_6,\ell-e_6-e_7,\\
          3\ell-e_1-e_2-e_3-e_4-e_5-2e_6-e_7,\\
          3\ell-e_1-e_2-e_3-e_4-2e_6-e_7-e_8,\\
          3\ell-e_1-e_2-e_3-e_4-2e_5-e_6-e_7.
        \end{array}$ }
  \end{center}
  \end{minipage} &
    \begin{minipage}[m]{.1\linewidth}\setlength{\unitlength}{.25mm}

\end{minipage}  \\[6ex]
      \hline

     \emph{10} & $\mathsf{D_6+A_1}$ & \begin{minipage}[m]{.65\linewidth}\setlength{\unitlength}{.25mm}
      \begin{center}
      \footnotesize{$\begin{array}{cc}
         e_4,e_6,e_7,e_8,\ell-e_5-e_7,\ell-e_5-e_8,\\
		  3\ell-e_1-e_2-e_3-e_4-2e_5-e_7-e_8,\\
		  3\ell-e_1-e_2-e_3-e_4-e_5-e_6-2e_7,\\
		  3\ell-e_1-e_2-e_3-e_4-e_5-e_6-2e_8.
        \end{array}$}
  \end{center}
  \end{minipage}  &
    \begin{minipage}[m]{.1\linewidth}\setlength{\unitlength}{.25mm}
\end{minipage}\\[4ex]
      \hline

     \emph{11} & $\mathsf{D_7}$ & 
     \begin{minipage}[m]{.65\linewidth}\setlength{\unitlength}{.25mm}
      \begin{center}
      \footnotesize{$\begin{array}{cc}
        e_7,e_8,\ell-e_1-e_2,\ell-e_1-e_8,\\
		  3\ell-2e_1-e_2-e_3-e_4-e_5-e_6-e_7.
        \end{array}$ }
  \end{center}
  \end{minipage} &
    \begin{minipage}[m]{.1\linewidth}\setlength{\unitlength}{.25mm}

\end{minipage}\\[2ex]
      \hline

     \emph{12} & $\mathsf{E_6+A_1}$ &  \begin{minipage}[m]{.65\linewidth}\setlength{\unitlength}{.25mm}
      \begin{center}
      \footnotesize{$\begin{array}{cc}
       e_6,e_7,e_9,\ell-e_8-e_{10},\ell-e_7-e_{10},\\
          5\ell-2e_1-2e_2-e_3-e_4-e_5-e_6-2e_7-e_8-3e_{10},\\
          5\ell-2e_1-2e_2-e_3-e_4-e_5-e_6-e_7-2e_8-3e_{10},\\   
          5\ell-2e_1-2e_2-e_3-e_4-e_5-e_6-2e_8-e_9-3e_{10}.
        \end{array}$ }
  \end{center}
  \end{minipage}&
    \begin{minipage}[m]{.1\linewidth}\setlength{\unitlength}{.25mm}

\end{minipage}\\[4ex]
      \hline

     \emph{13} & $\mathsf{E_7}$ & \begin{minipage}[m]{.65\linewidth}\setlength{\unitlength}{.25mm}
      \begin{center}
      \footnotesize{ $\begin{array}{cc}
        e_6,e_8,\ell-e_7-e_8\\
		  2\ell-e_1-e_2-e_3-e_4-e_5\\
		  5\ell-2e_1-2e_2-2e_3-2e_4-2e_5-2e_6-e_7-e_8.
        \end{array}$}
\end{center}
\end{minipage}  & 
    \begin{minipage}[m]{.1\linewidth}\setlength{\unitlength}{.25mm}
    
\end{minipage}\\[3ex]
    \hline\hline
    \end{longtable}
    \end{center}

Finally, the ample polar cylindricity of Du Val del Pezzo surfaces of degree 1 with Picard rank 2 is summarized in Table \ref{table:cyl} below.

\begin{center}
    \begin{longtable}{c|c|c|c|c|c}
    \caption{Ample Polar Cylindricity of Du Val del Pezzo Surfaces}\label{table:cyl}\\
    \hline\hline
    No. & \begin{minipage}[m]{.2\linewidth}\setlength{\unitlength}{.10mm}
\begin{center}
Singularity Type
\end{center}
\end{minipage}  & 
    \begin{minipage}[m]{.15\linewidth}\setlength{\unitlength}{.25mm}
      \begin{center}
    $\camp(S)$
      \end{center}
    \end{minipage} 
    & No. & \begin{minipage}[m]{.2\linewidth}\setlength{\unitlength}{.10mm}
      \begin{center}
      Singularity Type
      \end{center}
    \end{minipage} 
    & \begin{minipage}[m]{.15\linewidth}\setlength{\unitlength}{.25mm}
      \begin{center}
    $\camp(S)$
      \end{center}
    \end{minipage} 
    \\
    \hline
    \emph{1} & $\mathsf{A_4+A_2+A_1}$  & 
    \begin{minipage}[m]{.15\linewidth}\setlength{\unitlength}{.25mm}
      \begin{center}
      $\amp(S)$
  \end{center}
  \end{minipage} & 
\emph{11}
 & \begin{minipage}[m]{.2\linewidth}\setlength{\unitlength}{.10mm}
      \begin{center}
      $\mathsf{D_7}$
      \end{center}
    \end{minipage} 
    & \begin{minipage}[m]{.2\linewidth}\setlength{\unitlength}{.25mm}
      \begin{center}
    $\amp(S)$
      \end{center}
    \end{minipage} 
  \\[0.5ex]
    \hline

     \emph{2} & $\mathsf{A_4+A_3}$  &  \begin{minipage}[m]{.15\linewidth}\setlength{\unitlength}{.25mm}
      \begin{center}
      $\amp(S)$
      \end{center}
  \end{minipage} &
\emph{12}
 & \begin{minipage}[m]{.2\linewidth}\setlength{\unitlength}{.10mm}
      \begin{center}
      $\mathsf{E_6+A_1}$
      \end{center}
    \end{minipage} 
    & \begin{minipage}[m]{.2\linewidth}\setlength{\unitlength}{.25mm}
      \begin{center}
    $\amp(S)$
      \end{center}
    \end{minipage} 
  \\[0.5ex]
    \hline

     \emph{3} & $\mathsf{A_5+2A_1}$  & \begin{minipage}[m]{.15\linewidth}\setlength{\unitlength}{.25mm}
      \begin{center}
      $\amp(S)$
  \end{center}
  \end{minipage}  &
\emph{13}
 & \begin{minipage}[m]{.2\linewidth}\setlength{\unitlength}{.10mm}
      \begin{center}
      $\mathsf{E_7}$
      \end{center}
    \end{minipage} 
    & \begin{minipage}[m]{.2\linewidth}\setlength{\unitlength}{.25mm}
      \begin{center}
    $\amp(S)$
      \end{center}
    \end{minipage} 
  \\[0.5ex]
    \hline

     \emph{4} & $\mathsf{A_5+A_2}$  & \begin{minipage}[m]{.15\linewidth}\setlength{\unitlength}{.25mm}
      \begin{center}
     $\amp(S)$
  \end{center}
  \end{minipage}    &
\emph{14}
 & \begin{minipage}[m]{.2\linewidth}\setlength{\unitlength}{.10mm}
      \begin{center}
      $\mathsf{3A_2+A_1}$
      \end{center}
    \end{minipage} 
    & \begin{minipage}[m]{.2\linewidth}\setlength{\unitlength}{.25mm}
      \begin{center}
    $-K \not\in \text{Amp}^{\text{cyl}}$
      \end{center}
    \end{minipage} 
  \\[0.5ex]
    \hline

     \emph{5} & $\mathsf{A_6+A_1}$ & \begin{minipage}[m]{.15\linewidth}\setlength{\unitlength}{.25mm}
      \begin{center}
      $\amp(S)$
  \end{center}
  \end{minipage}   &
\emph{15}
 & \begin{minipage}[m]{.2\linewidth}\setlength{\unitlength}{.10mm}
      \begin{center}
      $\mathsf{A_3+4A_1}$
      \end{center}
    \end{minipage} 
    & \begin{minipage}[m]{.2\linewidth}\setlength{\unitlength}{.25mm}
      \begin{center}
    $-K \not\in \text{Amp}^{\text{cyl}}$
      \end{center}
    \end{minipage} 
  \\[0.5ex]
    \hline

     \emph{6} & $\mathsf{A_7'}$ &  \begin{minipage}[m]{.15\linewidth}\setlength{\unitlength}{.25mm}
      \begin{center}
$\amp(S)$
  \end{center}
  \end{minipage}  &
\emph{16}
 & \begin{minipage}[m]{.2\linewidth}\setlength{\unitlength}{.10mm}
      \begin{center}
      $\mathsf{A_3+A_2+2A_1}$
      \end{center}
    \end{minipage} 
    & \begin{minipage}[m]{.2\linewidth}\setlength{\unitlength}{.25mm}
      \begin{center}
    $-K \not\in \text{Amp}^{\text{cyl}}$
      \end{center}
    \end{minipage} 
  \\[0.5ex]
    \hline

     \emph{7} & $\mathsf{A_7''}$  & \begin{minipage}[m]{.15\linewidth}\setlength{\unitlength}{.25mm}
      \begin{center}
      $\amp(S)$
  \end{center}
  \end{minipage}   &
\emph{17}
 & \begin{minipage}[m]{.2\linewidth}\setlength{\unitlength}{.10mm}
      \begin{center}
      $\mathsf{2A_3+A_1}$
      \end{center}
    \end{minipage} 
    & \begin{minipage}[m]{.2\linewidth}\setlength{\unitlength}{.25mm}
      \begin{center}
    $-K \not\in \text{Amp}^{\text{cyl}}$
      \end{center}
    \end{minipage} 
  \\[0.5ex]
    \hline

    \emph{8} & $\mathsf{D_5+2A_1}$ &  
    \begin{minipage}[m]{.15\linewidth}\setlength{\unitlength}{.25mm}
      \begin{center}
    $\amp(S)$
  \end{center}
  \end{minipage}  &
\emph{18}
 & \begin{minipage}[m]{.2\linewidth}\setlength{\unitlength}{.10mm}
      \begin{center}
      $\mathsf{D_4+3A_1}$
      \end{center}
    \end{minipage} 
    & \begin{minipage}[m]{.2\linewidth}\setlength{\unitlength}{.25mm}
      \begin{center}
    $-K \not\in \text{Amp}^{\text{cyl}}$
      \end{center}
    \end{minipage} 
  \\[0.5ex]
    \hline

     \emph{9} & $\mathsf{D_5+A_2}$ & \begin{minipage}[m]{.15\linewidth}\setlength{\unitlength}{.25mm}
      \begin{center}
      $\amp(S)$
  \end{center}
  \end{minipage}  &
\emph{19}
 & \begin{minipage}[m]{.2\linewidth}\setlength{\unitlength}{.10mm}
      \begin{center}
      $\mathsf{D_4+A_3}$
      \end{center}
    \end{minipage} 
    & \begin{minipage}[m]{.2\linewidth}\setlength{\unitlength}{.25mm}
      \begin{center}
    $-K \not\in \text{Amp}^{\text{cyl}}$
      \end{center}
    \end{minipage} 
  \\[0.5ex]
    \hline

     \emph{10} & $\mathsf{D_6+A_1}$ & \begin{minipage}[m]{.15\linewidth}\setlength{\unitlength}{.25mm}
      \begin{center}
     $\amp(S)$
  \end{center}
  \end{minipage} & & 
  \\[0.5ex]
    \hline\hline
    \end{longtable}
    \end{center}

\begin{ack}
The first and second authors were supported by the National Research Foundation of Korea (NRF-2021R1A6A1A10039823), the second author is also partially supported by Basic Science Research Program through the National Research Foundation of Korea (NRF) funded by the Ministry of Education (No. RS-2023-00237440) and by Samsung Science and Technology Foundation under Project Number SSTF-BA2302-03, and the third author is supported by JSPS KAKENHI Grant Number JP24K22823 and JP25K17222.  
\end{ack}

\end{document}